\documentstyle[12pt,amssymb]{amsart}
\newcommand{\BBB}[1]{{\Bbb{#1}}}
\newcommand{\MM}{{\cal M}_3(\theta)}
\newcommand{\CC}{{\cal C}_5(\theta)}
\newcommand{\PP}{{\cal P}_3^{cm}(\theta)}
\newcommand{\BB}{{\cal B}_5^{cm}(\theta)}
\newcommand{\QQ}{Q_{\theta}:z\mapsto e^{2 \pi i \theta}z+z^2}
\newcommand{\FT}{f_{\theta}}
\newcommand{\iso}{\stackrel{\simeq}{\longrightarrow}}

\renewcommand{\marginpar}[1]{}
\catcode`\@=12

\def\Empty{}
\newcommand\oplabel[1]{
  \def\OpArg{#1} \ifx \OpArg\Empty {} \else
  	\label{#1}
  \fi}
		
\long\def\realfig#1#2#3#4{
\begin{figure}[htp]
\centerline{\psfig{figure=#2,width=#4}}
\caption[#1]{#3}
\oplabel{#1}
\end{figure}}

\newcommand{\comm}[1]{}
    
\input{psfig}
\newtheorem{thm}{Theorem}[section]
\newtheorem{cor}[thm]{Corollary}
\newtheorem{lem}[thm]{Lemma}
\newtheorem{prop}[thm]{Proposition}

\newcommand{\thmref}[1]{Theorem~\ref{#1}}
\newcommand{\propref}[1]{Proposition~\ref{#1}}

\newcommand{\lemref}[1]{Lemma~\ref{#1}}
\newcommand{\corref}[1]{Corollary~\ref{#1}} 
\newcommand{\figref}[1]{Fig.~\ref{#1}}

\def\SBIMSMark#1#2#3{
 \font\SBF=cmss10 at 10 true pt
 \font\SBI=cmssi10 at 10 true pt
 \setbox0=\hbox{\SBF Stony Brook IMS Preprint \##1}
 \setbox2=\hbox to \wd0{\hfil \SBI #2}
 \setbox4=\hbox to \wd0{\hfil \SBI #3}
 \setbox6=\hbox to \wd0{\hss
             \vbox{\hsize=\wd0 \parskip=0pt \baselineskip=10 true pt
                   \copy0 \break%
                   \copy2 \break%
                   \copy4 \break}}
 \dimen0=\ht6   \advance\dimen0 by \vsize \advance\dimen0 by 8 true pt
 \dimen2=\hsize \advance\dimen2 by .25 true in
 \ht6=0pt \dp6=0pt
 \setbox8=\vbox to \dimen0{\vfill \hbox to \dimen2{\hss \copy6}}
 \ht8=0pt \dp8=0pt \wd8=0pt
 \copy8
}

\begin{document}
\SBIMSMark{1998/4}{received December 1997}{and in revised form May 1998} 
\title{On Dynamics of Cubic Siegel Polynomials} 
\author{Saeed Zakeri}
\address{Department of Mathematics, SUNY at Stony Brook, NY 11794}
\email{zakeri@math.sunysb.edu}
\maketitle

{\large
\tableofcontents
}

\section{Introduction}
\label{sec:introduction}

Let $f$ be a polynomial of degree $d\geq 2$ in the complex plane and consider the following statements:

({\bf A}$_d$) ``If $f$ has a fixed Siegel disk $\Delta$ of bounded type rotation number,
then $\partial \Delta$ is a quasicircle passing through some critical point of $f$.''

({\bf B}$_d$) ``If $f$ has a fixed Siegel disk $\Delta$ such that $\partial \Delta$ is a quasicircle passing through some critical point of $f$, then the rotation number of $\Delta$ is bounded type.''

Statement ({\bf A}$_2$) is a theorem of Douady, Ghys, Herman, and Shishikura,
({\bf B}$_d$) is open, even for $d=2$, and the main object of this work is to prove ({\bf A}$_3$):\\ \\
{\bf Theorem 14.7.} {\it Let $P$ be a cubic polynomial which has a fixed Siegel disk $\Delta$ of rotation number $\theta$. Let $\theta$ be of bounded type. Then the boundary of $\Delta$ is a quasicircle which contains one or both critical points of $P$.}\\ \\
Along the way, we prove several results about the dynamics of cubic Siegel polynomials. In fact, we study the one-dimensional slice $\PP$ in the cubic parameter space which consists of all cubics with a fixed Siegel disk of a given rotation number $\theta$. Many of the results apply to general $\theta$ of Brjuno type.  

Siegel disks are examples of quasiperiodic motion in holomorphic dynamical systems. Let $p$ be an {\it irrationally indifferent} fixed point of a rational map $f:\overline{\BBB C}\rightarrow \overline{\BBB C}$. This means that $f(p)=p$ and the {\it multiplier} $f'(p)=\lambda$ is of the form $e^{2 \pi i \theta}$, where the {\it rotation number} $0< \theta<1$ is irrational. $p$ is called {\it linearizable} if there exists a holomorphic change of coordinates near $p$ which conjugates $f$ to the rigid rotation $z\mapsto \lambda z$. The largest domain on which this linearization is possible is a simply-connected domain $\Delta$ which is called the {\it Siegel disk} of $f$ centered at $p$. In other words, there exists a conformal isomorphism $h:(\BBB D,0)\iso (\Delta,p)$ such that $h(\lambda z)=f(h(z))$ for all $z\in \BBB D$, and $\Delta$ is not contained in any larger domain with this property. While the Siegel disk $\Delta$ is a component of the Fatou set of $f$, the boundary of $\Delta$ is a subset of the Julia set. 

Every punctured Siegel disk $\Delta \smallsetminus \{ p\}$ is foliated by dynamically-defined real-analytic invariant curves. However, as we get close to $\partial \Delta$, these invariant curves usually become more and more wiggly, and in the limit, we lose the control over the distortion of them. So, a priori, we do not even know if the $\partial \Delta$ is a Jordan curve. This question is difficult to answer partially because of the delicate analytic issues which arise in the study of the boundary behavior of the (essentially unique) linearizing map $h$.
 
It was conjectured by Douady and Sullivan in the early $80$'s that the boundary of every Siegel disk for a rational map has to be a Jordan curve (see \cite{Douady1}). This has still remained an open problem, even for polynomials, even when the degree is $2$. Even worse, there are very few explicit examples of polynomials for which we can effectively verify the conjecture. For instance, it is easy to see that local-connectivity of the Julia set implies the boundary of a Siegel disk to be a Jordan curve, but except for one case in the quadratic family \cite{Petersen}, we do not know how to check local-connectivity of the Julia set of a rational map which has a Siegel disk (and even in that single case, the boundary being a Jordan curve is proved without any reference to local-connectivity!). On the other hand, there are examples of quadratic Siegel polynomials whose Julia sets are not locally-connected, yet the boundaries of the Siegel disks are quasicircles \cite{Herman3} or even smooth Jordan curves \cite{Perez3}.

It is known that in any counterexample to this conjecture, the boundary of the Siegel disk must either be very complicated (an indecomposable continuum) or very simple (a circle with infinitely many topologist's sine curves implanted on it) \cite{Rogers}.

Let $\theta=[a_1, a_2, \ldots, a_n, \ldots ]$ be the continued fraction
expansion of $\theta$ and let $p_n/q_n=\discretionary{}{}{}
[a_1, a_2, \ldots, a_n]$ be its $n$-th 
rational approximation, with every $a_i$ being a positive integer. 
According to the theorem of Brjuno-Yoccoz \cite{Yoccoz}, every holomorphic germ with an indifferent fixed point of multiplier $\lambda=e^{2 \pi i \theta}$ is linearizable if and only if $\theta$ satisfies
$$\sum_{n=1}^{\infty} \frac{ \log q_{n+1}}{q_n}< +\infty.$$ 
Such $\theta$, or the corresponding $\lambda$, is called of {\it Brjuno type}. It is not hard to show that this set has full measure on the unit circle. The set of irrational numbers of Brjuno type contains two important arithmetic subsets: (1) numbers of {\it Diophantine type}, the set of all $0<\theta<1$ for which there exist positive constants $C$ and $\nu$ such that $|\theta-p/q|>C/q^{\nu}$ for every rational number $0\leq p/q<1$; and (2) numbers of {\it bounded type}, the set of all $0<\theta<1$ for which $\sup_n a_n< +\infty$.  

By the above discussion, every rational map with an indifferent fixed point whose multiplier is of Brjuno type has a Siegel disk. However, whether the multiplier of every Siegel disk of a rational map has to be of Brjuno type is only known to be true for quadratic polynomials by a theorem of Yoccoz \cite{Yoccoz} (see also \cite{Perez1} for a partial generalization).

Another issue is the existence of critical points on the boundary of Siegel disks. This problem was first studied by Ghys \cite{Ghys} under the assumption that the boundary is a Jordan curve and the rotation number is Diophantine. Later Herman improved the result by showing that when the rotation number is Diophantine and the action on the boundary is injective, there must be a critical point on the boundary \cite{Herman1}. The idea is to extract a circle diffeomorphism from the action on the boundary when there is no critical point there, and then to use the condition on the rotation number to extend the linearization to a neighborhood of the boundary, which gives a contradiction. A very short proof of this theorem is now possible with the knowledge of the ``Siegel compacts'' as recently introduced by Perez-Marco \cite{Perez2}. In the case of quadratic polynomials, no critical point on the boundary of the Siegel disk automatically implies that the map acts injectively on this boundary. Hence one concludes that for $\theta$ of Diophantine type, the critical point of $\QQ$ is on the boundary of the Siegel disk centered at $0$.

Later Herman gave the first example of a $\theta$ of Brjuno type for which the boundary of the Siegel disk for $Q_{\theta}$ is disjoint from the entire orbit of the critical point \cite{Herman3}. 

The most significant example in which one can explicitly show that the boundary of a Siegel disk is a Jordan curve containing a critical point is the quadratic map $\QQ$, with $\theta$ being of bounded type. The idea is to consider the degree $3$ Blaschke product
$$\FT (z)=e^{2 \pi i t(\theta)} z^2 \left ( \frac{z-3}{1-3z} \right )$$
which has a double critical point at $1$ and $t(\theta)$ is chosen such that the rotation number of the restriction of $\FT$ to the unit circle is $\theta$. 
Using a theorem of Swiatek and Herman on quasisymmetric linearization of critical circle maps \cite{Swiatek},\cite{Herman2}, one can redefine $\FT$ on the unit disk to make it quasiconformally conjugate to the rigid rotation. After modifying the conformal structure on the unit disk and all its preimages, one applies the Measurable Riemann Mapping Theorem of Ahlfors and Bers to prove that the resulting topological picture is quasiconformally conjugate to a quadratic polynomial $Q$. But the image of the unit disk has to be a Siegel disk for $Q$ of rotation number $\theta$, and there is only one such quadratic, so $Q=Q_{\theta}$ up to an affine conjugacy, which proves ({\bf A}$_2$).

In any attempt to generalize this result to higher degrees, one must address
several questions. In fact, the main difficulty is not the surgery which can be performed in all degrees in a similar way, provided that one has the appropriate Blaschke products in hand. Instead, we have to face a different set of questions, such as the parametrization of the candidate Blaschke products by their critical points, the combinatorics of various ``drops'' of their Julia sets, the continuity of the surgery, and the injectivity of this operation. None of these questions arises in degree $2$, where the corresponding parameter spaces are single points.    
 
Let us briefly sketch the organization of this paper: In  Section \ref{sec:cubpar} we introduce a normal form for the critically marked cubic polynomials with a Siegel disk of a given rotation number $\theta$ centered at the origin. We show that the associated parameter space $\PP$ is homeomorphic to the punctured plane and has a symmetry induced by the inversion $c\mapsto 1/c$ through the unit circle. We then study elementary topological properties of the connectedness locus $\MM$ in $\PP$. 

In  Section \ref{sec:stab} we show that the Julia sets of cubics in $\PP$ move holomorphically away from the boundary of the connectedness locus $\MM$ where various bifurcations do occur. In particular, if some cubic has an indifferent cycle other than the center of the Siegel disk at the origin, it must belong to the boundary of $\MM$. Both facts resemble well-known properties of the Mandelbrot set. 

In  Section \ref{sec:intcomp} we study components of the interior of $\MM$. In the interior of $\MM$, we can observe two possibilities: (1) The free critical point approaches an attracting cycle. In this case the cubic is called {\it hyperbolic-like} and is renormalizable in the sense of the definition in  Section \ref{sec:rencub}. (2) The free critical point eventually maps into the Siegel disk centered at the origin. This is called a {\it capture} cubic. The hyperbolic-like and capture components are the only possibilities one expects. However, as in the case of the Mandelbrot set, it is not known if these cases in fact cover all possibilities. As a third possibily, a cubic in the interior of $\MM$ which is neither hyperbolic-like nor capture is called {\it queer}. The most significant property of these cubics is that their Julia sets support invariant line fields and in particular have positive Lebesgue measure. To show this, unlike the quadratic family where the holomorphic motion of the basin of infinity extends automatically to the whole plane by the $\lambda$-lemma, here we must do some extra work. In fact, when the rotation number $\theta$ is an arbitrary number of Brjuno type, we do not know if the boundary of the Siegel disks are Jordan curves. Hence it is difficult to extend the holomorphic motions to the grand orbits of the Siegel disks in order to construct deformations of a queer cubic. Following McMullen and Sullivan, we overcome this difficulty by an application of the so-called ``harmonic $\lambda$-lemma'' of Bers and Royden. This section ends with a static version of an extension lemma, which turns out to be useful later in the construction of the surgery map.

 Section \ref{sec:rencub} studies the class of renormalizable cubics in $\PP$. From every such cubic one can extract the quadratic Siegel polynomial $\QQ$ using the theory of polynomial-like maps. As an easy byproduct, using a theorem of Petersen \cite{Petersen}, we show that the Julia set of a hyperbolic-like cubic or a cubic with disconnected Julia set has measure zero when $\theta$ is bounded type. 

 Section \ref{sec:connect} supplies a proof of connectivity of $\MM$. The standard Douady-Hubbard map on the exterior component of $\BBB C^{\ast}\smallsetminus \MM$ turns out to be proper holomorphic of degree $3$, so in order to prove that this component is homeomorphic to a punctured disk one needs to show that the map only branches over infinity. Again, this difficulty does not appear in the proof of connectivity of the Mandelbrot set, where the similar map has degree $1$.

In  Section \ref{sec:qc} we characterize the quasiconformal conjugacy classes in $\PP$. The most important feature is the quasiconformal rigidity of the cubics on the boundary of $\MM$, which will be the crucial step in the proof of continuity of the surgery map in  Section \ref{sec:continuity}. The material here is standard, except possibly for the existence of centers for capture components, which follows from the fact that the condition of being capture is open.

 Section \ref{sec:blaschke} is the beginning of our study of an auxiliary family of degree $5$ critically marked Blaschke products. To define the parameter space here, we need to show that our Blaschke products can be parametrized by their critical points, a fact that is trivial in the polynomial case. We prove that such a ``critical parametrization'' is always possible. 

In  Section \ref{sec:blapar} we use this critical parametrization to define the parameter space $\BB$. Every $B\in \BB$ is a degree $5$ Blaschke product with superattracting fixed points at $0$ and $\infty$, a double critical point on the unit circle $\BBB T$, and a pair $\{ c, 1/\overline{c} \}$ of symmetric critical points which may or may not be on $\BBB T$. The space of all such critically marked Blaschke products is homeomorphic to the punctured plane. For a $B\in \BB$, we study the structure of the preimages of the unit disk which are called the {\it drops} of $B$. 

Section \ref{sec:surgery} describes a surgery on Blaschke products in $\BB$, with $\theta$ being of bounded type, in order to obtain critically marked cubics in $\PP$. The theorem of Swiatek-Herman on quasisymmetric linearization of critical circle maps is the main tool. Unlike the quadratic case, here we must address a new question: Does the result of the surgery depend on various choices we make along the way? The answer turns out to be negative by an application of the Bers Sewing Lemma.

Section \ref{sec:newcon} defines the connectedness locus $\CC$ in $\BB$. We present a dynamical meaning for this locus by finding an alternative description for the Julia sets of elements in $\BB$. This description turns out to be useful in the study of injectivity of the surgery map.

In  Section \ref{sec:continuity} we show that the surgery map ${\cal S}:\BB \rightarrow \PP$ is continuous. The proof depends strongly on the fact that the parameter spaces have one complex dimension. One expects the similar map in higher dimensions to be discontinuous. 

 Section \ref{sec:renbla} introduces the notion of a renormalizable Blaschke product.
We show that from every such map one can extract the standard degree $3$ Blaschke product $\FT$ introduced by Douady, Ghys, Herman and Shishikura. This will be a very useful fact in  Section \ref{sec:inject}, mainly because of the simple observation that $\FT$ is quasiconformally rigid.

In  Section \ref{sec:inject} we prove that the surgery map is injective on the set of all Blaschke products which map to hyperbolic-like cubics or cubics with disconnected Julia sets. We actually prove a stronger result: Any two Blaschke products in the fiber over a noncapture cubic must be quasiconformally conjugate, with the conjugacy being conformal except on the Julia set. The fact that the surgery map is proper and restricts to a homeomorphism between the complementary components of $\CC$ and $\MM$ allows us to deduce that it is surjective. This proves \thmref{main}.

Finally, in  Section \ref{sec:twocrit} we study the set of all cubics in $\PP$ which have both critical points on the boundary of their Siegel disks. We prove that this is a Jordan curve in the boundary of $\MM$, which in some sense is parametrized by an ``angle'' between the two critical points.\\ \\           
{\bf Acknowledgements.} I am grateful to Jack Milnor for many inspiring discussions and his moral support. Among other things, he provided me with his computer programs to create pictures of various dynamically defined objects and showed me with great patience how to write them. Further thanks are due to Misha Lyubich and Dierk Schleicher for very useful conversations during the Spring and Fall semesters of 1997 at Stony Brook.
\vspace{0.17in}

\section{A Cubic Parameter Space}
\label{sec:cubpar}

We would like to parametrize the space of all cubic polynomials which have a fixed Siegel disk of multiplier $\lambda=e^{2 \pi i \theta}$ centered at the origin, where $0< \theta <1$ is an irrational number of Brjuno type. By the theorem of Brjuno-Yoccoz, every holomorphic germ $w\mapsto e^{2 \pi i \theta} w+O(w^2) $ with $\theta $ of Brjuno type is holomorphically linearizable near 0 \cite{Yoccoz}. Therefore, any such cubic polynomial has to be of the form 
$$w \mapsto \lambda w+a_2w^2+a_3w^3,$$
where $(a_2,a_3)\in \BBB C \times \BBB C^{\ast}$. We can mark the critical points of this polynomial by assuming that they 
are located at the points $c$ and $1$ with $c\neq 0$. In fact, one can conjugate the
above cubic by the linear map $w\mapsto z=\alpha w$, and the new cubic in the $z$-plane will have the form
$$z\mapsto \lambda z+\frac{a_2}{\alpha}z^2+\frac{a_3}{\alpha ^2}z^3.$$
It is easy to see that a critical point of this map is located at $1$
if we choose $\alpha$ to be any root of the equation $\lambda \alpha^2 +2 a_2 \alpha +3 a_3=0$. In this case, the other critical point $c$ will satisfy  
$$c=\frac{\lambda \alpha ^2}{3a_3}$$
so that the map gets the form 
\begin{equation}
\label{eqn:normform}
P_c:z\mapsto \lambda z \left ( 1-\frac{1}{2}(1+\frac{1}{c})z+\frac{1}{3c}z^2
\right )
\end{equation}
with $c\in \BBB C^{\ast}$. 

We denote the space of all critically marked cubic polynomials of the 
form (\ref{eqn:normform}) by $\PP$. In other words, $\PP\simeq \BBB C^{\ast}$
is parametrized by the invariant $c$. By an abuse of notation, we often identify the cubic $P_c$ with the parameter $c$. Note that $P_c$ and $P_{1/c}$ are affinely conjugate as maps, but certainly their critical points have different marking. Hence they will be regarded as distinct elements of $\PP$.  

In the topology of $\PP$, a sequence $P_n$ converges to some $P$ if there exist
$c_n,c\in \BBB C^{\ast}$, with $P_n=P_{c_n}$ and $P=P_c$, such that $c_n \rightarrow c$ as $n\rightarrow \infty$. In other words, the topology is given by uniform convergence of cubics on compact sets respecting the convergence of the marked critical points. \\ \\
{\bf Notation and Terminology.} Throughout this paper, the Siegel disk of the cubic $P_c$ centered at the origin is denoted by $\Delta_c$. When we do not want to emphasize the dependence on $c$, we denote the Siegel disk of a cubic $P$ by $\Delta_P$. By the {\it grand orbit} $GO(\Delta_P)$ we mean the set of all points in the plane which eventually map to the Siegel disk under the iteration of $P$. In other words,
$$GO(\Delta_P)=\bigcup_{k\geq 0} P^{-k}(\Delta_P).$$
{\bf Remark.} From classical Fatou-Julia theory (\cite{Milnor1}, Corollary 11.4), we know that every point on the boundary of the Siegel disk $\Delta_c$ must be in the closure of the orbit of either $c$ or $1$. According to Herman \cite{Herman1}, $P_c|_{\partial \Delta_c}$
has a dense orbit. It follows that the orbit of either $c$ or $1$ must accumulate on the entire boundary $\partial \Delta_c$.\\ 

The ``size'' of the Siegel disk $\Delta_c$ can be measured by the following invariant:\\ \\
{\bf Definition (Conformal Capacity).} Consider the Siegel disk $\Delta_c$ for $c \in {\BBB C}^{\ast}$ and the unique linearizing map $h_c:{\BBB D}(0,r_c) \iso \Delta_c$, with $h_c(0)=0$ and $h_c'(0)=1$. The radius $r_c>0$ of the domain of $h_c$ is called the {\it conformal capacity} of $\Delta_c$ and is denoted by $\kappa (\Delta_c)$. \\ \\
Alternatively, $\kappa (\Delta_c)$ can be described as the derivative $\varphi_c'(0)$ of the unique linearizing map $\varphi_c: {\BBB D} \iso \Delta_c$ normalized by $\varphi_c(0)=0$ and $\varphi_c'(0)>0$. Naturally, one is interested in the behavior of the function $c \mapsto \kappa (\Delta_c)$. The following lemma gives a basic result in this direction (compare \cite{Yoccoz}):

\begin{lem}
\label{upper}
The conformal capacity function $c \mapsto \kappa (\Delta_c)$ is upper semicontinuous.
\end{lem}

\begin{pf}
\comm{Consider the power series expansion of $h_c$ near $z=0$ for $c \in {\BBB C}^{\ast}$:
$$h_c(z)=z+a_2(c)z^2+a_3(c)z^3+ \cdots$$
Since $h_c(\lambda z)=P_c(h_c(z))$ for $z$ sufficiently close to zero and $P_c$ depends holomorphically in $c$, one can explicitly calculate the coefficients $a_j(c)$ to conclude that for every $j$, $c\mapsto a_j(c)$ is holomorphic. We claim that the conformal capacity $\kappa (\Delta_c)$ is equal to the radius of convergence of the above power series. If not, there exists an $r>\kappa (\Delta_c)$ such that $h_c$ can be continued analytically on the disk ${\BBB D}(0,r)$. The conjugacy relation $h_c(\lambda z)=P_c(h_c(z))$ persists over this larger disk. This implies that $P_c$ has bounded orbits on an open neighborhood of the closed disk $\overline{\Delta}_c$, which is a contradiction. We conclude that 
$$\kappa (\Delta_c)=\frac{1}{\limsup_{j \rightarrow \infty} |a_j(c)|^{1/j}},$$
or 
$$\log (\kappa (\Delta_c))= - \limsup_{j \rightarrow \infty} \frac{\log |a_j(c)|}{j}.$$
Since $\log|a_j(c)|$ is harmonic, ... }
Let $c_n \rightarrow c$ and $\kappa(\Delta_{c_n}) \geq r$. We would like to prove that $\kappa(\Delta_{c}) \geq r$ as well. The sequence of normalized univalent maps $h_{c_n}: {\BBB D}(0,\kappa(\Delta_{c_n})) \rightarrow {\BBB C}$ is normal on ${\BBB D}(0,r)$, so we may assume that a subsequence converges locally uniformly to a univalent function $h:{\BBB D}(0,r) \rightarrow {\BBB C}$. Since $h(\lambda z)=P_c(h(z))$ trivially, $h({\BBB D}(0,r))$ must be contained in the Siegel disk $\Delta_c$. Hence $\kappa(\Delta_{c}) \geq r$. 
\end{pf}

Since the conformal capacity $\kappa (\Delta_c)$ is upper semicontinuous by the above lemma, a priori it can jump to a {\it lower} value by a small perturbation. Intuitively, this means that the size of the Siegel disk $\Delta_c$ can become much smaller by a very small perturbation of the cubic $P_c$. Later we will see that for $\theta$ of bounded type, this cannot happen. In fact, in this case the closed Siegel disk $\overline{\Delta}_c$ is a quasidisk which moves continuously in the Hausdorff topology on compact subsets of the plane (see \thmref{move}). Therefore $\kappa (\Delta_c)$ is actually continuous as a function of $c$. On the other hand, for arbitrary $\theta$ of Brjuno type, I do not know if $c \mapsto \kappa (\Delta_c)$ is continuous. However, we have the following general theorem of Yoccoz \cite{Yoccoz}:

\begin{thm}
\label{size of disk}
Let $0< \theta <1$ be an irrational number of Brjuno type, and set $W(\theta)=\sum_{n=1}^{\infty} (\log q_{n+1})/q_n < \infty$. Let $S(\theta)$ be the space of all univalent functions $f:{\BBB D} \rightarrow {\BBB C}$ with $f(0)=0$ and $f'(0)=e^{2 \pi i \theta}$, with the maximal Siegel disk $\Delta_f \subset {\BBB D}$. Finally, define $\kappa(\theta)=\inf_{f \in S(\theta)} \kappa(\Delta_f)$. Then, there is a universal constant $C>0$ such that $|\log(\kappa(\theta))+W(\theta)|<C$.
\end{thm}

As an immediate corollary of the above theorem, we have:

\begin{cor}
\label{lower bound}
In the family $\{ P_c\} $ of cubic polynomials in $($\ref{eqn:normform}$)$, the conformal capacity function $c \mapsto \kappa (\Delta_c)$ is locally bounded away from $0$. 
\end{cor}

\noindent
{\bf Definition.} We define the {\it Cubic Connectedness Locus} $\MM$ as the set of all critically marked cubics $P\in \PP$ whose Julia sets
$J(P)$ are connected. It follows from classical Fatou-Julia theory (\cite{Milnor1}, Theorem 17.3) that $P\in \MM$ if and only if both critical points of $P$ have bounded orbits. We can formally set
$$\begin{array}{rl}
\MM & =\{ c \in {\BBB C}^{\ast}: \mbox{The Julia set $J(P_c)$
is connected} \} \\
                   & = \{ c \in {\BBB C}^{\ast}: \mbox{Both sequences $ \{P_c^{\circ k}(c)\}$  and $\{ P_c^{\circ k}(1)\}$ are bounded} \}. 
\end{array}$$
 
Since $P_c$ and $P_{1/c}$ are affinely conjugate as maps, neglecting the marking of the critical points, $\MM$ as a subset of the $c$-plane is invariant under the inversion $c\mapsto 1/c$ with respect to the unit circle. \figref{m3} shows the connectedness locus $\MM$ for the golden mean 
$\theta=(\sqrt 5 -1)/2=0.61803399...$ and \figref{m3zoom} shows the details of the same set near the unit circle.

\realfig{m3}{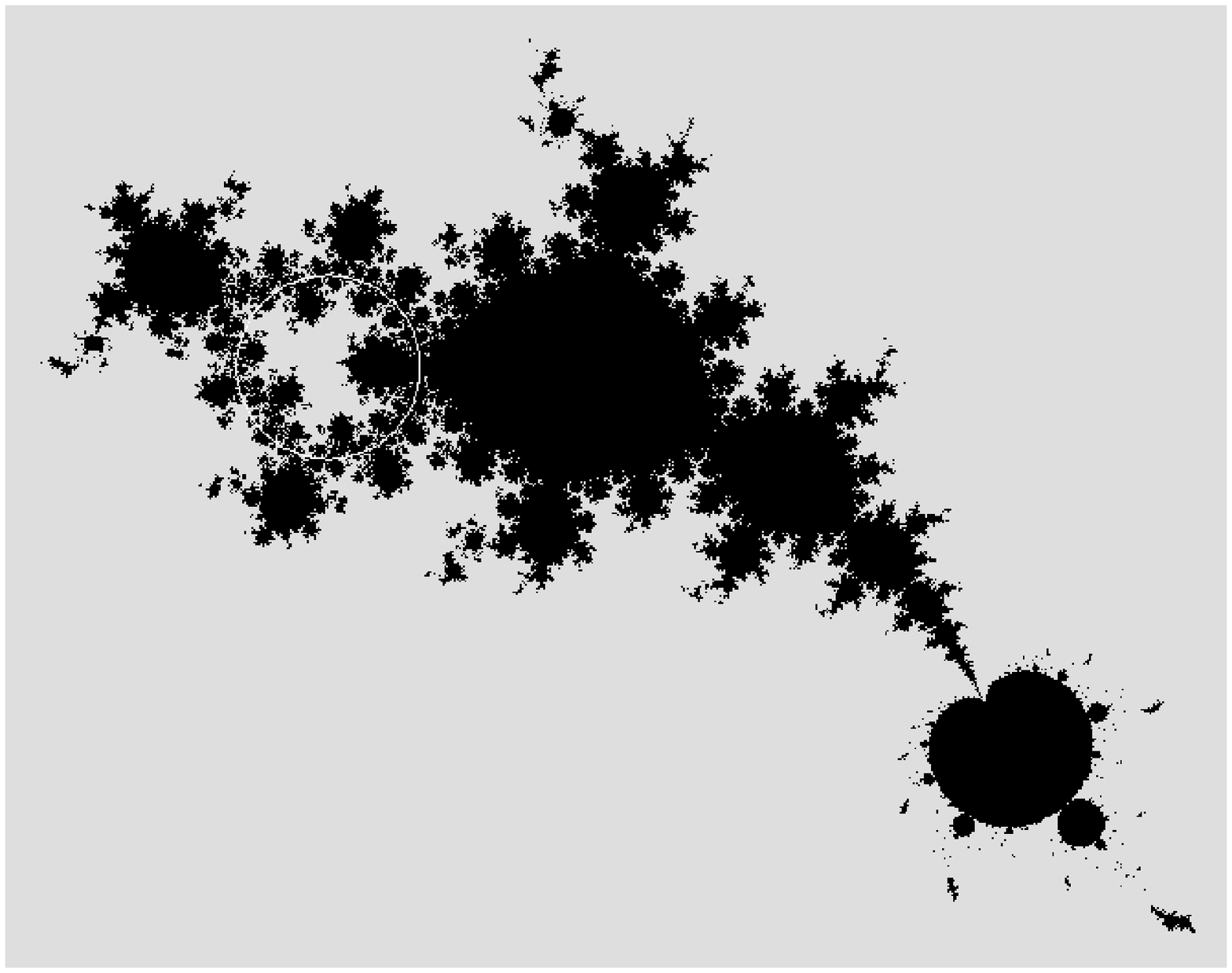}{{\sl The cubic connectedness locus $\MM$ for the golden mean $\theta=(\sqrt 5 -1)/2$. The circle in white is the unit circle centered at the origin.}}{9cm}

\realfig{m3zoom}{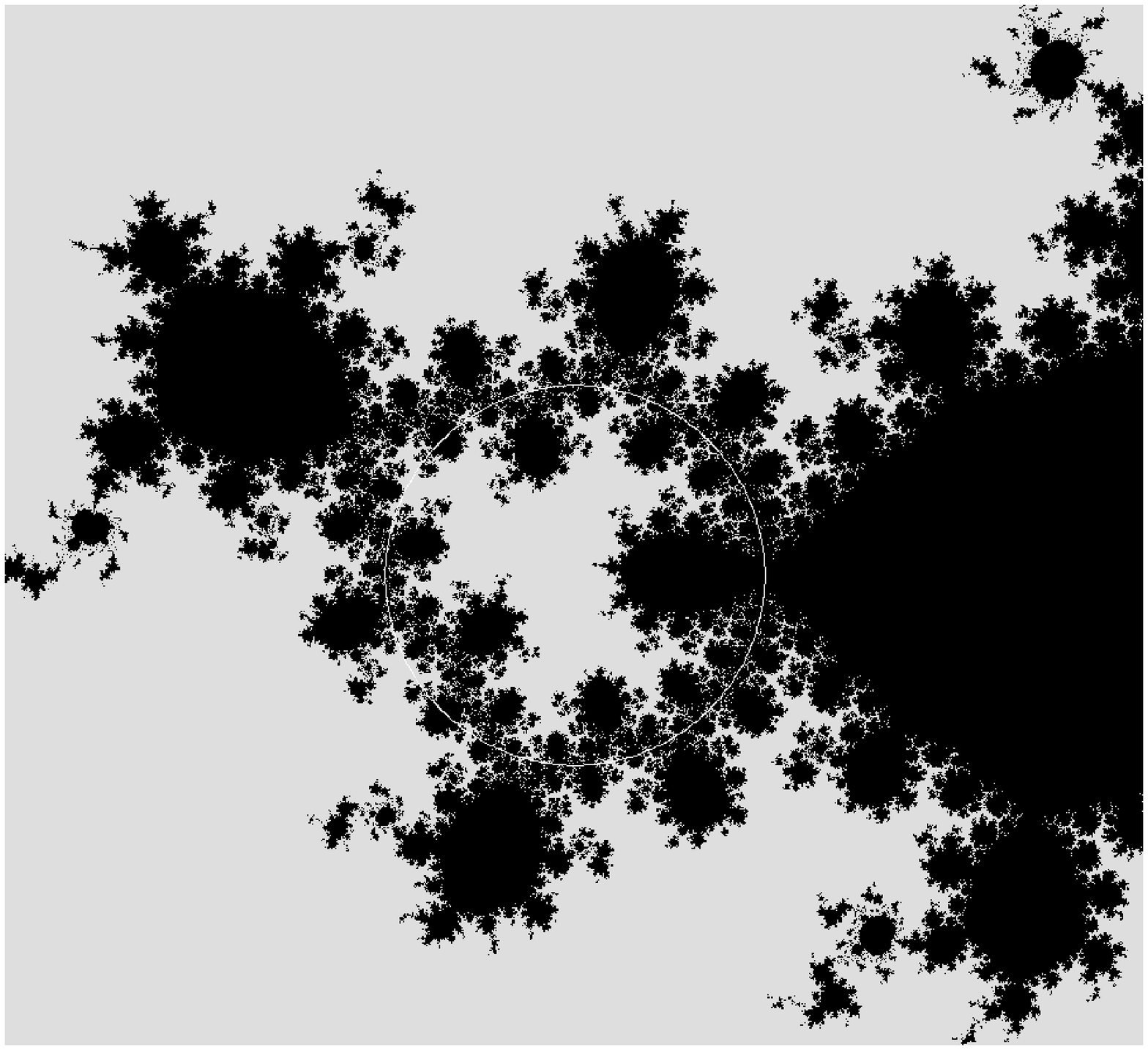}{{\sl Details of the connectedness locus $\MM$ near the unit circle centered at the origin.}}{13cm}
\begin{prop}
\label{m3com}

\noindent
\begin{enumerate}
\item[(a)]
$\MM$ is compact and contained in the open annulus $\BBB A(\frac{1}{30},30)$.
\item[(b)]
The complement $\BBB C^{\ast}\smallsetminus \MM $ has two connected components $\Omega_{ext}$ and $\Omega_{int}$ which are mapped to one another by the inversion $c\mapsto 1/c$. 
\end{enumerate}
\end{prop}

\begin{pf}
(a) $\MM$ is clearly closed. To see that it is bounded, we note that
$$\begin{array}{rl}
|P_c(z)| & = |\frac{1}{3c}z^2-\frac{1}{2}(1+1/c)z+1||z| \\ \\
         & \geq ( \frac{1}{|c|} | \frac{1}{3}z-\frac{1}{2}(c+1)||z|-1)|z|.
\end{array}$$
Let
\begin{equation}
\label{eqn:rcbound} 
m_c=(4.38)\max \{ |c|, 1 \}.
\end{equation}
If $|z|\geq m_c$, then 
$$\begin{array}{rl}
|P_c(z)| & \geq ( \frac{1}{|c|} ( \frac{1}{3}|z|-\frac{1}{4.38}|z|)|z|-1)|z|  \\ \\
         & \geq (0.46|z|-1) |z|\\ \\
         & \geq 1.0148\ |z|,
\end{array}$$
from which it follows that
\begin{equation}
\label{eqn:kcbound} 
K(P_c)\subset {\BBB D}(0,m_c),
\end{equation}
where $K(P_c)$ is the filled Julia set of $P_c$. Now let $|c|\geq 30$. Then 
$$|P_c(c)|=|\frac{1}{6}c-\frac{1}{2}||c|\geq (4.5)|c|>m_c,$$
which implies $P_c^{\circ k}(c)\rightarrow \infty$ as $k\rightarrow \infty$.
Therefore $\MM \subset {\BBB D}(0,30)$, hence by symmetry $\MM \subset \BBB A(\frac{1}{30},30)$.\\

(b) Let $\Omega_{ext}$ be the unbounded connected component of $\BBB C^{\ast} \smallsetminus \MM$. Since $\MM$ is invariant under the inversion $c\mapsto 1/c$, there exists a corresponding component $\Omega_{int}$ of the complement of $\MM$ containing a punctured neighborhood of the origin. By the proof of (a), we have 
$$\Omega_{ext}=\{ c\in \BBB C^{\ast}: P_c^{\circ k}(c)\rightarrow \infty\ \mbox{as}\ k \rightarrow \infty \},$$
and similarly
$$\Omega_{int}=\{ c\in \BBB C^{\ast}: P_c^{\circ k}(1)\rightarrow \infty\ \mbox{as}\ k \rightarrow \infty \}.$$
Suppose that there exists a bounded connected component $U$ of $\BBB C^{\ast} \smallsetminus \MM$ which is not $\Omega_{int}$. Then 
$$0<\sup_{c\in {\partial U}}|c|=R< + \infty.$$
If $c\in \partial U$, it follows from (\ref{eqn:kcbound}) that for each $k\geq 0$, $|P_c^{\circ k}(c)|$ and $|P_c^{\circ k}(1)|$ are not greater than $m_c$, and 
$$\sup_{c\in \partial U} m_c \leq (4.38)\max \{ R, 1 \} < + \infty.$$
Since $U\neq \Omega_{int}$, we have $\partial U \subset \partial \MM$ and both $P_c^{\circ k}(c)$ and $P_c^{\circ k}(1)$ are holomorphic in $U$. It follows from the Maximum Principle that the iterates $P_c^{\circ k}(c)$ and
$P_c^{\circ k}(1)$ are uniformly bounded throughout $U$, which is a contradiction. 
\end{pf}
{\bf Remarks.} 

(1) The bound $30$ in (a) is not sharp. Computer experiments show that 
it can actually be replaced by 11.266519.

(2) Later we will prove that $\Omega_{ext}$ (hence $\Omega_{int}$) is
homeomorphic to a punctured disk. This will show that $\MM$ is a connected set (see \thmref{m3con}).
\vspace{0.17in}

\section{Stability of Cubics} 
\label{sec:stab}

In this section we prove the following result, which is reminiscent of the
similar fact about the Mandelbrot set. For terminology and basic results
on holomorphic motions and $J$-stability, see \cite{McMullenbook} .

\begin{thm}[Boundary of $\MM$ is Unstable]
\label{unstable}
The complement $\BBB C^{\ast} \smallsetminus \partial \MM$ is the set of parameters for which the corresponding cubics are J-stable in $\PP$.
\end{thm}

\begin{pf} A polynomial $P_{c_0}\in \PP$ is $J$-stable if and only if both
sequences $\{ P_c^{\circ k}(c)\} $ and $\{ P_c^{\circ k}(1)\} $ are 
normal for $c$ in a neighborhood of $c_0$ (\cite{McMullenbook}, Theorem 4.2). If $c_0 \in \Omega_{ext}$, then $c_0$ escapes to infinity under $P_{c_0}$, while $1$ has bounded orbit. For $c$ close to $c_0$, the orbit of $c$ under $P_c$ will still converge to infinity while $1$
will have bounded orbit, with a bound given by $m_c$ in (\ref{eqn:rcbound}). It follows from the Montel's theorem that both sequences are normal throughout a neighborhood of $c_0$. Hence $c_0$ is $J$-stable. Similarly, every $P_{c_0}$
with $c_0 \in\Omega_{int}$ is $J$-stable. If $c_0$ belongs to the interior of $\MM$, then 
both $c_0$ and $1$ will have orbits contained in $\BBB D(0,m_{c_0})$ and
the same holds for all $c$ sufficiently close to $c_0$. Again by Montel, both sequences
$\{ P_c^{\circ k}(c)\} $ and $\{ P_c^{\circ k}(1)\} $ are normal in a neighborhood of $c_0$. Finally, if $c_0$ belongs to the boundary of $\MM$, then
a small perturbation will make either $c$ or $1$ escape to infinity. Hence
at least one of the sequences $\{ P_c^{\circ k}(c)\} $ or $\{ P_c^{\circ k}(1)\} $ fails to be normal in any neighborhood of $c_0$.
\end{pf}
 
\begin{thm}
\label{indiff}
Let $P_{c_0}\in \PP$ have an indifferent periodic orbit other than the fixed point at the origin. Then $c_0 \in \partial \MM$.
\end{thm}

\begin{pf}
Otherwise $c_0$ will be a $J$-stable parameter by the above
theorem. But any stable indifferent cycle has to be persistent (\cite{McMullenbook}, Theorem 4.2). This means that the indifferent cycle $z(c_0)\mapsto P_{c_0}(z(c_0))\mapsto \cdots \mapsto
P_{c_0}^{\circ k-1}(z(c_0))\mapsto z(c_0)$ can be continued analytically as a function of $c$ in a neighborhood of $c_0$ and
the multiplier function $c\mapsto (P_c^{\circ k})'(z(c))$ will be constant in
this neighborhood. But this cycle can be continued analytically to the whole $c$-plane except for a finite number of singular points by the implicit function theorem, and the multiplier has to remain constant during the continuation. It follows that for every parameter $c$, the cubic $P_c$ has an indifferent cycle other than $0$. This is clearly impossible, since for example when $c=3-6 \overline{\lambda}$, $P_c(c)=c$ is a superattracting fixed point, hence there cannot be any indifferent periodic point other than 0. 
\end{pf}

To prove the next corollary, we use the following lemma in \cite{Kiwi} which is a much sharpened version of an earlier result of Goldberg and Milnor (\cite{Goldberg-Milnor}, Theorem 3.3). This useful lemma will also be applied in Section \ref{sec:rencub} below to extract quadratic-like maps out of renormalizable cubics.

\begin{lem}[Separation Lemma]
\label{seplem}
Let $P$ be a polynomial with connected Julia
set. Then there exists a finite collection of closed preperiodic
external rays, separating the plane into disjoint open simply-connected sets $\{ U_j\}$, such that:
\begin{enumerate}
\item[$\bullet$]
Each $U_j$ contains at most one non-repelling periodic point or periodic Fatou component of $P$.
\item[$\bullet$]
If $z_1\mapsto \cdots \mapsto z_p\mapsto z_1$ is a non-repelling cycle meeting $U_{i_1}\mapsto \cdots \mapsto U_{i_p}\mapsto U_{i_1}$, 
then $\bigcup _{j=1}^p U_{i_j}$ contains the entire orbit of at least one critical point of $P$. 
\end{enumerate}
\end{lem}

\begin{cor}
\label{second ind}
If $P\in \PP$ has an indifferent periodic point other than
the fixed point at the origin, then a critical point of $P$, other than the one which accumulates on the boundary of $\Delta_P$, accumulates on the extra indifferent point (in case it is not linearizable)
or on the boundary of the extra Siegel disk (in case it is linearizable).\ \ $\Box$
\end{cor}

\vspace{0.17in}

\section{Components of the Interior of ${\cal M}_3(\theta)$}
\label{sec:intcomp} 
\noindent
{\bf Definition (Types of Components).}\ A component $U$ of the interior of $\MM$ is called {\it hyperbolic-like} if for every $c\in U$, the orbit of either $c$ or $1$
under $P_c$ converges to an attracting cycle. $U$ is called a {\it capture}
component if for every $c\in U$, either $c$ or $1$ eventually maps to the
Siegel disk $\Delta_c$. In case $U$ is neither hyperbolic-like nor capture, we call it a {\it queer} component.
We say that $P_c$ is hyperbolic-like, capture, or queer if the corresponding parameter $c$ belongs to such a component.\\

For example, there is a hyperbolic-like component in the form of the main cardioid of a large copy of the Mandelbrot set on the lower right corner of \figref{m3}. For every $c$ in this component, the orbit of the critical point $c$ of $P_c$ converges to an attracting fixed point. On the other hand, the large component
which is attached on the right side of the unit circle to $c=1$ is a capture,
consisting of all $c$ for which $P_c(c)$ belongs to $\Delta_c$. \figref{hypcub}-\figref{hypext} show examples of the filled Julia sets of cubics in $\PP$.

In the above definition, we tacitly assumed that hyperbolic-like
or capture cubics define components of the interior of $\MM$. The condition of being hyperbolic-like is clearly open. So to justify the definition in this case, we have to show that it is also closed in the interior of $\MM$.

\realfig{hypcub}{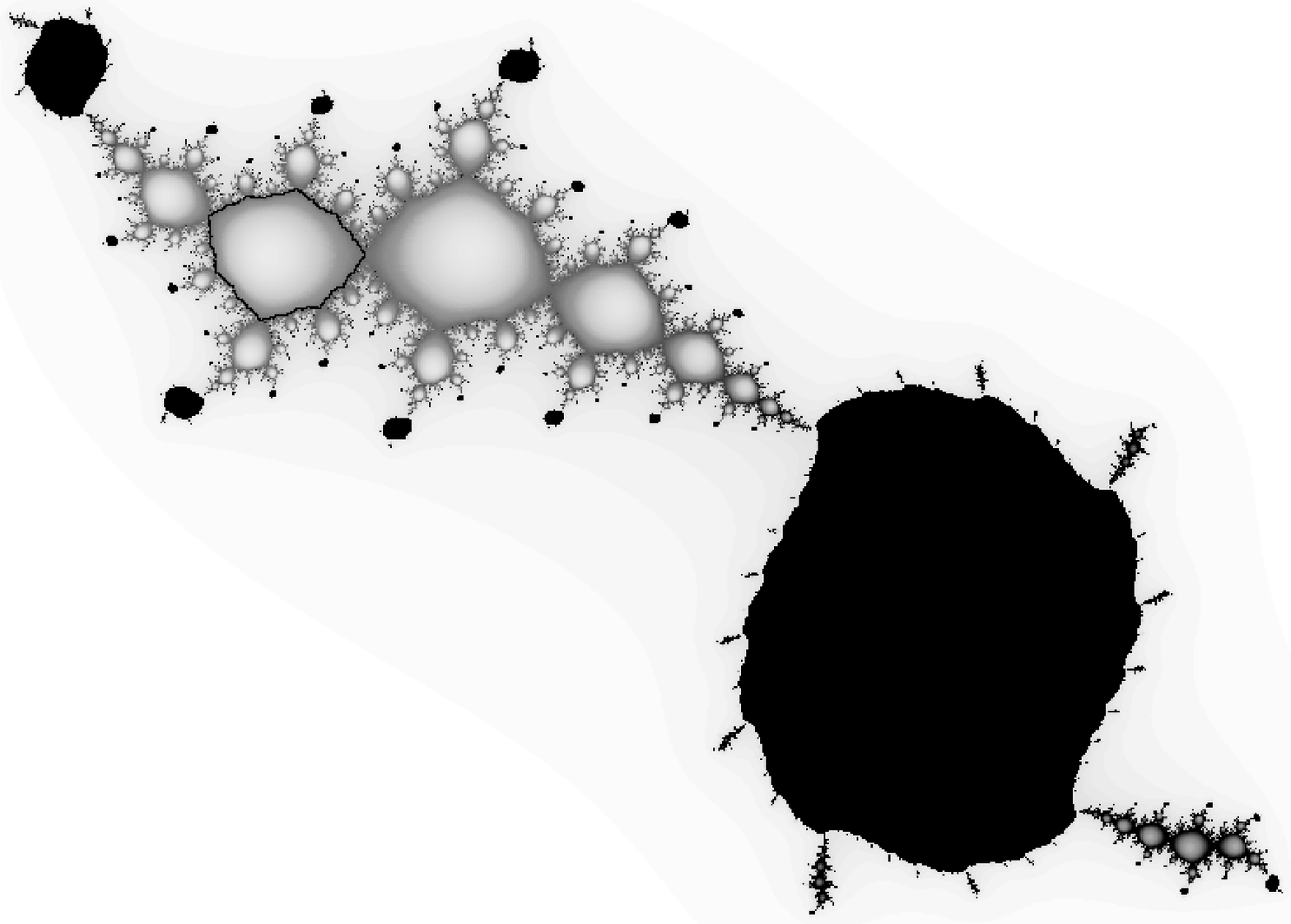}{{\sl The filled Julia set of a hyperbolic-like cubic in $\PP$ with $\theta=(\sqrt {5}-1)/2$. The large open topological disk on the right is the immediate basin of attraction of an attracting fixed point.}}{9cm}
\realfig{captcub}{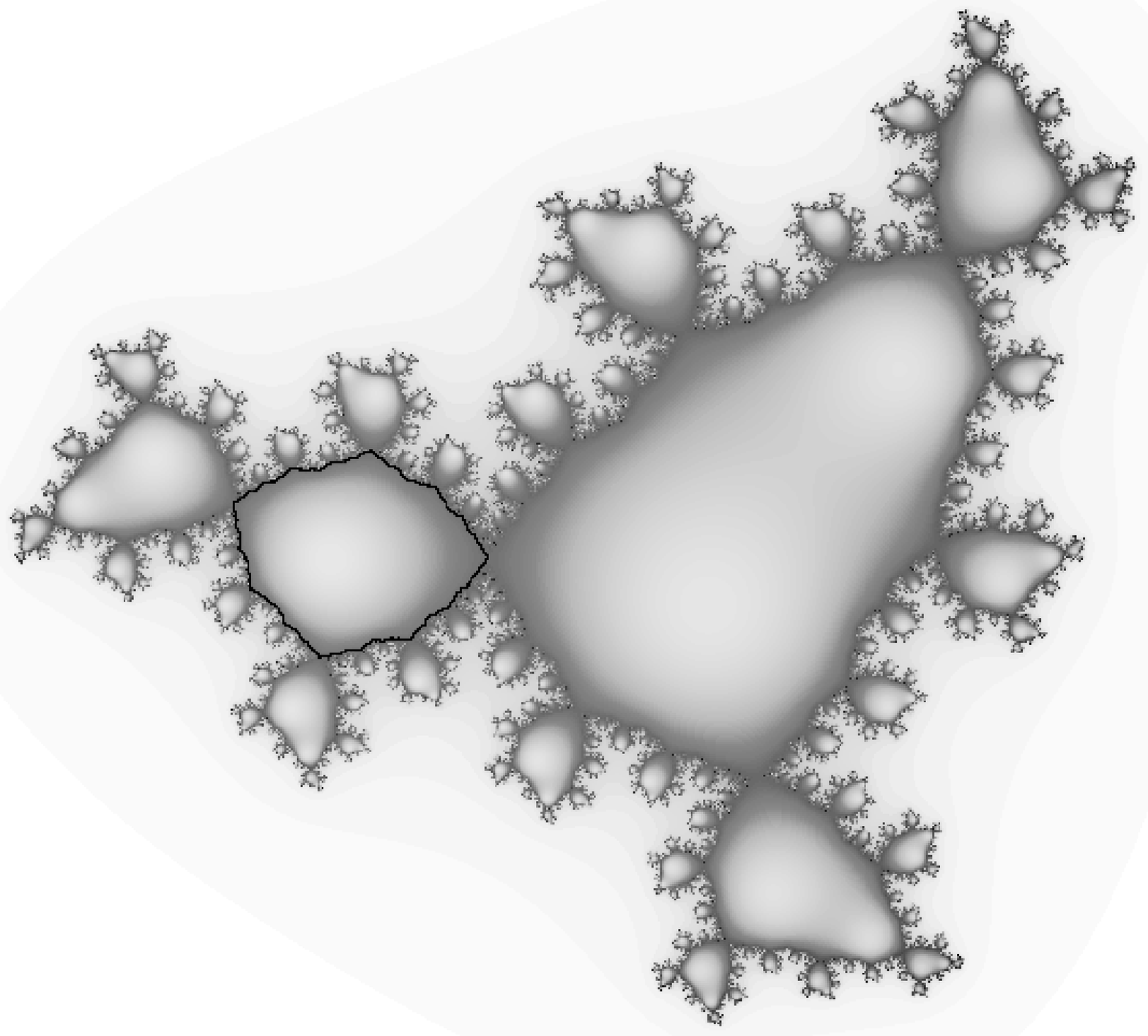}{{\sl The filled Julia set of a capture cubic in
$\PP$ with $\theta=(\sqrt {5}-1)/2$. Every bounded Fatou component eventually maps to the Siegel disk
centered at the origin. There is a critical point in the large preimage of the Siegel disk on the right. Hence this preimage maps to the Siegel disk by a 2-to-1 branched covering.}}{9cm}

\realfig{paracub}{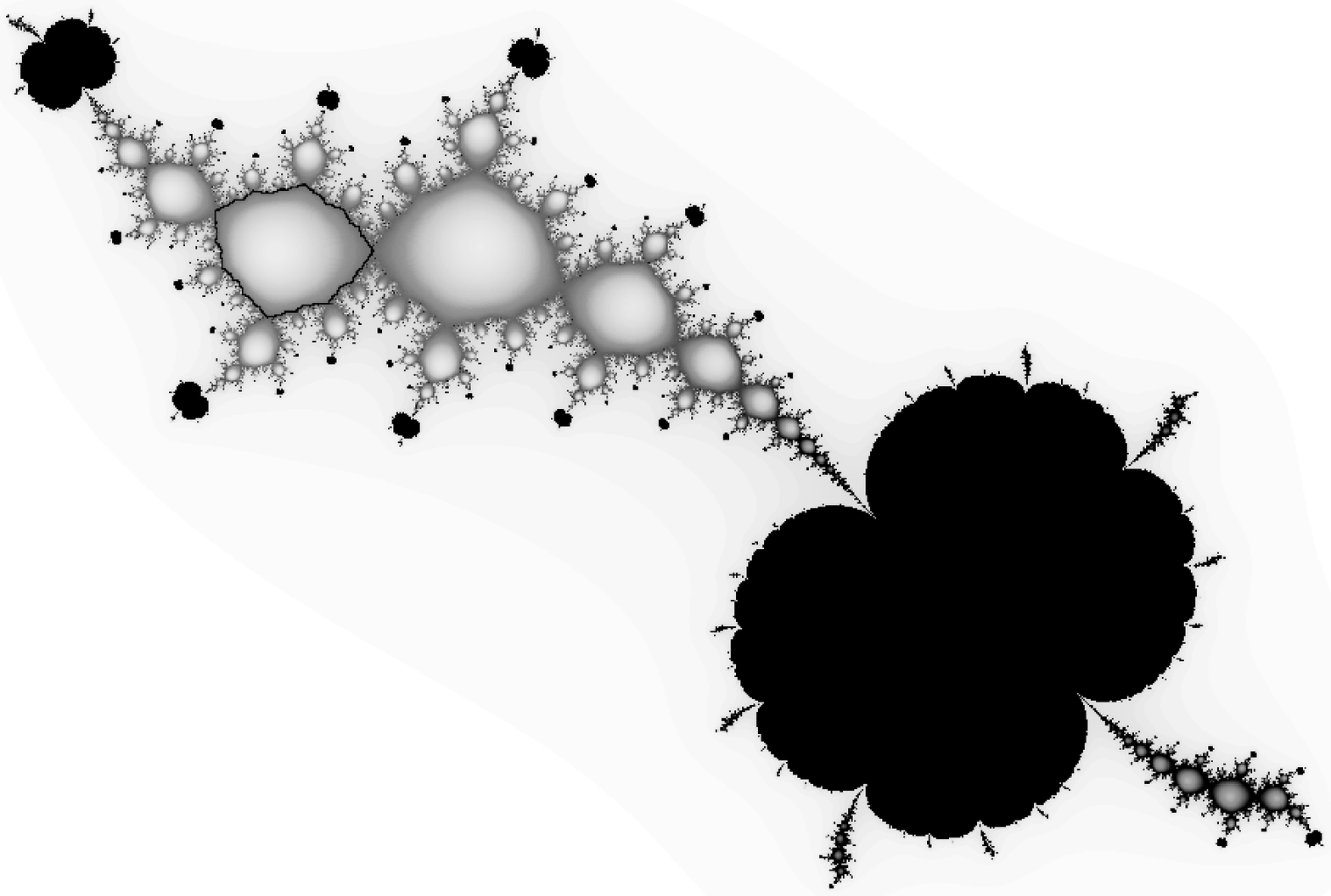}{{\sl The filled Julia set of a cubic on the
boundary of $\MM$ with $\theta=(\sqrt {5}-1)/2$. The large open topological
disk on the right is the immediate basin of attraction of a parabolic fixed
point. Here the parameter $c$ is located at the ``cusp'' of the large
cardioid on the right lower corner of \figref{m3}.}}{11cm} 
\realfig{2crit}{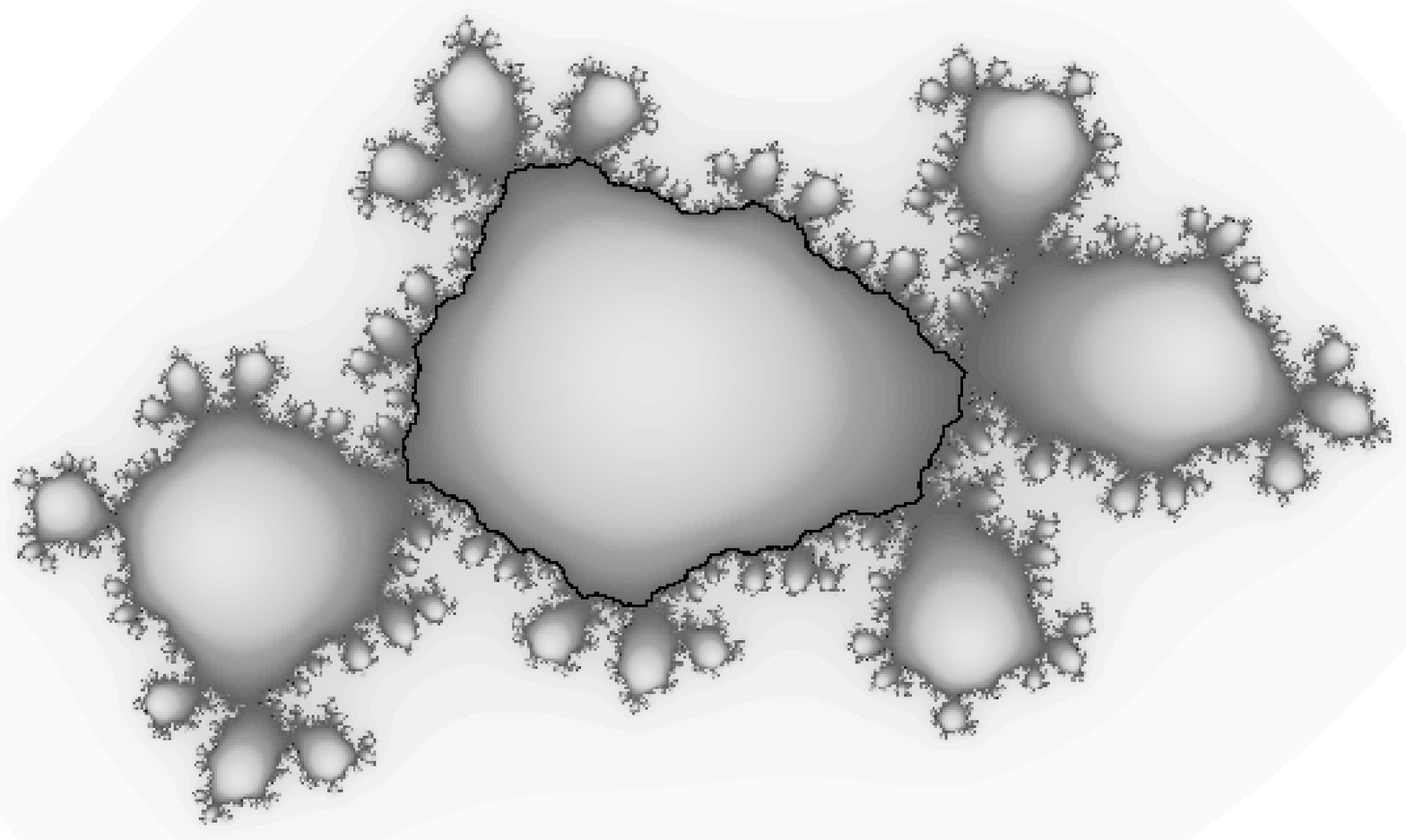}{{\sl Another example of the filled Julia set of a
cubic on the boundary of $\MM$ with $\theta=(\sqrt {5}-1)/2$. In this
example there are two distinct critical points on the boundary of the Siegel
disk centered at the origin. }}{11cm} 
\realfig{hypext}{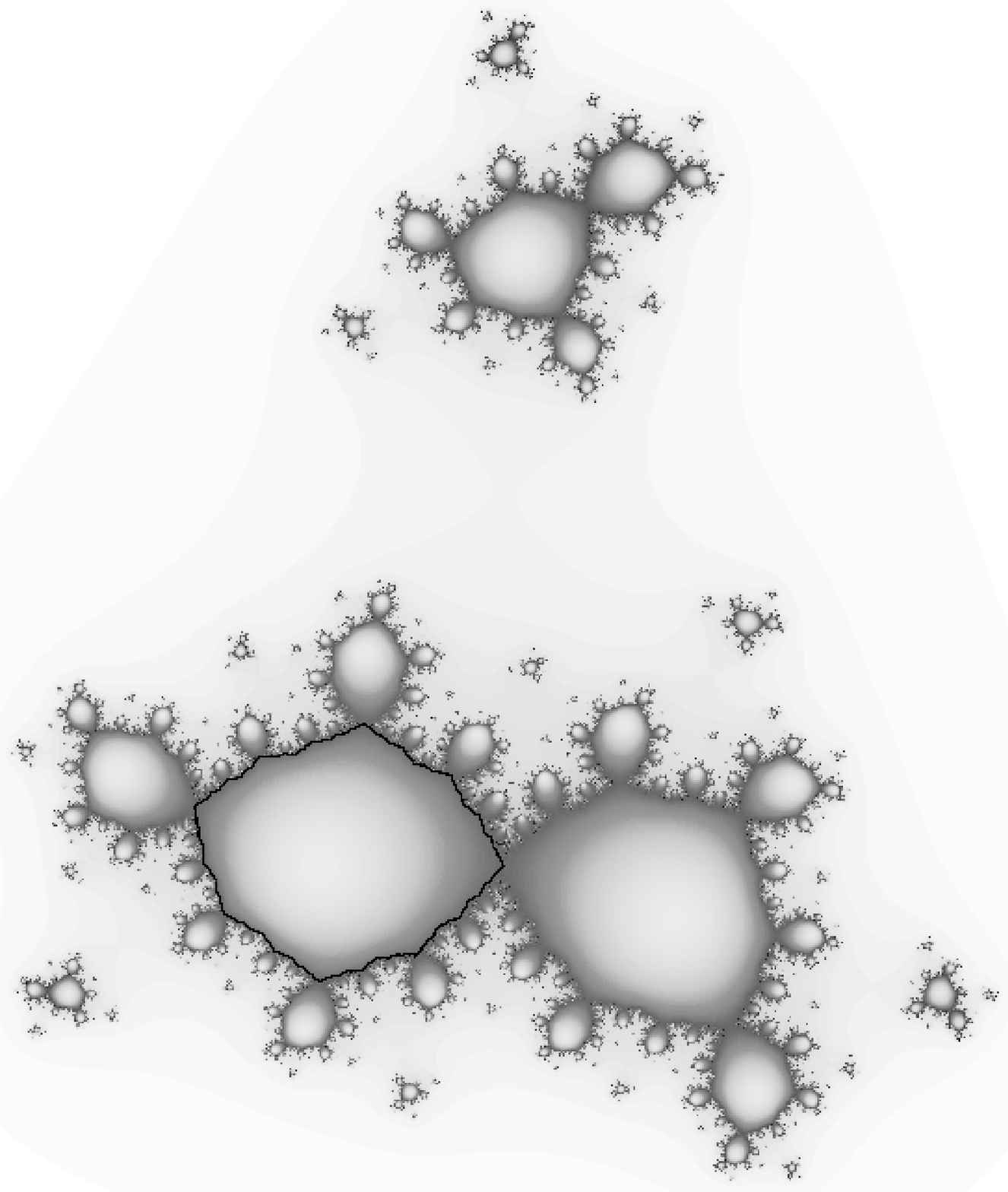}{{\sl The filled Julia set of a cubic in
$\Omega_{ext}$ for $\theta=(\sqrt {5}-1)/2$. It consists of uncountably many connected components. There are countably many components each quasiconformally homeomorphic to the quadratic Siegel filled Julia set with the same rotation number. Each remaining component is a single point.}}{11cm}

Let $P$ be hyperbolic-like and let $V$ and $U$ be the connected components containing $P$ of the hyperbolic-like cubics and the interior of $\MM$ respectively. Clearly $V\subset U$. If $V\neq U$, we can choose a $Q\in \partial V\cap U$. Since $Q$ is $J$-stable by \thmref{unstable}, the number of attracting cycles remains constant for all $Q'$ in a small neighborhood of $Q$. This number 
is $1$ if $Q'\in V$, hence $Q$ has to have an attracting cycle, contradicting $Q\in \partial V$.

Now consider the property of being capture for $P\in \PP$. It follows from 
\thmref{unstable} that when $P$ belongs to the interior of $\MM$, the connected component $V$ of the capture cubics containing $P$ has nonempty interior. Define $U$ as the component of the interior of $\MM$ containing $P$. Clearly $int(V)\subset U$. If they are not equal, let $Q\in \partial V \cap U$. Since $Q$ is
$J$-stable, for all $Q'$ in a small neighborhood of $Q$, a critical point of $Q'$ belongs to the Fatou set of $Q'$ if and only if the corresponding critical point of $Q$ belongs to the Fatou set of $Q$. If we choose $Q'\in V$, there is a critical point of $Q'$ which hits the Siegel disk $\Delta_{Q'}$. It follows that the same is true for $Q$, hence $Q$ is capture, which contradicts $Q\in \partial V$. This proves $int(V)=U$. In other words, when a capture cubic $P$ belongs to the interior of $\MM$, the whole interior component containing $P$ consists of capture cubics, hence the name ``capture component.''

However, since there is no a priori reason why the boundary of the Siegel disk of $P\in \PP$ should move continuously, being capture is not trivially seen to be an open condition. Hence, the above argument does not rule out the possibility of a capture being on the boundary of the connectedness locus $\MM$. In fact, this follows from a different type of argument which is standard in deformation theory for rational maps (see \thmref{capop}). On the other hand, when $\theta$ is of bounded type, we will show that the boundary of the Siegel disk of a cubic in $\PP$ moves continuously (see \thmref{move}), hence in this case the condition of being capture is automatically open. 

Conjecturally, queer components do not exist. But if they do, every cubic in a queer component exhibits an outstanding property: It admits an invariant line field on its Julia set, and in particular, its Julia set 
has positive Lebesgue measure. The proof of this fact depends on the harmonic $\lambda$-lemma of Bers and Royden \cite{Bers-Royden} as well as the elementary observation of Sullivan \cite{Sullivan2} that if the boundary of a Siegel disk moves holomorphically in a family of rational maps, then there is a choice of holomorphically varying Riemann maps for the Siegel disks (also see the new expanded version \cite{Mc-Sul}).

There is a technical difficulty showing up in the proof: For a general $\theta$ of Brjuno type, it is not known whether the boundary of the Siegel disk of a $P\in \PP$ is a Jordan curve. For this reason, the extension of holomorphic motions to the grand orbits of Siegel disks will require some extra work. On the other hand, we will prove later that for $\theta$ of bounded type, the boundary of the Siegel disk $\Delta_P$ of a $P\in \PP$ is a Jordan curve (see \thmref{main}). In this case the following theorem can be proved using the more elementary argument of \lemref{ext} with much less effort.

First we need the following useful lemma of L. Bers \cite{Bers}, \cite{Douady-Hubbard2}:

\begin{lem}[Bers Sewing Lemma]
\label{BSL}
Let $E\subset \BBB C$ be closed and $U$ and $V$ be two open neighborhoods of $E$. Let $\varphi :U \iso \varphi(U)$ and $\psi :V \iso \psi(V)$ be two homeomorphisms such that 
\begin{enumerate}
\item[$\bullet$]
$\varphi$ is $K_1$-quasiconformal,
\item[$\bullet$]
$\psi|_{V\smallsetminus E}$ is $K_2$-quasiconformal,
\item[$\bullet$]
$\varphi|_{\partial E}=\psi|_{\partial E}$.
\end{enumerate}
Then the map $\varphi \amalg \psi$ defined on $V$ by 
$$(\varphi \amalg \psi)(z)= \left \{ 
\begin{array}{ll}
\varphi(z) & z\in E \\
\psi(z)    & z\in V\smallsetminus E
\end{array}
\right. $$
is a $K$-quasiconformal homeomorphism with $K=\max \{ K_1, K_2 \} $. Moreover, $\overline{\partial}(\varphi \amalg \psi)=\overline{\partial} \varphi$ almost everywhere on $E$.
\end{lem}  

Throughout this paper, we say that two critically marked cubics $P_{c_1},P_{c_2}\in \PP$ are {\it quasiconformally conjugate} if there exists a quasiconformal homeomorphism $\varphi$ of the plane such that $\varphi \circ P_{c_1}=P_{c_2} \circ \varphi$ {\it and} $\varphi(c_1)=c_2$, In other words, all conjugacies must respect the markings of the critical points.
 
\begin{thm}[Invariant Line Fields for Queer Cubics]
\label{invline} 
Let $U$ be a queer component of the interior of $\MM$. Then for any $c\in U$, the Julia set $J(P_c)$ has positive Lebesgue measure and supports an invariant line field.
\end{thm}
\begin{pf}
The beginning of the argument is similar to the case of the Mandelbrot set. Fix some $c_0\in U$. We first note that every Fatou component of $P_{c_0}$ eventually maps to the Siegel disk $\Delta_{c_0}$ and the mapping is a conformal isomorphism: There cannot be further attracting cycles (since $P_{c_0}$ is not hyperbolic-like) or indifferent periodic orbits (see \thmref{indiff}). In particular, $K(P_{c_0})=\overline{GO(\Delta_{c_0})}$.
 
Choose some $c \in U$ with $c\neq c_0$, and let
$$\varphi_c:\BBB C\smallsetminus K(P_{c_0})\iso \BBB C\smallsetminus K(P_c)$$
be the conformal conjugacy given by composition of the B\"{o}ttcher maps of $P_{c_0}$ and $P_c$ (see Section \ref{sec:connect}). A brief computation shows that $\varphi_c(z)=\sqrt{c/c_0} z+O(1)$ and we can choose the branch of the square root near $c_0$ for which $\varphi_{c_0}(z)=z$. Since $\varphi_c$ depends holomorphically on $c$, 
it defines a holomorphic motion of $\BBB C\smallsetminus K(P_{c_0})$. By the harmonic $\lambda$-lemma \cite{Bers-Royden}, this motion extends to a {\it unique} holomorphic motion $\hat{\varphi}_c$ of the entire plane, which is now defined only for $c$ in a small neighborhood $V$ of $c_0$, with the following properties:
\begin{enumerate}
\item[$\bullet$]
For every $c\in V$, $\hat{\varphi}_c$ is a quasiconformal homeomorphism of the plane.
\item[$\bullet$] 
For every $c\in V$, the Beltrami differential $\displaystyle{\frac{\overline{\partial}\hat{\varphi}_c}{\partial \hat{\varphi}_c}\frac{d\overline z}{dz}}$ is harmonic in $GO(\Delta_{c_0})$. 
\end{enumerate}
It is easy to see that uniqueness of this extended motion implies that $\hat{\varphi}_c$ conjugates $P_{c_0}$ to $P_c$ on the entire plane (compare \cite{Mc-Sul}). In fact, one can replace $\hat{\varphi}_c$ by $P_c^{-1}\circ \hat{\varphi}_c \circ P_{c_0}$ on $GO(\Delta_{c_0})$, which also extends $\varphi_c$, where the branch of $P_c^{-1}$ is determined uniquely by the values of $\hat{\varphi}_c$ on the Julia set $J(P_{c_0})$. Hence $\hat{\varphi}_c=P_c^{-1}\circ \hat{\varphi}_c \circ P_{c_0}$ by uniqueness.

Next, we want to show that the restriction $\hat{\varphi}_c: GO(\Delta_{c_0})\rightarrow GO(\Delta_c)$ is a {\it conformal} conjugacy. As Sullivan observes in \cite{Sullivan2}, the fact that the boundary of $\Delta_{c_0}$ moves holomorphically (\thmref{unstable}) implies that there is a choice of the Riemann map $\zeta_c: \BBB D\rightarrow \Delta_c$ such that $\zeta_c(0)=0$ and $c\mapsto \zeta_c$ is holomorphic in $c$. Define a conformal conjugacy $\psi_c:\Delta_{c_0}\rightarrow \Delta_c$ by $\psi_c=\zeta_c\circ \zeta_{c_0}^{-1}$, which can be extended to a conformal conjugacy
$\psi_c : GO(\Delta_{c_0})\rightarrow GO(\Delta_c)$ by taking pull-backs as follows. Take any component $W$ of $P_{c_0}^{-n}(\Delta_{c_0})$ and let $W_c=\hat{\varphi}_c(W)$ be the corresponding component of $P_c^{-n}(\Delta_c)$. Define $\psi_c:W\rightarrow W_c$ by $\psi_c=P_c^{-n}\circ \psi_c \circ P_{c_0}^{\circ n}$. Since $c\mapsto \psi_c$ is holomorphic and $\psi_c=id$ when $c=c_0$, it follows that $\psi_c$ defines a holomorphic motion of $GO(\Delta_{c_0})$. By the harmonic $\lambda$-lemma, it extends to a unique holomorphic motion $\hat{\psi}_c$ of the entire plane which is defined for $c$ in a neighborhood $V'$ of $c_0$ and has harmonic Beltrami differential on $\BBB C\smallsetminus K(P_{c_0})$. By an argument similar to the one we used for $\hat{\varphi}_c$, it follows that $\hat{\psi}_c$ respects the dynamics, i.e., it conjugates $P_{c_0}$ to $P_c$ on the entire plane. In particular, it sends the marked critical point $c_0$ of $P_{c_0}$ to the marked critical $c$ of $P_c$. Let us assume for example that the forward orbit of $c_0$ accumulates on the boundary of $\Delta_{c_0}$. Then the same is true for $c$ and $\Delta_c$. Since $\hat{\varphi}_c$ was also a conjugacy to begin with, for all $c\in V \cap V'$ we have
$\hat{\psi}_c(c_0)=c=\hat{\varphi}_c(c_0)$, and by induction $\hat{\psi}_c(P_{c_0}^{\circ k}(c_0))=P_c^{\circ k}(c)=\hat{\varphi}_c(P_{c_0}^{\circ k}(c_0))$ for all $k$. Since every point on the boundary of $\Delta_{c_0}$ is in the closure of the forward orbit of $c_0$, we conclude that $\hat{\psi}_c$ and $\hat{\varphi}_c$ agree on $\partial \Delta_{c_0}$. Evidently this shows that $\hat{\psi}_c$ and $\hat{\varphi}_c$ agree on the boundary of every bounded Fatou component of $P_{c_0}$, hence on the entire Julia set $J(P_{c_0})$. It follows then from the Bers Sewing \lemref{BSL} that $\varphi_c \amalg \psi_c$ defined by 
$$(\varphi_c \amalg \psi_c)(z)=  \left \{ \begin{array}{ll} 
\hat{\varphi}_c(z) &  z\in \BBB C \smallsetminus GO(\Delta_{c_0})  \\
\hat{\psi}_c(z) & z\in GO(\Delta_{c_0})
\end{array}
\right.$$
is a quasiconformal homeomorphism which trivially has harmonic Beltrami differential in $\BBB C \smallsetminus J(P_{c_0})$. Note that $\varphi_c \amalg \psi_c$ is an extension of both $\varphi_c$ and $\psi_c$. By uniqueness, we conclude that $\hat{\varphi}_c \equiv \hat{\psi}_c$. In particular, when $c \in V \cap V'$, $\hat{\varphi}_c$ is conformal away from the Julia set $J(P_{c_0})$.

Now, if the Julia set $J(P_{c_0})$ had measure zero, $\hat{\varphi}_c$ would have been conformal, contradicting $c\neq c_0$. So $J(P_{c_0})$ has positive measure. The desired invariant line field is then given by $\hat{\varphi}_c^{\ast}(\sigma_0)$, the pull-back of the standard conformal structure $\sigma_0$ on the plane by $\hat{\varphi}_c$.
\end{pf}

It must be apparent that the existence of holomorphic motions in the above proof was the crucial fact which made the conformal extensions possible. In the case we have ``static'' quasiconformal conjugacies rather than holomorphic motions, such  conformal extensions are still possible once we assume that the boundaries of Siegel disks are Jordan curves. To show this, we first need the following definition: \\ \\
{\bf Definition (Conformal Position).} Let $\Delta $ be a Jordan domain containing the origin, with a marked point $b$ on its boundary. Consider the unique conformal isomorphism $\zeta:\Delta \iso \BBB D$ such that $\zeta(0)=0$ and $\zeta(b)=1$. By the {\it conformal position} of a point $z\in \Delta $ we mean the image $\zeta(z)\in \BBB D$. Note that this notion is well-defined once the boundary marking is given.\\
   
Let $R_t:z\mapsto e^{2 \pi i t}z$ denote the rigid rotation on the unit circle. Let $\zeta: \Delta \rightarrow \BBB D$ be any conformal isomorphism with $\zeta(0)=0$. Then a homeomorphism $h_t:\partial \Delta \rightarrow \partial \Delta$ of the form $h_t=\zeta^{-1}\circ R_t \circ \zeta$ is called an {\it intrinsic rotation} of $\partial \Delta $. By the Schwarz Lemma, $h_t$ is independent of the choice of $\zeta$.

Now consider two Jordan domains $\Delta_1$ and $\Delta_2$ containing the origin and a homeomorphism $\varphi: \partial \Delta_1 \rightarrow \partial \Delta_2$. Suppose that for some irrational angle $t$, $\varphi \circ h_t^1=h_t^2 \circ \varphi$, where the $h_t^j$ denote the intrinsic rotations of $\partial \Delta_j$. Then we can talk about two points in $\Delta_1$ and $\Delta_2$ having the same conformal position {\em even if there is no preferred choice for the marked points as before}. In fact, we may choose any $b\in \partial \Delta_1$ and $\varphi(b)\in \partial \Delta_2$ as the marked points on the boundaries and define the conformal positions accordingly. It is easy to check that the notion of having the same conformal position for two points in $\Delta_1$ and $\Delta_2$ does not depend on the particular choice of $b$.

For our purposes, $\Delta_1$ and $\Delta_2$ will be the Siegel disks of two cubics in $\PP$ and the homeomorphism $\varphi:\partial \Delta_1 \rightarrow \partial \Delta_2$ comes from a conjugacy between the cubics on the boundaries of these Siegel disks. 

The following result, which will turn out to be useful later (see \propref{independent}), tells us how to extend a quasiconformal conjugacy between two cubics in $\PP$ to the grand orbits of their Siegel disks.

\begin{lem}[Extending QC Conjugacies]
\label{ext}
 Let $P$ and $Q$ be two cubics in
$\PP$ such that the boundary of the Siegel disk $\Delta_P$ of $P$ is a Jordan curve. Let $\varphi:\BBB C \rightarrow \BBB C$ be a quasiconformal homeomorphism  whose restriction $\BBB C \smallsetminus GO(\Delta_P) \rightarrow \BBB C \smallsetminus GO(\Delta_Q)$ conjugates $P$ to $Q$. Then
\begin{enumerate}
\item[(a)]
If $P$ is not capture, there exists a quasiconformal homeomorphism $\psi:\BBB C\rightarrow \BBB C$ which conjugates $P$ and $Q$, which is conformal on $GO(\Delta_P)$ and agrees with $\varphi$ on $\BBB C \smallsetminus GO(\Delta_P)$.
\item[(b)]
If $P$ is capture, we can construct a $\psi$ as in (a) if and only if 
the captured images of the critical points of $P$ and $Q$ in $\Delta_P$ and $\Delta_Q$ have the same conformal position.
\end{enumerate}
\end{lem}
\begin{pf}
(a) When $P$ is not capture, the extension is easy to define.
Fix some $b_1\in \partial \Delta_P$ and let $b_2=\varphi(b_1)$. Consider conformal isomorphisms $\zeta_1: \Delta_P \iso \BBB D$ and $\zeta_2: \Delta_Q \iso
\BBB D$, with $\zeta_1(0)=0=\zeta_2(0)$ and $\zeta_1(b_1)=1=\zeta_2(b_2)$, which conjugate $P$ on $\Delta_P$
and $Q$ on $\Delta_Q$ to the rigid rotation $R_{\theta}:z\mapsto e^{2\pi i 
\theta}z$ on $\BBB D$. Since the boundaries of $\Delta_P$ and $\Delta_Q$ are Jordan curves, $\zeta_1$ and 
$\zeta_2$ extend homeomorphically to the closures. The composition $\psi=
\zeta_2^{-1}\circ \zeta_1:\Delta_P\rightarrow \Delta_Q$ is conformal and conjugates 
$P$ on $\Delta_P$ to $Q$ on $\Delta_Q$. Also $\psi(b_1)=\varphi(b_1)=b_2$ and by induction $\psi(P^{\circ k}(b_1)))=Q^{\circ k}(b_2)=\varphi(P^{\circ k}(b_1))$. Since the orbit of $b_1$ is dense on the boundary of $\Delta_P$, we have $\psi |_{\partial \Delta_P}=\varphi |_{\partial \Delta_P}$. Therefore, $\psi$ gives the required extension of $\varphi$ to the Siegel disk $\Delta_P$. It is now easy to extend $\psi$ to the grand orbit $GO(\Delta_P)$ as follows: $P^{\circ k}$ maps any component of $P^{-k}(\Delta_P)$ isomorphically onto $\Delta_P$. Hence we can define $\psi$ on any such component as the composition $Q^{-k}\circ
\psi |_{\Delta_P}\circ P^{\circ k}$, where the branch of $Q^{-k}$ is determined by the values of $\varphi$ on the Julia set $J(P)$. Clearly this composition is conformal inside this component and agrees with $\varphi$ on its boundary. The fact that $\psi$ defined this way is a quasiconformal homeomorphism follows from the Bers Sewing \lemref{BSL}, with $U=V=\BBB C$ and $E=\BBB C \smallsetminus GO(\Delta_P)$.
\realfig{confpos}{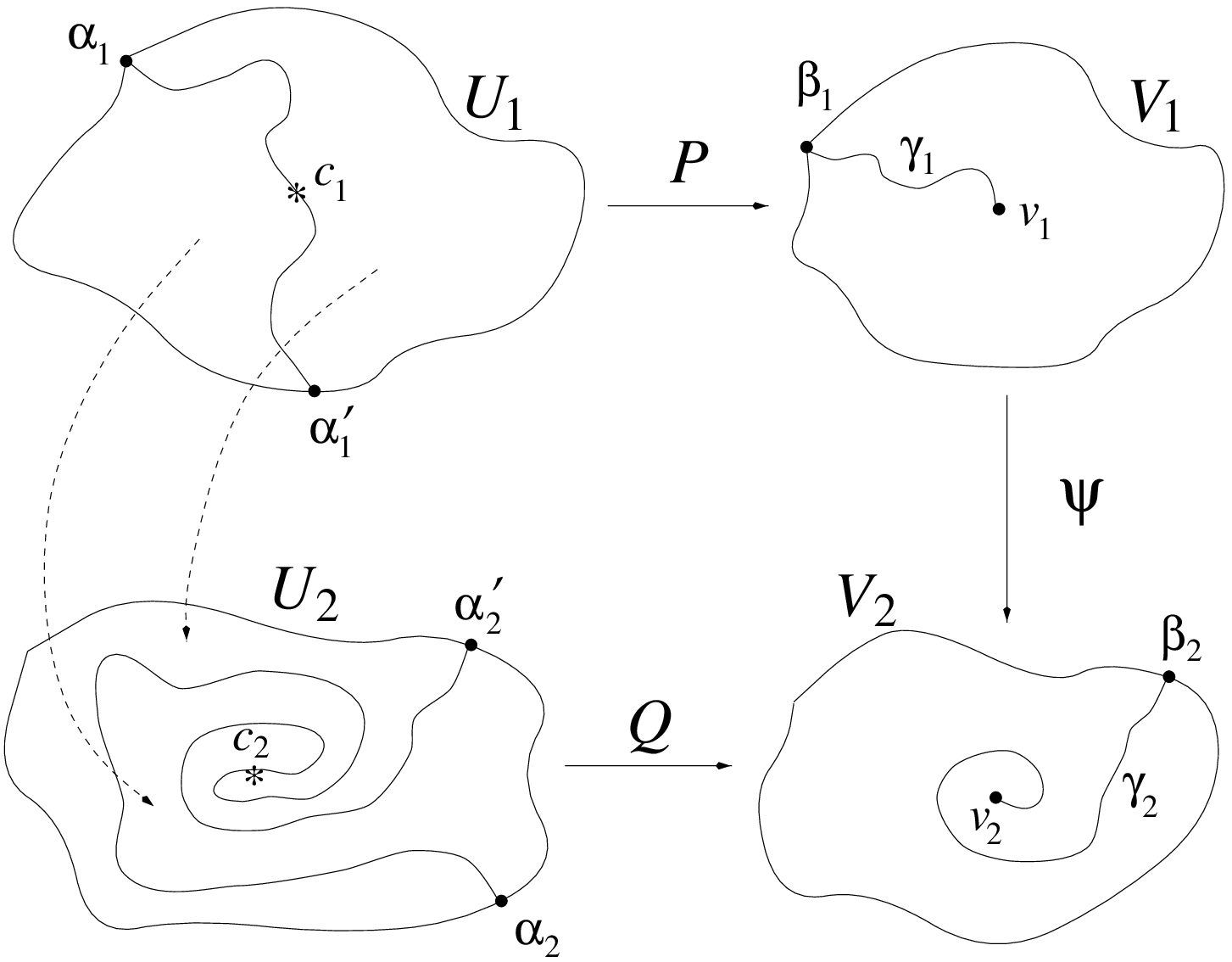}{{\sl Extending $\varphi$ in the capture case.}}{9cm}

(b) Now let $P$ be capture. The construction of $\psi$ goes through as in
case (a) except for the last part where we want to extend $\psi$ by taking pull-backs. Suppose that there exists a positive integer $k$ such that 
the critical point $c_1$ of $P$ belongs to the component $U_1$ of $P^{-k}(\Delta_P)$.
Let $V_1=P(U_1)$ and let $v_1=P(c_1)$ be the critical value in $V_1$. Since $P:\partial U_1 \rightarrow \partial V_1$ is a double covering and $\varphi$ conjugates $P$ to $Q$ on the Julia sets, there must be a critical point $c_2$ of $Q$ in a component $U_2$ of $Q^{-k}(\Delta_Q)$, with $\partial U_2=\varphi ( \partial U_1)$. Similarly define $V_2$ and $v_2$. By the proof of part (a) we can define $\psi$ inductively up to the $(k-1)$-th preimages of $\Delta_P$, including $V_1$. This gives us a conformal isomorphism $\psi:V_1\rightarrow V_2$ which necessarily maps $v_1$ to $v_2$, because by our assumption $P^{\circ k}(c_1)$ and $Q^{\circ k}(c_2)$ have the same conformal position in $\Delta_P$ and $\Delta_Q$ and so one gets
mapped to the other by $\psi |_{\Delta_P}$. Choose any simple arc $\gamma_1$ in $V_1$ connecting $v_1$ to some boundary point $\beta_1$. The simple arc $\gamma_2=\psi(\gamma_1)$ in $V_2$ connects $v_2$ to the boundary point 
$\beta_2=\psi(\beta_1)$. Pull $\gamma_1$ back by $P$ to get two branches of a simple arc passing through the critical point $c_1$ with two distinct endpoints $\alpha_1$ and $\alpha'_1$ on the boundary of $U_1$. Similarly we consider the pull-back of $\gamma_2$ by $Q$ and we get two endpoints on the boundary of $U_2$, which we label as $\alpha_2=\varphi(\alpha_1)$ and $\alpha'_2=\varphi(\alpha'_1)$ (see \figref{confpos}). Now the inverse $Q^{-1}$ can be defined analytically over $V_2 \smallsetminus \gamma_2$ and has two branches which take 
values in two different connected components of $U_2\smallsetminus Q^{-1}(\gamma_2)$. Define $\psi$ on $U_1$ as the composition $Q^{-1}\circ
\psi \circ P$, where the boundary orientation tells us which of the two branches of $Q^{-1}$ has to be taken. This way we extend $\psi$ to $U_1$ and  $\psi$ can then be defined on further preimages of $\Delta_P$ similar to the case (a).
\end{pf}
\vspace{0.17in}

\section{Renormalizable Cubics}
\label{sec:rencub}

This section briefly studies the class of renormalizable cubics in $\PP$. These are the cubics with disjoint critical orbits out of which one can extract the quadratic $\QQ$ by straightening. From a different point of view, one may consider a renormalizable cubic with connected Julia set as the result of ``intertwining'' the quadratic $Q_{\theta}$ with another quadratic with connected Julia set (compare \cite{Epstein-Yam}).

For background on polynomial-like maps, straightening, and hybrid classes, see for example \cite{Douady-Hubbard2}.\\ \\
{\bf Definition.} Let $P\in \PP$. We call $P$ {\it renormalizable} if there exists a pair of Jordan domains $U$ and $V$, with $0\in U\Subset V$, such that the restriction $P|_U: U\rightarrow V$ is a quadratic-like map hybrid equivalent to $\QQ$.\\

When $\theta$ is irrational of bounded type, it follows from \cite{Douady2} that the boundary of the Siegel disk of $Q_{\theta}$ is a quasicircle passing through 
the critical point. Hence the same is true for the Siegel disk $\Delta_P$ when $P$ is renormalizable. 
\begin{thm}
\label{ren}
A cubic $P\in \PP$ is renormalizable if one of the following conditions holds:
\begin{enumerate}
\item[(a)]
$P$ has a non-repelling periodic orbit other than $0$ which is not parabolic.\footnote{That the parabolic case must be excluded was pointed out to me by M. Yampolsky.}
\item[(b)]
$P$ has disconnected Julia set.
\end{enumerate}
\end{thm}

\begin{pf}
We use the Separation \lemref{seplem}. First assume that we are in case (a) so that $J(P)$ is connected. Let $\cal R$ be the finite collection of the closed preperiodic external rays given by the Separation Lemma. Let $V$ be the component of $\BBB C\smallsetminus \cal R$ which contains $0$, cut off by an equipotential of $K(P)$. Finally, let $U$ be the component of $P^{-1}(V)$ containing $0$. Since all the rays in $\cal R$ are preperiodic, $P(\cal R)\subset \cal R$, hence $U\subset V$.
$U$ necessarily contains a critical point of $P$ since otherwise the Schwarz lemma and $|P'(0)|=1$ would imply that $U=V$ and $P|_U:U\rightarrow V$ is a conformal isomorphism conjugate to a rotation. This would contradict the fact that $U$ intersects the basin of attraction of infinity for $P$. The other critical point of $P$ has to stay away from $V$ because by the second part of the Separation Lemma its entire orbit lives in the cycle of components of $\BBB C\smallsetminus \cal R$ which contains the non-repelling periodic orbit of $P$. 

Since by our assumption the non-repelling cycle of $P$ is not parabolic, the landing points of the external rays in $\cal R$ must all be repelling. Therefore, by a simple ``thickening'' procedure (see for example \cite{Milnor2}), we can assume that $\overline{U}\subset V$, so that $P|_U:U\rightarrow V$ is a quadratic-like map. Since up to affine conjugation there is only one quadratic polynomial which has a fixed Siegel disk of rotation number $\theta$, this quadratic-like map has to be hybrid equivalent to $\QQ$. This proves the theorem in the case $J(P)$ is connected. 

Now suppose that we are in situation (b) so that $J(P)$ is disconnected. For $\epsilon >0$, let $U_{\epsilon}$ be the connected component of $\{ z\in \BBB C: G_P(z)< \epsilon \}$ containing the Siegel disk $\Delta_P$, where $G_P:\BBB C \rightarrow \{ x\in \BBB R: x\geq 0 \}$ is the Green's function of $K(P)$. It is not hard to see that for small $\epsilon$, $P|_{U_{\epsilon}}: U_{\epsilon}\rightarrow U_{3\epsilon}$ is a quadratic-like map, necessarily hybrid equivalent to $Q_{\theta}$. 
\end{pf} 

\figref{hypcub} and \figref{hypext} demonstrate the above theorem. In each example, there is a piece of the filled Julia set which is quasiconformally homeomorphic to the filled Julia set of $\QQ$ in \figref{Q}. This piece is just the filled Julia set of the quadratic-like restriction $P|_U:U\rightarrow V$ given by the above theorem.
\realfig{Q}{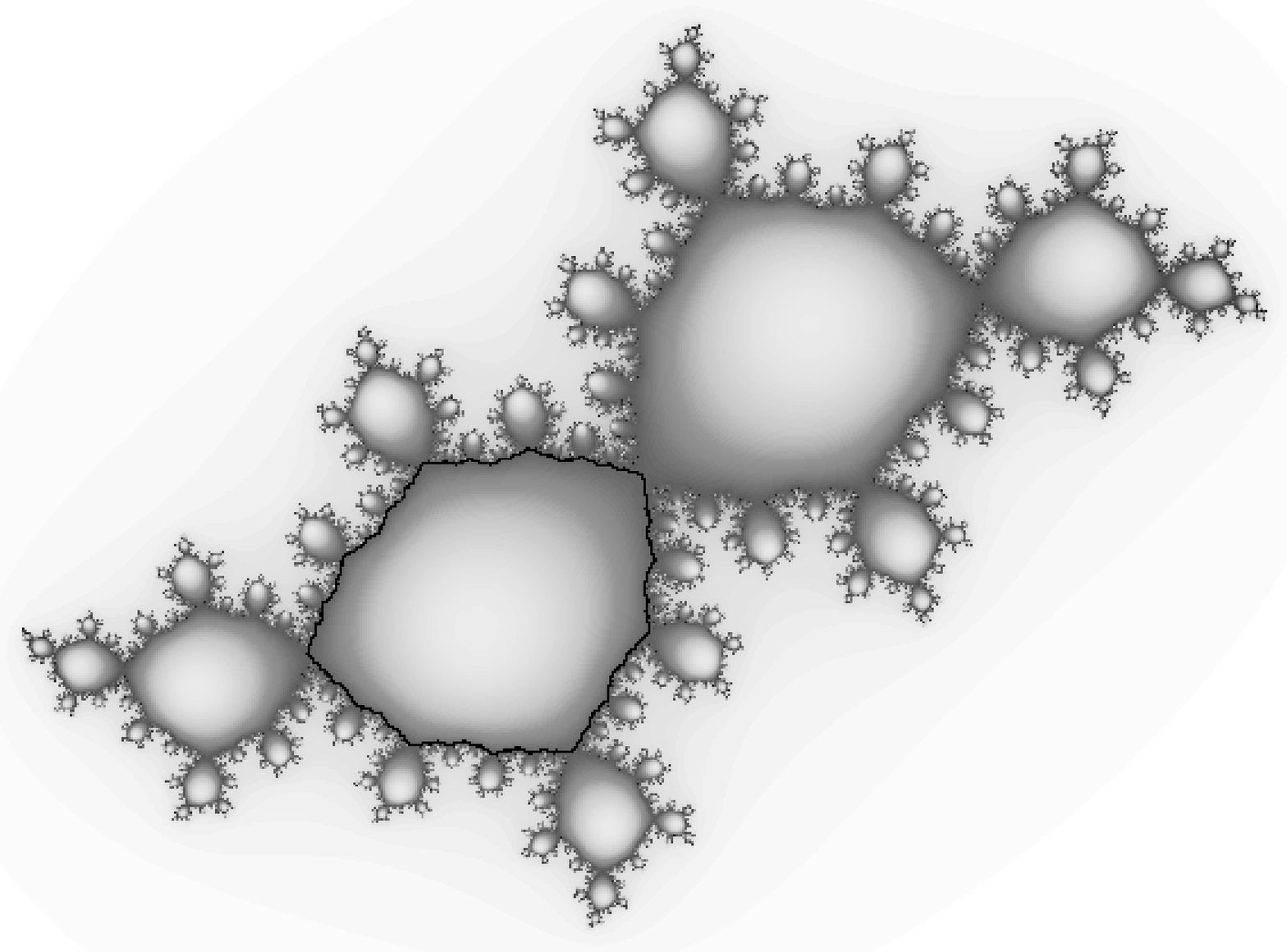}{{\sl The filled Julia set of the quadratic $\QQ$ for $\theta=(\sqrt{5}-1)/2$.}}{9cm}

\begin{cor}
\label{measure}
Let $\theta$ be an irrational number of bounded type. Let $P\in \PP$ be hyperbolic-like or have disconnected Julia set $J(P)$. Then $J(P)$ has Lebesgue measure zero.
\end{cor}
\begin{pf}
Let $P|_U:U\rightarrow V$ be the quadratic-like restriction given by \thmref{ren} with the filled Julia set $K$. Since this restriction is hybrid equivalent to $\QQ$ whose Julia set has measure zero when $\theta$ is bounded type \cite{Petersen}, we simply conclude that $\partial K$ has Lebesgue measure zero.

It is well-known that the forward orbit of almost every point $z\in J(P)$ accumulates on the $\omega$-limit set of the critical points of $P$ (\cite{Lyubich}, Proposition 1.14), which in this case is just $\partial \Delta_P$ union the attracting periodic orbit (resp. $\partial \Delta_P$) if $P$ is hyperbolic-like (resp. with disconnected Julia set). So the orbit of almost every $z\in J(P)$ accumulates on $\partial \Delta_P$. This implies that for all $n\geq N=N(z)$, $P^{\circ n}(z)\in V$. This can happen only if $P^{\circ N}(z)\in \partial K$ or equivalently $z\in P^{-N}(\partial K)$. We conclude that, up to a set of measure zero, $J(P)=\bigcup_{N\geq 0}P^{-N}(\partial K)$. But the right-hand side has measure zero because 
$ \partial K$ does. This proves that $J(P)$ has Lebesgue measure zero as well.
\end{pf} 

The next supplementary result will be useful later in the proof of connectivity of $\MM$ (\thmref{m3con}). I am indebted to M. Lyubich for pointing out that every quasiconformal self-conjugacy of the map $z\mapsto z^2$ near the unit circle $\BBB T$ extends to the identity map on $\BBB T$. This fact is the heart of the following lemma.
\begin{lem}
\label{EF}
Let $f:U\rightarrow V$ and $g:U'\rightarrow V'$ be quadratic-like maps both hybrid equivalent to the same quadratic polynomial $Q:z\mapsto z^2+c$ with connected Julia set. Let $E$ and $F$ be two subsets of $U$ and $U'$ respectively, such that
\begin{enumerate}
\item[$\bullet$]
$E\cap K(f)=\emptyset$ and $F\cap K(g)=\emptyset$,
\item[$\bullet$]
$E\cup K(f)$ and $F\cup K(g)$ are closed in $U$ and $U'$ respectively, and
\item[$\bullet$]
$f^{-1}(E)\subset E$ and $g^{-1}(F)\subset F$.
\end{enumerate}
Then any quasiconformal homeomorphism $\varphi:U\smallsetminus (E\cup K(f))\rightarrow U'\smallsetminus (F\cup K(g))$ which conjugates $f$ and $g$ extends to a quasiconformal homeomorphism $\varphi:U\smallsetminus E\rightarrow U'\smallsetminus F$. Moreover, we can arrange $\overline{\partial}\varphi=0$ on $K(f)$. 
\end{lem}
\begin{pf}
By straightening, we may assume without loss of generality that both $f$ and $g$ are the quadratic $Q$. Under this assumption, we prove that $\varphi$ extends continuously to the identity on the filled Julia set $K(Q)$. The last part of the theorem will follow because the $\overline{\partial}$-derivative of every hybrid equivalence vanishes on the corresponding filled Julia set.

Consider the B\"{o}ttcher map $\beta:\BBB C\smallsetminus K(Q)\iso \BBB C\smallsetminus \overline{\BBB D}$ which conjugates $Q$ to $z\mapsto z^2$ near $K(Q)$. Put $\tilde{U}=\beta(U\smallsetminus K(Q))$ and $\tilde{E}=\beta(E)$, and similarly define $\tilde{U'}$ and $\tilde{F}$. The induced map $\tilde{\varphi}=\beta \circ \varphi \circ \beta^{-1}:\tilde{U}\smallsetminus \tilde{E}\rightarrow \tilde{U'}\smallsetminus \tilde{F}$ is then a quasiconformal homeomorphism which satisfies $\tilde{\varphi}(z^2)=(\tilde{\varphi}(z))^2$. 

Consider the universal covering map $\zeta:\BBB H\rightarrow \BBB C\smallsetminus \overline{\BBB D}$ defined on the upper-half plane by $\zeta(z)=e^{-2 \pi i z}$. Let $\hat{U}=\zeta^{-1}(\tilde{U})$, $\hat{E}=\zeta^{-1}(\tilde{E})$, etc. Lift $\tilde{\varphi}$ to a quasiconformal homeomorphism $\hat{\varphi}:\hat{U}\smallsetminus \hat{E} \rightarrow \hat{U'}\smallsetminus \hat{F}$ which satisfies $\hat{\varphi}(z+1)=\hat{\varphi}(z)+1$ and $\hat{\varphi}(2z)=2\hat{\varphi}(z)$. Without loss of generality we can assume that $\hat{U}$ contains the horizontal strip $\{ z: 1\leq \Im (z) \leq 2 \}$. Clearly
$$\sup \{ d_{\BBB H}(z,\hat{\varphi}(z)): z\in \hat{U}\smallsetminus \hat{E}, 1\leq \Im (z)\leq 2 \} =C < +\infty, $$
where $d_{\BBB H}$ is the hyperbolic distance in $\BBB H$. Now given any point $z\in \hat{U}\smallsetminus \hat{E}$, choose $n\in \BBB Z$ so that $1\leq 2^n \Im (z)\leq 2$. Then 
$$d_{\BBB H}(z,\hat{\varphi}(z))=d_{\BBB H}(2^nz,2^n\hat{\varphi}(z))=d_{\BBB H}(2^nz,\hat{\varphi}(2^nz))\leq C.$$
By the Schwarz lemma applied to the composition $\beta^{-1}\circ \zeta$, we have
$$\sup \{ d(z,\varphi (z)): z\in U\smallsetminus (E\cup K(Q)) \}\leq C,$$
where $d$ is the hyperbolic distance in $\BBB C\smallsetminus K(Q)$. Hence, as $z\rightarrow J(Q)$ in $U\smallsetminus (E\cup K(Q))$, $|z-\varphi(z)|\rightarrow 0$. This means that we can define $\varphi(z)=z$ throughout $K(Q)$, and the extension will be a quasiconformal homeomorphism by the Bers Sewing Lemma.
\end{pf}
\vspace{0.17in}

\section{Connectivity of ${\cal M}_3(\theta)$}
\label{sec:connect}

In this section we prove that $\MM$ is connected. It will be more convenient to work with the double cover $\hat \MM$, which by definition is the set of all $s\in \BBB C^{\ast}$ such that $s^2\in \MM$. \propref{m3com} shows that the complement of $\hat \MM$ in $\BBB C^{\ast}$ has two connected components $\hat{\Omega}_{ext}$ and $\hat{\Omega}_{int}$ which are double covers of $\Omega_{ext}$ and $\Omega_{int}$ and are mapped to one another by the inversion $s\mapsto 1/s$. We would like to show that these open sets are homeomorphic to punctured disks. Connectivity of $\hat \MM$, hence $\MM$, will follow immediately. The strategy of the proof is more or less similar to the proof of connectivity of the Mandelbrot set with one additional
difficulty: We construct a holomorphic branched covering $\Phi:\hat{\Omega}_{ext}\rightarrow {\BBB C} \smallsetminus \overline {\BBB D}$ which extends holomorphically to infinity with $\Phi^{-1}(\infty)=\infty$. The degree of this map is 3, so to prove that $\hat{\Omega}_{ext}$ is a punctured disk one has to show that $\Phi$ has no critical point other than $\infty$. This additional difficulty does not show up in the case of the Mandelbrot set, where the similar map has degree 1, so it automatically becomes a conformal isomorphism (see \cite{Douady-Hubbard1}).

Recall that the {\it B\"{o}ttcher map} $\beta$ associated to a polynomial $$P:z\mapsto a_dz^d+a_{d-1}z^{d-1}+\cdots+a_1z+a_0, \ \ \ \ a_d\neq 0$$ 
is a conformal isomorphism  defined near $\infty$, with $\beta(\infty)=\infty$, which conjugates $P$ to 
the map $z\mapsto z^d$; that is, $\beta(P(z))=\beta(z)^d$. This map is unique up to multiplication by a $(d-1)$-th root of unity, so it can be normalized so that the derivative at infinity $\beta'(\infty)$ becomes any $(d-1)$-th root of $1/a_d$.

There is a classical formula for $\beta$ in terms of the iterates of the polynomial $P$ (see for example \cite{Carleson-Gamelin}):
\begin{equation}
\label{eqn:zarb}
\beta(z)=\lim_{n\rightarrow \infty} \left ( P^{\circ n}(z) \right ) ^{d^{-n}} = z \prod_{n=1}^{\infty} \left ( \frac{P^{\circ n}(z)}{(P^{\circ n-1}(z))^d} \right )^{d^{-n}}.
\end{equation}
The infinite product converges uniformly outside a sufficiently large disk centered at the origin.

For each $s\in {\BBB C}^{\ast}$, consider the polynomial 
$$P^s: z \mapsto \lambda z \left ( 1-\frac{1}{2}(s+\frac{1}{s})z+\frac{1}{3} z^2 \right ) $$
which has critical points at $s$ and $1/s$. The dilation $z\mapsto s z$ conjugates $P^s$ to $P_c$ in (\ref{eqn:normform}) with $c=s^2$. Hence $P^s \in \hat \MM$ if and only if $P_c \in \MM$.
 
\begin{thm}[Connectivity of $\MM$]
\label{m3con}
The open set $\hat{\Omega}_{ext}$ is homeomorphic to ${\BBB C}\smallsetminus \overline {\BBB D}$. Therefore, $\hat \MM$, hence $\MM$, is connected.
\end{thm}
\begin{pf}
Let $s\in \hat{\Omega}_{ext}$ and let $\beta_s$ be the B\"{o}ttcher
map which conjugates $P^s$ to $z\mapsto z^3$ near infinity, with $\beta'_s(\infty)=\sqrt{3/\lambda}$. It is a standard argument to show that 
$\beta_s$ depends holomorphically on $s$ and can be extended conformally down to the equipotential $\gamma_s$ passing through the escaping critical point $s$ of $P^s$ and it maps outside of $\gamma_s$ to the outside of some closed disk $\overline{\BBB D}(0,r)$, where $r>1$. Note that $\gamma_s$ is topologically a figure eight with $s$ as double point. Define a map $\Phi:\hat{\Omega}_{ext}\rightarrow {\BBB C} \smallsetminus \overline {\BBB D}$ by
$$\Phi (s)=\beta_s (P^s(s)).$$
This is a holomorphic map which extends holomorphically to infinity. It is not hard to show that $\Phi$ is proper, i.e., $|\Phi(s)|\rightarrow 1$ as $s\rightarrow \partial \hat\MM$. Hence $\Phi$ is a finite-degree branched covering from $\hat{\Omega}_{ext}\cup \{ \infty \}$ to the topological disk 
$\overline {\BBB C} \smallsetminus \overline {\BBB D}$. Let us compute the mapping degree of $\Phi$. By (\ref{eqn:zarb}), we have
$$\begin{array}{rl}
\beta_s(z) & =z\ \displaystyle{ \prod_{n=1}^{\infty} \left [ \frac{\lambda}{3}-\frac{\lambda}{2}(s+\frac{1}{s})\frac{1}{(P^s)^{\circ n-1}(z)}+\frac{\lambda}{((P^s)^{\circ n-1}(z))^2} \right ] ^{3^{-n}} } \\ \\
           & =: z\ \displaystyle{ \prod_{n=1}^{\infty} \beta_n(z,s)^{3^{-n}} },
\end{array}$$
hence
\begin{equation}
\label{eqn:fi}
\Phi(s)=P^s(s) \prod_{n=1}^{\infty} \beta_n(P^s(s),s) ^ {3^{-n}},
\end{equation}
where $P^s(s)=-\frac{\lambda}{6}s^3+\frac{\lambda}{2}s$.
By considering the logarithm of $\Phi$, we see that near infinity the infinite product in (\ref{eqn:fi}) is of the form $\sqrt{\frac{\lambda}{3}}(1+O(1/s))$. Hence $\Phi(s)=\sqrt{\frac{\lambda}{3}}(-\frac{\lambda}{6}s^3+\frac{\lambda}{2}s)(1+O(1/s))$. 
Since $\Phi^{-1}(\infty)=\infty$, this means that the mapping degree of $\Phi$ is 3. 
In particular, $\infty$ is a double critical point of $\Phi$. By the Riemann-Hurwitz formula, the Euler characteristic of $\hat{\Omega}_{ext}$ is equal to $-n$, where $n$ is the number of critical points of $\Phi$ in $\hat{\Omega}_{ext}$. Therefore,
$\hat{\Omega}_{ext}$ is homeomorphic to a punctured open disk if and only if $\Phi$ has no critical point in $\hat{\Omega}_{ext}$. In what follows, we prove that $\Phi$ is locally injective in $\hat{\Omega}_{ext}$. Since $\Phi$ is also holomorphic, this will prove that there are no critical points other than $\infty$.

So assume $\Phi(s_1)=\Phi(s_2)=w$ for some $s_1,s_2\in \hat{\Omega}_{ext}$.
To simplify the notation, put $P_1=P^{s_1}, P_2=P^{s_2}$. Let $v_1=P_1(s_1)$ and $v_2=P_2(s_2)$ be the critical values and $a_1$ and $a_2$
be the {\it co-critical} points, i.e., $a_i\neq s_i$ and $P_i(a_i)=v_i$ for 
$i=1,2$. Finally, let $\gamma_1$ and $\gamma_2$ be the equipotentials of the corresponding B\"{o}ttcher maps $\beta_1$ and $\beta_2$ which pass through the critical points $s_1$ and $s_2$ (see \figref{conn}).
\realfig{conn}{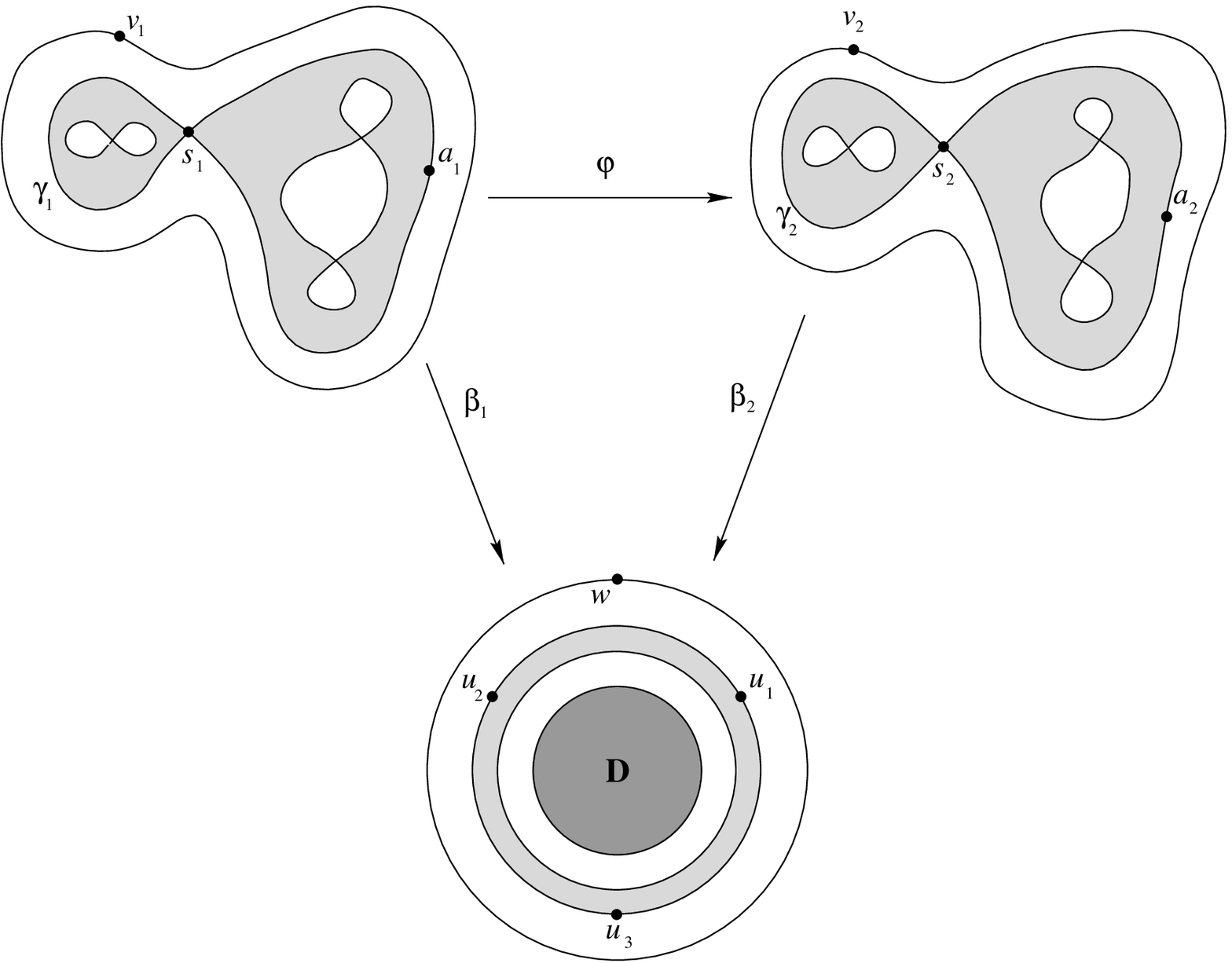}{{\sl Extending conformal conjugacies.}}{10cm}

Define a conformal map $\varphi$ from the outside of $\gamma_1$ to the outside of $\gamma_2$ by $\varphi=\beta_2^{-1}\circ \beta_1$. We would like to extend $\varphi$ to the entire basin of attraction of $\infty$ for $P_1$. Let $u_1, u_2, u_3$ be the cube roots of $w$. Under $\beta_1$, every connected component of $\gamma_1\smallsetminus \{ s_1, a_1 \}$ maps homeomorphically to 
one of the 1/3-circles joining $u_1, u_2, u_3$. If $s_1$ is sufficiently close to $s_2$, it follows by continuity that the corresponding components of  $\gamma_2\smallsetminus \{ s_2, a_2 \}$ will map homeomorphically to the {\it same} circular segments joining $u_1,u_2,u_3$ (see \figref{conn}). This allows us to extend $\varphi$ to a homeomorphism $\gamma_1\stackrel{\simeq}{\rightarrow}\gamma_2$.
Now it is straightforward to extend $\varphi$ further: The annulus bounded
by $\gamma_1$ and $P_1(\gamma_1)$ has two preimages which are mapped onto it in a 1-to-1 and 2-to-1 fashion. We can extend $\varphi$ to these
preimages by taking pull-backs, i.e., we define $\varphi$ to be 
$P_2^{-1} \circ \varphi \circ P_1$, where the boundary values of $\varphi|_{\gamma_1}:\gamma_1\rightarrow \gamma_2$ tell us which branch of $P_2^{-1}$ must be taken. It is not hard to see that this process of taking pull-backs can be continued until $\varphi$ is defined on the entire basin of attraction of $\infty$ for $P_1$. (One formal way to keep track of various preimages of these annuli is to consider the {\it pattern} associated with each cubic as introduced by Branner and Hubbard \cite{Branner-Hubbard}. In their language, $P_1$ and $P_2$ have ``homeomorphic patterns of infinite depth.'') The extension of $\varphi$ defined this way is conformal, since it is a homeomorphism which is holomorphic except on a disjoint countable union of piecewise analytic curves.

By \thmref{ren}, $P_1$ and $P_2$ are both renormalizable. Hence there are quadratic-like restrictions $f=P_1|_U:U\rightarrow V$ and $g=P_2|_{U'}:U'\rightarrow V'$ of both polynomials which are hybrid equivalent to $\QQ$. Note that $K(f)$ and $K(g)$ are just the components of $K(P_1)$ and $K(P_2)$ which contain the Siegel disks $\Delta_{P_1}$ and $\Delta_{P_2}$, respectively. The filled Julia set $K(P_1)$ decomposes as $\bigcup_{n\geq 0}P_1^{-n}(K(f)) \cup G_1$, where $G_1 \subset J(P_1)$ is the uncountable union of trivial components (see \figref{hypext}). Similarly, we have $K(P_2)=\bigcup_{n\geq 0}P_2^{-n}(K(g)) \cup G_2$. Now we are exactly in the situation of \lemref{EF}, with $E=(K(P_1) \smallsetminus K(f)) \cap U$ and $F=(K(P_2) \smallsetminus K(g)) \cap U'$. By \lemref{EF}, $\varphi$ can be extended to $K(f) \iso K(g)$ and then to $\bigcup_{n\geq 0}P_1^{-n}(K(f)) \iso \bigcup_{n\geq 0}P_2^{-n}(K(g))$ by taking pull-backs. Note that this extension has zero $\overline{\partial}$-derivative on this union by \lemref{EF} and the Bers Sewing Lemma. It is not hard to see that $\varphi$ also extends homeomorphically to $G_1 \iso G_2$. Therefore, we obtain a homeomorphism $\varphi: \BBB C \rightarrow \BBB C$ conjugating $P_1$ to $P_2$, which is quasiconformal at least on $\BBB C \smallsetminus J(P_1)$.

 We would like to show that $\varphi$ is quasiconformal everywhere. One way to do this proceeds as follows. By \thmref{qcclass} below, $P_1$ and $P_2$ are conjugate by a quasiconformal homeomorphism $\psi: \BBB C \rightarrow \BBB C$. The proof of \lemref{EF} shows that $\psi^{-1}\circ \varphi$ is the identity map on $\partial K(f)$. It easily follows that $\varphi=\psi$ on the entire Julia set $J(P_1)$. Therefore, $\varphi$ is identically equal to the homeomorphism $\varphi \amalg \psi$ defined by 
$$(\varphi \amalg \psi)(z)= \left \{ 
\begin{array}{ll}
\varphi(z) & z\in \BBB C \smallsetminus J(P_1) \\
\psi(z)    & z\in J(P_1)
\end{array}
\right. $$
which is quasiconformal by the Bers Sewing Lemma.

Finally, we show that $\varphi$ is in fact a conformal homeomorphism. Just as in the proof of \corref{measure}, up to a set of measure zero, $K(P_1)=\bigcup_{n\geq 0}P_1^{-n}(K(f))$. Hence the measure of $G_1$ is zero. It follows that $\varphi$ is conformal on $\BBB C\smallsetminus K(P_1)$ and has zero $\overline{\partial}$-derivative almost everywhere on $K(P_1)$. Hence $\varphi$ is conformal everywhere, which means $P_1=P_2$. 
\end{pf}
\vspace{0.17in}

\section{Cubic Quasiconformal Conjugacy Classes}
\label{sec:qc}

In this section we prove that quasiconformal conjugacy classes in $\PP$ are either open and connected or single points. This result, together with the fact that any holomorphic family of rational maps with constant critical orbit relations forms a quasiconformal conjugacy class, enables us to completely characterize the quasiconformal conjugacy classes in $\PP$.
\begin{thm}[Parametrization of QC Conjugacy Classes] 
\label{qcpar}
Let $P_{c_0},P_{c_1}$ be distinct cubics in $\PP$ and let $\varphi:\BBB C \rightarrow \BBB C$ be a $K$-quasiconformal homeomorphism which conjugates $P_{c_0}$ to $P_{c_1}$, i.e., $\varphi\circ P_{c_0}=P_{c_1}\circ \varphi$ and $\varphi(c_0)=c_1$. Then there exists a nonconstant holomorphic map $t\mapsto c_t$ from an open disk $\BBB D(0,r)\ (r>1)$ into $\BBB C^{\ast}$ which maps $0$ to $c_0$ and $1$ to $c_1$, such that for every $t\in \BBB D(0,r)$, $P_{c_0}$ is conjugate to $P_{c_t}$ by a $K_t$-quasiconformal homeomorphism $\varphi_t: \BBB C\rightarrow \BBB C$. Moreover, $K_t\rightarrow 1$ as $t\rightarrow 0$.
\end{thm}
\begin{pf}
The idea of the proof goes back to Douady and Hubbard \cite{Douady-Hubbard2}: Define a conformal structure $\sigma$ on $\BBB C$ by $\sigma=\varphi^{\ast}\sigma_0$, where, as usual, $\sigma_0$ is the standard conformal structure on $\BBB C$. (To simplify the notation, in what follows we identify a conformal structure on $\BBB C$ with its associated Beltrami differential.) Since $P_{c_1}$ is holomorphic, $P_{c_0}$ has to preserve $\sigma$. Since $\varphi$ is quasiconformal, $\| \sigma \|_{\infty}<1$. Define a one-parameter family $\{ \sigma_t \}$ of complex-analytic deformations of $\sigma$ by $\sigma_t=t\sigma$,
where $t\in \BBB D(0,r)$ for some $r>1$ such that $r\| \sigma \| _{\infty}<1$.
By the Measurable Riemann Mapping Theorem of Ahlfors and Bers \cite{Ahlfors-Bers}, there exists a unique quasiconformal homeomorphism $\varphi_t$ of the plane which solves the Beltrami equation $\varphi_t^{\ast} \sigma_0=\sigma_t$ and fixes $0$, $1$ and $\infty$. Define $P^t=\varphi_t\circ P_{c_0}\circ \varphi_t^{-1}$. Since $P_{c_0}$ is holomorphic, it acts as a pure rotation on Beltrami differentials. Hence $P_{c_0}^{\ast}\sigma=\sigma$ implies $P_{c_0}^{\ast}\sigma_t=\sigma_t$ and therefore $P^t$ is a quasiregular self-map of the plane which preserves $\sigma_0$ and is conjugate to a cubic polynomial. It is then easy to see that $P^t$ itself is a cubic polynomial with a fixed Siegel disk of rotation number $\theta$ centered at $0$ with a marked critical point at $z=1$.

Note that $t\mapsto \sigma_t$ is holomorphic, so the same is true for $t\mapsto \varphi_t$ and hence $t\mapsto P^t$ by the analytic dependence of the solutions of the Beltrami equation on parameters \cite{Ahlfors-Bers}. Therefore the map $t\mapsto c_t$ which defines the second critical point of $P^t$ so that $P^t=P_{c_t}$ is holomorphic. It is easy to see that $c_t$ has all the required properties.      
\end{pf}
\begin{cor}
\label{qcopen}
Quasiconformal conjugacy classes in $\PP$ are either single
points or open and connected. In particular, cubics on the boundary $\partial \MM$ are quasiconformally rigid, i.e., their conjugacy classes are single points.\ $\Box$
\end{cor}

\begin{thm}[Capture is an open condition]
\label{capop}
Let $P_{c_0} \in \PP$ be a capture cubic. Then there is an open neighborhood $U$ of $c_0$ such that for every $c\in U$, $P_c$ is also capture.
\end{thm}

\begin{pf}
When $\theta$ is of bounded type, we will see that the boundary of the Siegel disks of cubics in $\PP$ move continuously (see \thmref{move}) and in this case the theorem follows immediately. The following proof uses a standard argument in quasiconformal deformation theory which is similar to the proof of \thmref{qcpar}. I am indebted to X. Buff who pointed out to me that a defomation approach would work in the general case: To fix the ideas, let us assume that $P_{c_0}^{\circ k}(c_0)\in \Delta_{c_0}$ and $k\geq 1$ is the smallest such integer. First assume that $P_{c_0}^{\circ k}(c_0) \neq 0$. Let $A\subset \Delta_{c_0}$ be the annulus bounded by $\partial \Delta_{c_0}$ and the analytic invariant curve in $\Delta_{c_0}$ passing through $P_{c_0}^{\circ k}(c_0)$. Take a conformal isomorphism $\psi: A \iso {\BBB A}(1,\epsilon)$, with $\epsilon=e^{2 \pi mod(A)}>1$, which conjugates $P_{c_0}$ on $A$ to the rotation on ${\BBB A}(1,\epsilon)$. Postcompose $\psi$ with a (nonconformal) dilation ${\BBB A}(1,\epsilon) \rightarrow {\BBB A}(1,\epsilon^2)$ to get a quasiconformal homeomorphism $\varphi: A \rightarrow {\BBB A}(1,\epsilon^2)$ conjugating $P_{c_0}$ to the rotation. Define a $P_{c_0}$-invariant conformal structure $\sigma$ on $\BBB C$ by putting $\sigma=\varphi^{\ast}\sigma_0$ on $A$ and pulling it back by the inverse branches of $P_{c_0}$ to the entire grand orbit of $A$. Set $\sigma=\sigma_0$ everywhere else. As in the proof of \thmref{qcpar}, we define $\sigma_t=t\sigma$ for $t\in {\BBB D}(0,r)$ for some $r>1$, solve the Beltrami equation $\varphi_t^{\ast}\sigma_0=\sigma_t$ and set $P^t=\varphi_t \circ P_{c_0} \circ \varphi_t^{-1}$. Then $P^t$ is a capture cubic in $\PP$ and $P^0=P_{c_0}$. The holomorphic mapping $t\mapsto P^t$ is not constant because $mod(\varphi_1(A))$ is the same as the modulus of $A$ equipped with the conformal structure $\sigma$, which in turn is $(1/2 \pi) \log(\epsilon^2)=2\ mod(A)$. Hence $P^1\neq P^0$ and the mapping $t\mapsto P^t$ is open.

Now consider the case where $P_{c_0}^{\circ k}(c_0)=0$. In this case, by \corref{lower bound}, the conformal capacity of $\Delta_c$ has a positive lower bound for all $c$ sufficiently close to $c_0$. It follows that there exists an $\epsilon >0$ such that for all $c$ close to $c_0$, $\Delta_c \supset {\BBB D}(0, \epsilon)$. Hence a small perturbation of $P_{c_0}$ will still be a capture cubic.     
\end{pf} 

By a {\it center} of a hyperbolic-like component $U\subset\MM$ we mean a cubic $P_c\in U$ with one of the critical points $c$ or $1$ being periodic. Similarly, a center of a capture component will be a cubic with one critical point eventually mapped
to the indifferent fixed point at the origin.
\begin{lem}[Existence of Centers]
\label{cen}
Every hyperbolic-like or capture component of the interior of $\MM$ has a center. 
\end{lem}

By the remark after the proof, centers of hyperbolic-like or capture components are unique when $\theta$ is of bounded type.

\begin{pf}
First let $U$ be a hyperbolic-like component. For every $c\in U$, consider the multiplier $m(c)$ of the unique attracting periodic orbit of $P_c$. The mapping $c\mapsto m(c)$ from $U$ into $\BBB D$ is easily seen to be proper and holomorphic. Hence it vanishes at a finite number of points in $U$.

Now let $U$ be capture. To be more specific, let us assume that for every $c\in U$, $P_c^{\circ k}(c)$ belongs to the Siegel disk $\Delta_c$, and let $k$ be the smallest such integer. Since $P_c$ is $J$-stable by \thmref{unstable}, the boundary of $\Delta_c$ moves holomorphically. Then, as in the proof of \thmref{invline}, there is a holomorphically varying choice of the Riemann maps $\zeta_c:\BBB D \rightarrow \Delta_c$ with $\zeta_c(0)=0$. Define a map $m:U\rightarrow \BBB D$ by
$$m(c)=\zeta_c^{-1}(P_c^{\circ k}(c)).$$
(In the language of the definition before \lemref{ext}, this is just the conformal position with respect to $\zeta_c$ of the captured image of the critical point $c$ of $P_c$.)
Clearly $m$ is holomorphic. Let $c_n \in U$ be any sequence which converges to $c\in \partial U$ as $n\rightarrow \infty$. For simplicity, put $\zeta_{c_n}=\zeta_n$. Let $z_n=P_{c_n}^{\circ k}(c_n)\in \Delta_{c_n}$ and $w_n=m(c_n)=\zeta_n^{-1}(z_n)\in \BBB D$. If $w_n$ does not converge to the unit circle, we can find a subsequence $w_{n(j)}$ such that $w_{n(j)}\rightarrow w\in \BBB D$ as $j\rightarrow \infty$. Since the family of univalent functions $\{ \zeta_n :\BBB D\rightarrow \BBB C \}$ is normal, by passing to a further subsequence if necessary, we may assume that $\zeta_{n(j)}\rightarrow \zeta$ locally uniformly on $\BBB D$. Clearly $\zeta (\BBB D)\subset \Delta_c$. Therefore, $\zeta(w)=\lim_j \zeta_{n(j)}(w_{n(j)})=\lim_j z_{n(j)}=P_c^{\circ k}(c)\in \Delta_c$. But this means that $P_c$ is capture, which contradicts $c\in \partial U$. This proves that $w_n$ converges to the unit circle. Hence $m$ is a proper map. Now, as before, $m^{-1}(0)$ has to be nonvacuous and finite.
\end{pf}
\vspace{0.15in}
\noindent
{\bf Remark.} When the rotation number $\theta$ is of bounded type, there is a simple proof for the {\it uniqueness} of centers. (Compare \cite{MCM} or \cite{JACK}, where this is shown for every hyperbolic component in a holomorphic family of polynomial maps.) We sketch such a proof briefly. By \corref{qcopen}, it is enough to prove that any two centers for a component are quasiconformally conjugate. First let $U$ be a hyperbolic-like component and $c_1$ and $c_2$ be centers of $U$. Let $P_i=P_{c_i}$. Then, as in the proof of \thmref{invline}, there is a conformal conjugacy $\varphi: \BBB C \smallsetminus K(P_1) \iso \BBB C \smallsetminus K(P_2)$ which extends quasiconformally to the whole plane. Let $z_1=c_1 \mapsto \cdots \mapsto z_p \mapsto z_1$ be the superattracting cycle of $P_1$ which is contained in the cycle $U_1 \mapsto \cdots \mapsto U_p \mapsto U_1$ of Fatou components. By an argument similar to the proof of \thmref{ren}, there exists a quadratic-like restriction $P_1^{\circ p}:U \rightarrow V$ with $U_1 \subset U$ which is hybrid equivalent to $z\mapsto z^2$. Similarly we get a quadratic-like restriction $P_2^{\circ p}:U' \rightarrow V'$ hybrid equivalent to $z\mapsto z^2$. This gives a quasiconformal conjugacy between $P_1$ and $P_2$ on $U_1$, and then on the grand orbit of $U_1$ by taking pull-backs, which extends $\varphi$ to this set. Since the boundary of $\Delta_{P_1}$ is a Jordan curve by \thmref{main}, \lemref{ext} allows us to extend $\varphi$ to a quasiconformal conjugacy on the whole plane. (It is easy to check that $\varphi$ is conformal away from the Julia set $J(P_1)$. But $J(P_1)$ has measure zero by \corref{measure}, hence $\varphi$ is in fact conformal.)

Now let $U$ be a capture component with $c_1$ and $c_2$ being two centers of $U$. As before, there is a conformal conjugacy $\varphi: \BBB C \smallsetminus \overline{GO(\Delta_{P_1})} \iso \BBB C \smallsetminus \overline{GO(\Delta_{P_2})}$ which extends quasiconformally to the entire plane. Again, by \thmref{main} and \lemref{ext}, $\varphi$ can be extended to a quasiconformal conjugacy $\BBB C \rightarrow \BBB C$.\\  

Now we can completely characterize the quasiconformal conjugacy classes in $\PP$.

\begin{thm}[QC Conjugacy Classes in $\PP$]
\label{qcclass}
Every quasiconformal conjugacy class in $\PP$ is one from the following list:
\begin{enumerate}
\item[(a)]
A hyperbolic-like or capture component of the interior of $\MM$ with the center(s) removed.
\item[(b)]
The two components $\Omega_{ext}$ and $\Omega_{int}$.
\item[(c)]
A queer component of the interior of $\MM$.
\item[(d)]
A center of a hyperbolic-like or capture component.
\item[(e)]
A single point on the boundary $\partial \MM$.
\end{enumerate}
\end{thm}
\begin{pf}
\corref{qcopen} shows that no conjugacy class intersects two distinct members of the above list. It also proves that (d) and (e) are in fact conjugacy classes. Also the proof of \thmref{invline} shows that every queer component is a conjugacy class. So it remains to prove that (a) and (b) are quasiconformal conjugacy classes.

Recall that a {\it critical orbit relation} for a rational map $f$ on the sphere
is a coincidence of the form $f^{\circ k}(c_1)=f^{\circ n}(c_2)$, where $c_1$ and $c_2$ are critical points of $f$ and $k$ and $n$ are nonnegative integers 
with $k+n>0$ (we may have $c_1=c_2$). A holomorphic family $\cal F$ of
rational maps has {\it constant critical orbit relations} if every critical orbit relation for $f\in \cal F$ persists under perturbation of $f$ in $\cal F$. Any two rational maps in a holomorphic family with constant critical orbit relations are quasiconformally conjugate (\cite{Mc-Sul}, Theorem 2.7). In other words, critical orbit relations are the only obstruction to constructing quasiconformal conjugacies. 

Now suppose that $U$ is a hyperbolic-like or capture component with the center(s) removed, or $U=\Omega_{ext}$ or $\Omega_{int}$. Then the family $\cal F =\{ P_c \}_{c\in U}$ has no critical orbit relation at all. Therefore, $U$ has to be a quasiconformal conjugacy class. 
\end{pf}
\vspace{0.17in}

\section{Critical Parametrization of Blaschke Products}
\label{sec:blaschke}

This section is the beginning of a digression in the study of cubic Siegel polynomials. We look at a somewhat different class of maps, i.e., certain Blaschke products which will serve as models for the cubics in $\PP$. We will introduce these model maps in Section \ref{sec:blapar} and return to their relation with the cubics in Section \ref{sec:surgery}. 

Let us consider the following space of degree $5$ normalized Blaschke products:
\begin{equation}
\label{eqn:blas}
\hat{ \cal B} =\{ B:z\mapsto \tau z^3 \left ( \frac{z-p}{1-\overline{p}z} \right )
\left ( \frac{z-q}{1-\overline{q}z} \right ) : B(1)=1\ \mbox{and}\ |p|>1, |q|>1 \}, 
\end{equation}
where the rotation factor $\tau \in \BBB T$ is chosen so as to achieve the normalization $B(1)=1$. 
Each $B\in  \hat{ \cal B}$ has superattracting fixed points at $0$ and $\infty$
and four other critical points counted with multiplicity. We are interested in
the open subset $\cal B \subset \hat{ \cal B}$ of those normalized Blaschke 
products of the form (\ref{eqn:blas}) whose four critical points other than $0$ and $\infty$ are of the form
$$c_1,\ c_2,\ \frac{1}{\overline c_1},\ \frac{1}{\overline c_2}$$
with $|c_1|>1, |c_2|>1$. Our goal is to parametrize elements of $\cal B$ by their critical points $c_1$ and $c_2$. The following theorem provides this ``critical parametrization'' for $\cal B$:
\begin{thm}[Critical Parametrization] 
\label{critpar}
Let $c_1$ and $c_2$ be two points outside the closed unit disk in the complex plane. Then there exists a unique normalized Blaschke product $B\in \cal B$ whose critical points are located at $0, \infty , c_1, c_2, \displaystyle{ \frac{1}{\overline c_1}, \frac{1}{\overline c_2} } $.
\end{thm} 

The proof of this theorem will be given after the following two supporting lemmas. We remark 
that we would like to find a direct proof of this fact which can be generalized to higher degrees, but we have not been able to find such a proof. (Compare a similar situation in \cite{Zakeri}, where a conceptual proof is possible.)

The space $\hat{ \cal B}$ of all Blaschke products of the form (\ref{eqn:blas}) can be identified with the set of all unordered pairs $\{ p, q \}$ of points outside 
the closed unit disk. This is canonically homeomorphic to the symmetric product of two copies of the punctured plane. The latter can be identified with the space of all degree $2$ monic polynomials 
$$w\mapsto (w-w_1)(w-w_2)=w^2-(w_1+w_2)w+w_1w_2$$
with $w_1w_2\neq 0$. It follows that $\hat{\cal B}$ is homeomorphic to $\BBB C
\times \BBB C^{\ast}$. In particular, it is an open topological manifold of
real dimension $4$. 

In the same way, we may consider the space $\cal C$ of all unordered pairs
$\{ c_1, c_2\} $ of points outside the closed unit disk, which has a completely similar description. 

We consider the continuous map 
$$\Psi :\cal B \rightarrow \cal C$$
which sends a normalized Blaschke product $B \simeq \{ p,q \}$ with critical points  
 $\{ 0, \infty , c_1, c_2, \\
\displaystyle {\frac{1}{\overline c_1}, \frac{1}{\overline c_2} } \} $ to the unordered pair $\{ c_1, c_2\}$. 
\begin{lem}
\label{psiprop}
$\Psi$ is a proper map.
\end{lem}
\begin{pf}
Let $B_n \simeq \{ p_n, q_n\}$ be a sequence of normalized Blaschke products in $\cal B$ which leaves every compact subset of $\cal B$. Then either
\begin{enumerate}
\item[$\bullet$]
Some critical point of $B_n$ accumulates on the unit circle, or
\item[$\bullet$]
After relabeling, $p_n$ goes to $\infty$, or
\item[$\bullet$]
After relabeling, $p_n$ accumulates on the unit circle.
\end{enumerate}
In the first two cases, it is easy to see that $\Psi (B_n)$ leaves every compact subset of $\cal C$. In the third case, there is a subsequence of $B_n$ which converges locally uniformly to a Blaschke product of degree $<5$. It follows that the corresponding subsequence of $\Psi (B_n)$ has to leave every compact subset of $\cal C$.
\end{pf}
\begin{lem}
\label{psiinj}
$\Psi$ is injective.
\end{lem}
\begin{pf}
Let $A$ and $B$ be two normalized Blaschke products in $\cal B$ with the same critical points $\{ 0, \infty , c_1, c_2, \displaystyle{ \frac{1}{\overline c_1}, \frac{1}{\overline c_2} } \} $ . Let  
$$ A:z\mapsto \tau_A z^3 \left ( \frac{z-p_1}{1-\overline{p_1}z} \right )
\left ( \frac{z-q_1}{1-\overline{q_1}z} \right ), $$
$$ B:z\mapsto \tau_B z^3 \left ( \frac{z-p_2}{1-\overline{p_2}z} \right )
\left ( \frac{z-q_2}{1-\overline{q_2}z} \right ), $$
and assume by way of contradiction that $p_1\neq p_2$ and $p_1\neq q_2$.
Consider the rational function
$$R(z)=\frac{A(z)}{B(z)}.$$
Clearly $\deg R=4$ and hence $R$ has $6$ critical points counted with multiplicity. We have
$$A'(z)=(\mbox{const.})\frac{z^2\ \prod (z-c_j)(1-\overline{c}_jz)}{(1-\overline{p}_1z)^2
(1-\overline{q}_1z)^2}\ \ ,\ \ B'(z)=(\mbox{const.})\frac{z^2\ \prod (z-c_j)(1-\overline{c}_jz)}{(1-\overline{p}_2z)^2 (1-\overline{q}_2z)^2}$$
from which it follows that 
$$R'(z)=(\mbox{const.})\frac{1}{z}\ \prod (z-c_j)(1-\overline{c}_jz)\left \{ \frac{
\sum (-1)^j(z-p_j)(z-q_j)(1-\overline{p}_jz)(1-\overline{q}_jz)}{(z-p_2)^2 (z-q_2)^2(1-\overline{p}_1z)^2(1-\overline{q}_1z)^2}\right \}.$$
(Note that all the sums and products are taken over $j=1,2$.)
From the above expression, $R$ has already $4$ critical points at the $c_j$
and $1/\overline{c_j}$. So the rational function in the braces should have exactly 
$2$ roots. Since this fraction is irreducible (by our assumption $p_1\neq 
p_2$ and $p_1\neq q_2$), the numerator should have degree $2$. But that implies
$$p_1q_1=p_2q_2,$$
$$\overline{p}_1(1+|q_1|^2)+\overline{q}_1(1+|p_1|^2)=\overline{p}_2(1+|q_2|^2)+\overline{q}_2(1+|p_2|^2)$$
from which it follows that $p_1=p_2$ or $p_1=q_2$, hence $q_1=q_2$ or $q_1=p_2$,
which contradicts our assumption.
\end{pf}

\noindent
{\it Proof of the Theorem (Critical Parametrization).} By \lemref{psiprop} and \lemref{psiinj},
$\Psi $ is a covering map of degree $1$. Hence, it is a homeomorphism $\cal B \stackrel{\simeq}{\longrightarrow} \cal C$. $\Box$\\ 

In particular, the theorem shows that $\cal B$ is also homeomorphic to the product $\BBB C \times \BBB C^{\ast}$. 
\begin{cor}
\label{critpar2}
Given any two points $c_1$ and $c_2$ in the plane,
with $|c_1|\geq 1$ and $|c_2|\geq 1$, there exists a unique normalized Blaschke product $B$ in the closure $\overline{\cal B}$ with critical points $\{ 0, \infty , c_1, c_2, \displaystyle{\frac{1}{\overline c_1}, \frac{1}{\overline c_2}\}} $.
\end{cor} 

In other words, critical parametrization is possible even if one or both
critical points $c_1,c_2$ belong to the unit circle.
\begin{pf}
Take a sequence $\{ c_1^n, c_2^n \}$ of pairs of points outside the closed unit disk such that $c_1^n \rightarrow c_1$ and $c_2^n \rightarrow c_2$
as $n\rightarrow \infty$. The corresponding sequence $\Psi ^{-1} (\{ c_1^n, c_2^n \} )$ of normalized Blaschke products has a subsequence which converges to a normalized Blaschke product which, by continuity of $\Psi$, has critical points at $\{ 0, \infty , c_1, c_2, \displaystyle{\frac{1}{\overline c_1}, \frac{1}{\overline c_2} } \} $.

To see uniqueness, it is enough to note that the proof of \lemref{psiinj} can be repeated word by word even if we assume $|c_1|=1$ or $|c_2|=1$.  
\end{pf}

\begin{prop}
\label{realhom}
Every $B\in {\cal B}$ induces a real-analytic diffeomorphism of the unit circle. Consequently, if $B\in \overline{\cal B}\smallsetminus \cal B$, the restriction of $B$ to the unit circle will be a real-analytic homeomorphism with one (or two) critical point(s).
\end{prop}
\begin{pf}
Let us consider $B\in \cal B$ as in (\ref{eqn:blas}) which has critical points at $0, \infty , c_1, c_2, \displaystyle{ \frac{1}{\overline c_1}, \frac{1}{\overline c_2} }$, with
$|c_1|>1$ and $|c_2|>1$ and let us prove that $B|_{\BBB T}$ is a real-analytic 
diffeomorphism. Since $B|_{\BBB T}$ has no critical points, it is a local 
diffeomorphism, hence a covering map of some degree $d\leq 5$. We will prove that $d=1$. 

$B$ induces a branched covering from every connected component $D$ of $B^{-1}(\BBB D)$ to $\BBB D$. Let $D$ be any such component other than the one
whose boundary is $\BBB T$ and contains the origin. Then $\partial D \cap \BBB T=\emptyset $ since otherwise
every point in the intersection would be a critical point of $B$. Since 
$p,q\in B^{-1}(\BBB D)$, either 
\begin{enumerate}
\item[(i)]
There are two components $D_1$ and $D_2$ of $B^{-1}(\BBB D)$ with $p\in D_1$
and $q\in D_2$ such that $B:D_j\rightarrow \BBB D$ is a conformal isomorphism 
for $j=1,2$; or
\item[(ii)]
Both $p$ and $q$ belong to the same component $D$ of $B^{-1}(\BBB D)$ and $B:D\rightarrow \BBB D$ is a 2-to-1 branched covering. 
\end{enumerate}
By the Maximum Principle and the fact that all poles of $B$ are inside $\BBB D$, these components have to be topological disks with piecewise analytic boundaries. It follows that in either case (i) or (ii) the boundaries of the corresponding components give two preimages for $\BBB T$ counted with multiplicity. Since $B^{-1}(\BBB T)$ is symmetric with respect to the unit circle, we have a total number of $4$ preimages for $\BBB T$ other than $\BBB T$ itself. Clearly this means that the degree of $B|_{\BBB T}$ is $1$.

Now let us assume that $B\in \overline{\cal B}\smallsetminus \cal B$. Then there exists a sequence $B_n\in \cal B$ which converges locally uniformly to $B$. Since $B$ has at least one double critical point on $\BBB T$, it follows that $B|_{\BBB T}$ is a real-analytic homeomorphism.
\end{pf}

\vspace{0.17in}

\section{A Blaschke Parameter Space}
\label{sec:blapar} 

Now we focus on a certain class of degree $5$ Blaschke products. These are the maps $B$ with the following two properties:\\
\begin{enumerate}
\item[(i)]
$B$ has the form
\begin{equation}
\label{eqn:blass} 
 B:z\mapsto e^{2 \pi i t} z^3 \left ( \frac{z-p}{1-\overline{p}z} \right )
\left ( \frac{z-q}{1-\overline{q}z} \right ), \ \ \ \ |p|>1, |q|>1 
\end{equation}
where $p$ and $q$ are chosen such that $B$ has a double critical point on the unit circle $\BBB T$ and a pair $(c,1/\overline{c})$ of symmetric critical points which may or may not be on $\BBB T$.
\item[(ii)]
$t$ is the unique number in $[0,1]$ for which the rotation number of $B|_{\BBB T}$ is equal to $\theta$, with $0<\theta<1$ being a given irrational number.
\end{enumerate}

The number $t$ in (ii) is unique because the rotation number of $B$ in (\ref{eqn:blass}) is a continuous and increasing function of $t$ which is strictly increasing at all irrational values (see for example \cite{Katok}, Proposition 11.1.9). 

From the above description, it follows that every $B$ which satisfies (i) and (ii) can be represented as a normalized Blaschke product in $\overline{\cal B}\smallsetminus \cal B$ followed by a unique rotation which adjusts the rotation number to $\theta$. As a consequence, \corref{critpar2} shows that every such $B$ is uniquely determined by the position of its critical points. 

The rotation group {\bf rot}$=\{ R_{\rho}:z\mapsto \rho z\ \ \mbox{with}\ \ |\rho|=1 \}$ acts on the set of all such Blaschke products by conjugation. In fact,
$$R_{\rho}^{-1}\circ B\circ R_{\rho}:z\mapsto e^{2 \pi i t} \rho^4 z^3 \left ( \frac{z-p \overline{\rho}}{1-\overline{p}\rho z} \right )
\left ( \frac{z-q \overline{\rho}}{1-\overline{q}\rho z} \right ). $$

\noindent
We would like to understand the topology of the space of all ``critically marked'' Blaschke products satisfying (i) and (ii) modulo the action of {\bf rot}. By \corref{critpar2}, the conjugacy class of such a Blaschke product is uniquely determined by the location of its critical points up to a rotation. In case there is only one double critical point on $\BBB T$, we can simply represent every conjugacy class by the unique Blaschke product which has a double critical point at $z=1$. Therefore, the quotient space of all critically marked Blaschke products satisfying (i) and (ii) is canonically homeomorphic to the space of all configurations of the two marked critical points $c_1$ and $c_2$ outside the unit disk $\BBB D$ with $c_1=1$ or $c_2=1$. This is just the disjoint union of two copies of $\BBB C\smallsetminus \BBB D$ glued together along the boundary circle by the identification
$$ (1,c_2)\sim (c_1,1) \Longleftrightarrow c_1=\frac{1}{c_2}.$$
It is not hard to see that the resulting space is topologically a punctured plane (see \figref{glue}).

\realfig{glue}{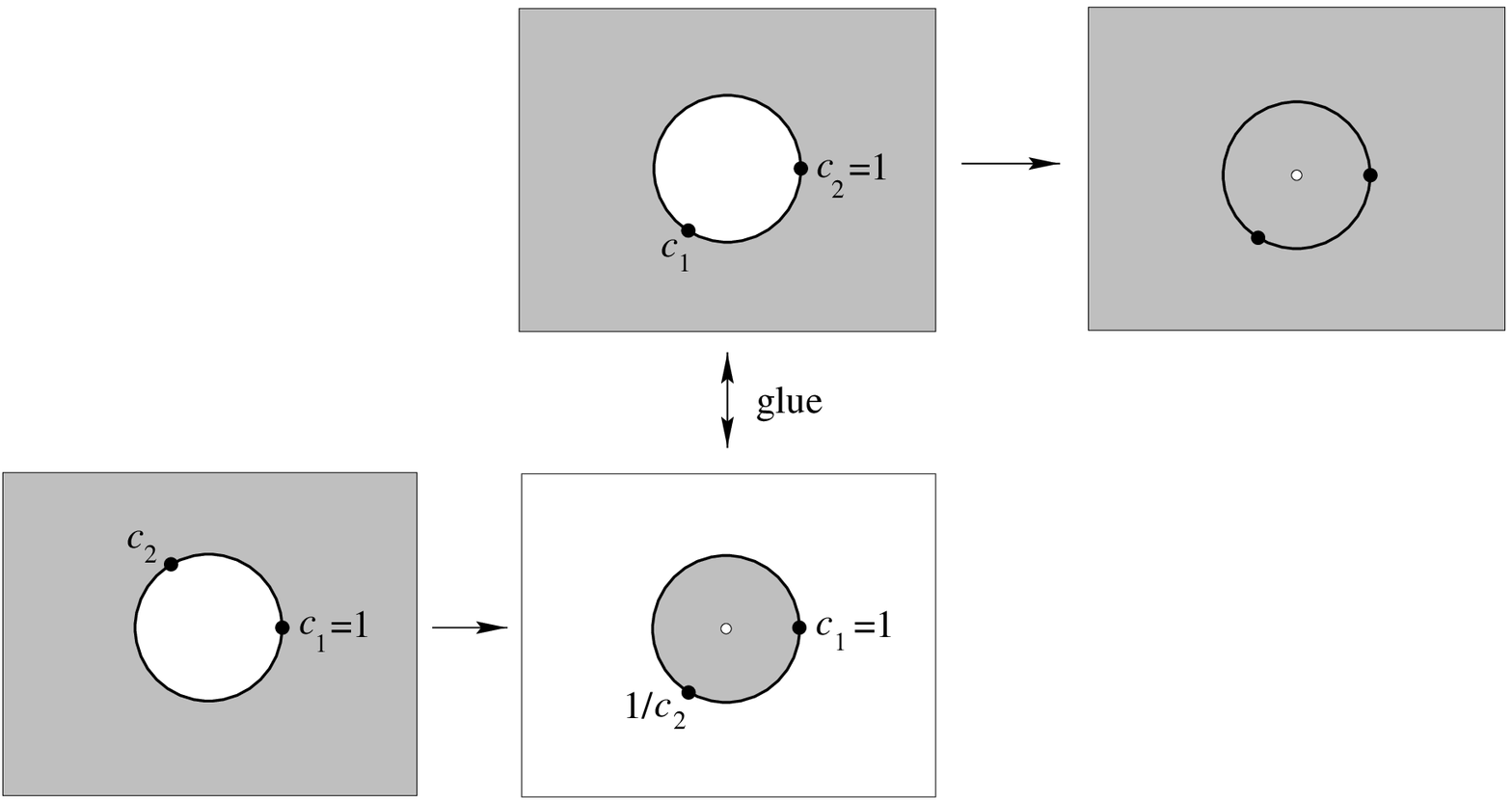}{{\sl Topology of the parameter space $\BB$.}}{13cm}

The space of all critically marked Blaschke products $B$ satisfying (i) and (ii) above modulo the action of {\bf rot} is denoted by $\BB$. The identification $\BB \simeq
\BBB C^{\ast}$ can be explained by introducing the uniformizing parameter $\mu:\BBB C^{\ast}\rightarrow \BB$ as follows: For $\mu \in \BBB C^{\ast}$ with $|\mu|> 1$, the corresponding Blaschke product $B_{\mu}$ has marked critical points at $\{ 0, \infty,
c_1=\mu, 1/{\overline \mu}, c_2=1 \}$. Similarly, if $|\mu |< 1$, $B_{\mu}$ is the unique Blaschke product with marked critical points at $\{ 0, \infty, c_1=1, c_2=1/\mu, {\overline \mu}\}$. Finally, when $|\mu|=1$, $B_{\mu}$ denotes either the unique Blaschke product $B$ with marked critical points at $\{ 0, \infty, c_1=\mu, c_2=1 \}$ or its conjugate $R_{\mu}^{-1}\circ B\circ R_{\mu}$ with marked critical points at $\{ 0, \infty, c_1=1, c_2=1/\mu \}$. Note that $B_{\mu}=B_{1/\mu}$ as maps, if we forget the marking of the critical points.

In the topology of $\BB$, the convergence of a sequence $\{ B_{\mu_n}\}$ to some $B_{\mu}$ has the following meaning: If $B_{\mu}$ has only one double critical point on $\BBB T$ so that $|\mu|\neq 1$, then $B_{\mu_n}\rightarrow B_{\mu}$ simply means $\mu_n\rightarrow \mu$, i.e., uniform convergence on compact subsets of the plane respecting the convergence of the marked critical points. On the other hand, if $B_{\mu}$ has two double critical points on the unit circle so that $|\mu|=1$, then $B_{\mu_n}\rightarrow B_{\mu}$ means that $\{ \mu_n \}$ can only accumulate on $\mu$ or $1/\mu=\overline{\mu}$. In other words, in the topology of local uniform convergence, $\{B_{\mu_n}\} $ can only accumulate on $B_{\mu}$ or its conjugate $R_{\mu}^{-1}\circ B_{\mu}\circ R_{\mu}$.\\

For future reference, we need a somewhat detailed analysis of the structure of the invariant set $\bigcup_{k\geq 0} B^{-k}(\BBB T)$ for a Blaschke product $B\in \BB$. For similar descriptions in a family of degree $3$ Blaschke products, see \cite{Petersen}.\\ \\
{\bf Definition (Skeletons).}\ Let $B\in \BB$. Define $T_0=\BBB T$ and $T_1=\overline{B^{-1}(T_0)\smallsetminus T_0}$. In general, for $k\geq 2$ we define
$T_k$ inductively as $T_k=B^{-1}(T_{k-1})$. We call the closed set $T_k$ the
$k$-{\it skeleton} of $B$. Note that $B$ commutes with the reflection $I:z\mapsto 1/\overline z$. Therefore, every $T_k$ is invariant under $I$.\\ 

Since $B$ is a holomorphic branched covering of the sphere, it is not hard to see that the preimage of every piecewise analytic Jordan curve under $B$ is a finite union of piecewise analytic Jordan curves intersecting one another at finitely many points which are necessarily among the critical points of $B$. Therefore, each $T_k$ decomposes into a finite number of piecewise analytic Jordan curves with this finite intersection property. 

The next proposition tells us what a $k$-skeleton looks like.
\begin{prop}[Structure of the $k$-Skeleton]
\label{skeleton}

\noindent
\begin{enumerate}
\item[(a)]
For $k\geq 1$, the $k$-skeleton $T_k$ is the union of finitely many piecewise analytic Jordan curves $\{ T_k^1, \cdots, T_k^m \}$ which intersect one another at finitely many points and do not cross the unit circle $\BBB T$. None of the $T_k^i$ encloses $\BBB T$. For any $T_k^i$ in this family, the reflected copy $I(T_k^i)$ also belongs to this family.
\item[(b)]
With the notation of (a), let $D_k^i$ denote the bounded component of $\BBB C\smallsetminus T_k^i$ for $k\geq 1$. For $k=0$, $D_0^i$ could mean either $\BBB D$ or $\overline{\BBB C} \smallsetminus \overline{\BBB D}$. Then for $k\geq 1$, $B$ maps $D_k^i$ onto some $D_{k-1}^j$. The mapping is either a conformal isomorphism or a 2-to-1 branched covering. As a result, $B^{\circ k}$ is a proper holomorphic map from  $D_k^i$ onto $\BBB D$ or $\overline{\BBB C} \smallsetminus \overline{\BBB D}$.
\item[(c)]
If $k\geq 1$ and $i\neq j$, we have $D_k^i\cap D_k^j = \emptyset$.
\item[(d)]
For $k>l\geq 1$, either $D_k^i$ and $D_l^j$ are
disjoint or $D_k^i\subset D_l^j$. Conversely, if $D_k^i\subset D_l^j$, we necessarily have $k\geq l$. 
\end{enumerate}
\end{prop}

Every $D_k^i$ is called a $k$-{\it drop} or simply a {\it drop} of $B$. In other words, $k$-drops are the open topological disks bounded by the Jordan curves 
in the decomposition of the $k$-skeleton of $B$. For $k=0$, we have slightly changed the notion of drops. The unit circle $\BBB T$ is the only Jordan curve in the 0-skeleton of $B$, {\it and we agree to call any of the two topological disks $\BBB D$ or $\overline{\BBB C}\smallsetminus \overline{\BBB D}$ a $0$-drop}. The integer $k$ is called the {\it depth} of $D_k^i$.

\begin{pf}
(a) $B^{-1}(\BBB T)$ is the union of $\BBB T$ and 2 or 4 piecewise
analytic Jordan curves which are symmetric with respect to the unit circle and
intersect it at at most one point. (In fact, none of them crosses the unit circle because a point of crossing would be a simple critical point of $B$ on $\BBB T$.) In particular, $T_1$ is the union of these 2 or 4 Jordan curves. It follows that $B^{-1}(\BBB D)$ consists of 1 or 2 open topological disks outside $\BBB D$ together with a subregion of $\BBB D$ which is bounded by $\BBB T$ and 1 or 2 preimages of $\BBB T$ in $\BBB D$ (see \figref{T1}).

As we mentioned earlier, from the fact that $B$ is a holomorphic branched covering of the sphere and by induction on $k$, it follows that $T_k$ is a finite union of piecewise analytic Jordan curves $\{ T_k^1, \cdots, T_k^m \}$ which intersect one another at finitely many points. These points are necessarily {\it precritical} points of $B$. The fact that none of the $T_k^i$ 
crosses the unit circle also follows easily by induction on $k$.
\realfig{T1}{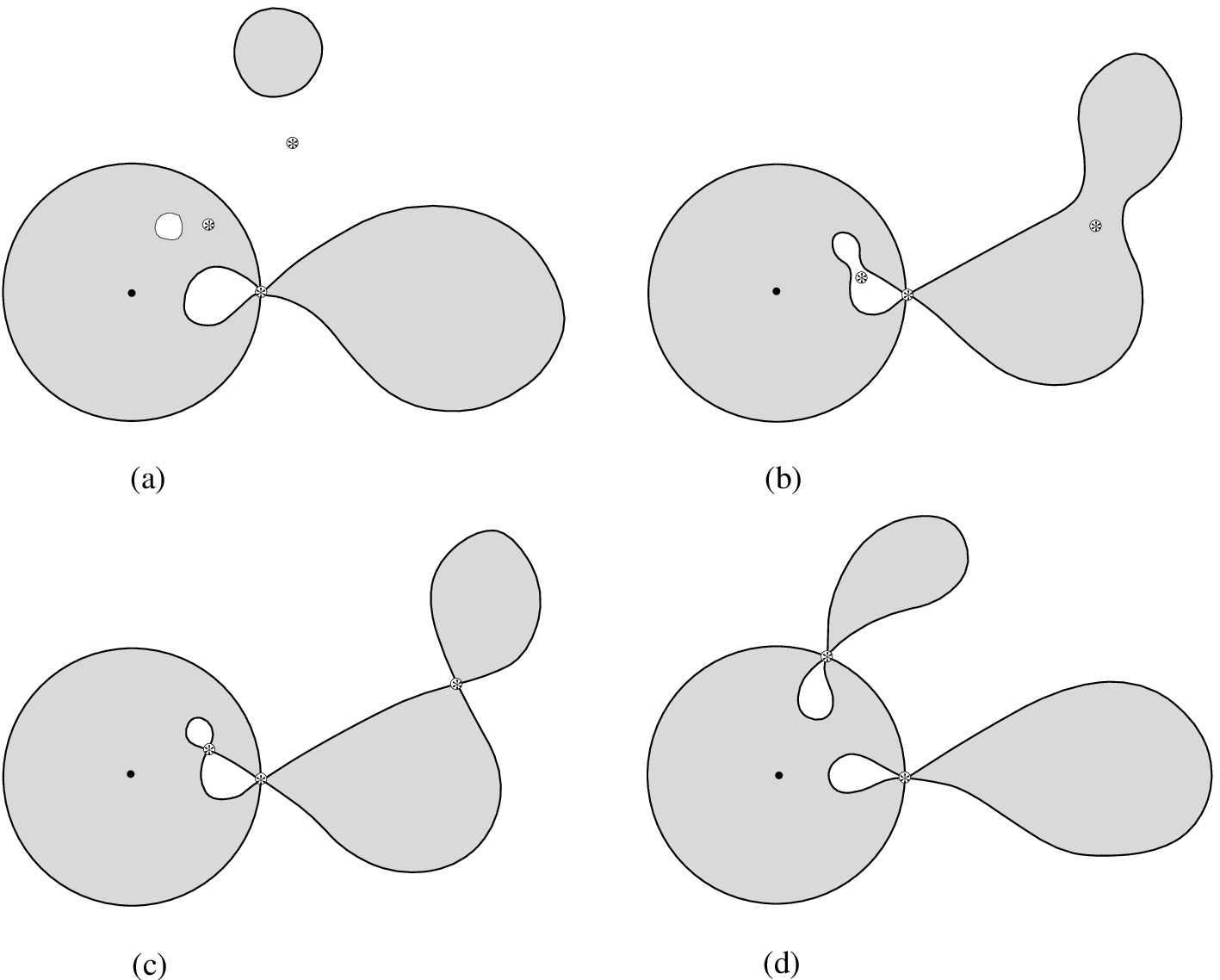}{{\sl Four different configurations for $B^{-1}(\BBB T)$, where $B\in \BB$. The shaded regions are components of $B^{-1}(\BBB D)$. The shaded subregion of $\BBB D$ is mapped to $\BBB D$ by a 3-to-1 branched covering with a superattracting fixed point at the origin. There is a critical point at $z=1$ and the other critical point(s) are symmetric with respect to the unit circle. They are marked by an asterisk. In (a) both components of $B^{-1}(\BBB D)$ outside $\BBB D$ are mapped isomorphically to $\BBB D$. In (b) there is only one component of $B^{-1}(\BBB D)$ outside $\BBB D$ which is mapped onto $\BBB D$ by a 2-to-1 branched covering. (c) is a limiting case of (a) or (b) and (d) is a limiting case of (a).}}{11cm}

(b) By the construction of $T_k$, $B$ maps every $T_k^i$ to some $T_{k-1}^j$. Let $k\geq 1$ and let us assume that $T_k^i$ is completely outside of $\BBB D$. Since all poles of $B$ are inside $\BBB D$, it follows that $B$ is holomorphic in $D_k^i$ and maps it in a proper way onto $D_{k-1}^j$. In case $T_k^i$ is inside $\BBB D$, it follows by symmetry that $B$ maps $D_k^i$ onto some $D_{k-1}^j$ (which is the reflection of the image of $I(D_k^i)$). Since every $D_k^i$ can contain at most one critical point of $B$, in either case the map $B:D_k^i\rightarrow D_{k-1}^j$ will be a conformal isomorphism or a 2-to-1 branched covering.

(c) We prove the claim by induction on $k$. This is obvious for $k=0$. Suppose that there exist two distinct $k$-drops $D_k^i$ and $D_k^i$ which intersect. By (b), $B$ maps both of them to some $(k-1)$-drops and the mapping is proper. It is easy to see that these two $(k-1)$-drops have to be distinct. Then every point in $D_k^i\cap D_k^i$ must map to a point in the intersection of the two $(k-1)$-drops. This contradicts the induction hypothesis.
  
(d) Let $k>l$ and $D_k^i\cap D_l^j \neq \emptyset$. If $D_k^i$ is not contained in $D_l^j$, then $D_k^i\cap T_l^j \neq \emptyset$. Applying $B^{\circ k}$ to $D_k^i$, it follows from (b) that $B^{\circ k}(T_l^j)$ intersects $\BBB D$ or  $\overline{\BBB C} \smallsetminus \overline{\BBB D}$. But $k>l$ implies $B^{\circ k}(T_l^j)=B^{\circ k-l}(B^{\circ l}(T_l^j))=B^{\circ k-l}(\BBB T)=\BBB T$. Conversely, if $D_k^i\subset D_l^j$, then $k$ has to be greater than $l$. This is because $D_k^i$ and $D_l^j$ are intersecting, so by the above argument $k<l$ would imply the reverse inclusion $D_l^j\subset D_k^i$. 
\end{pf}

\noindent
{\bf Definition (Nucleus of a Drop).} Let $D_k^i$ be a drop. We define the {\it nucleus}\footnote{Terminology suggested by A. Epstein.} $N_k^i$ of $D_k^i$ as the set of all points in $D_k^i$ which are not accumulated by any other drop of $B$. The nuclei of $k$-drops are said to have depth $k$.

It follows from \propref{skeleton}(c) that 
$$N_k^i=D_k^i\smallsetminus  \overline{\bigcup_{l\neq k}\bigcup_j D_l^j}.$$
Clearly every nucleus is open. It is also nonempty because every drop contains an open set which eventually maps to the immediate basin of attraction of $0$ or $\infty$, and this open set cannot intersect the closure of any other drop of $B$. 

We have two nuclei of depth zero: $N_0$, which is the nucleus of $\BBB D$ and contains the immediate basin of attraction of 0, and $N_{\infty}$, which is the nucleus of $\overline{\BBB C}\smallsetminus \overline{\BBB D}$ and contains the immediate basin of attraction of $\infty$. Obviously $N_{\infty}=I(N_0)$. It is not hard to see that both $N_0$ and $N_{\infty}$ are invariant under $B$:
\begin{equation}
\label{eqn:inv}
B(N_0)\subset N_0,\ \ \ \ B(N_{\infty})\subset N_{\infty}.
\end{equation}
This of course implies that $N_0$ and $N_{\infty}$ are subsets of the Fatou set of $B$.

It follows from \propref{skeleton}(b) that $B$ maps every nucleus of depth $k$ onto some nucleus of depth $k-1$ and the mapping is either a conformal isomorphism or a 2-to-1 branched covering. We include the following lemma for completeness:

\begin{lem}
\label{chain}
Let $N_k^i$ be the nucleus of a drop $D_k^i$ which eventually maps to the unit disk $\BBB D$. Then 
\begin{enumerate}
\item[(a)]
No point in the orbit
$$N_k^i=N_k^{i_0}\stackrel{B}{\longrightarrow}N_{k-1}^{i_1}\stackrel{B}{\longrightarrow}\cdots \stackrel{B}{\longrightarrow}N_1^{i_{k-1}}\stackrel{B}{\longrightarrow}N_0 $$
can intersect any of the reflected nuclei $I(N_{k-j}^{i_j}),\ 0\leq j\leq k$.
\item[(b)]
For $z\in N_k^i$, $B^{\circ k}$ is the first iterate of $B$ which sends $z$ to $N_0$.
\end{enumerate}
\end{lem}

\begin{pf}
(a) $B$ commutes with $I$, so there is a reflected orbit
$$I(N_k^i)=I(N_k^{i_0})\stackrel{B}{\longrightarrow}I(N_{k-1}^{i_1})\stackrel{B}{\longrightarrow}\cdots \stackrel{B}{\longrightarrow}I(N_1^{i_{k-1}})\stackrel{B}{\longrightarrow}N_{\infty}. $$
Now any point in both orbits would have to map to a point in $N_0$ and $N_{\infty}$ simultaneously, which is impossible since $N_0 \cap N_{\infty}=\emptyset$.

(b) This is obvious if $k=1$. Suppose that $k>1$ and that for some $0<l<k$, $B^{\circ l}(z)\in N_0$. Then by (\ref{eqn:inv}), $B^{\circ k-1}(z)\in N_0\subset \BBB D$. But $B^{\circ k-1}(z) \in B^{\circ k-1}(D_k^i)$ and $B^{\circ k-1}(D_k^i)$ is a 1-drop which does not intersect $\BBB D$.
\end{pf}

\noindent
{\bf Remark.} If $z\in N_k^i$, it is {\it not} true that $B^{\circ k}$ is the first iterate of $B$ which sends $z$ to the unit disk. In fact, the orbit of $z$ can pass through $\BBB D$ several times before it maps to $N_0$ (see \figref{several}).
\realfig{several}{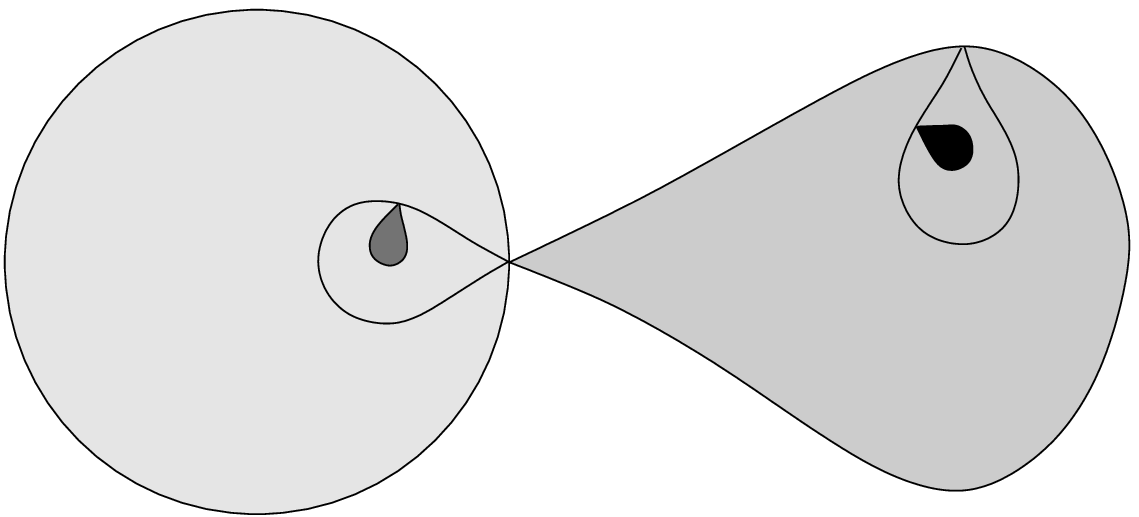}{{\sl The orbit of a 3-drop under the iteration of a Blaschke product $B \in \BB$. This dark drop on the right maps successively to lighter drops. It visits the unit disk once before it maps onto it.}}{5cm}
\begin{prop}
\label{nuc}

\noindent
\begin{enumerate}
\item[(a)] 
Distinct nuclei are disjoint.
\item[(b)]
The map $B^{\circ k}$ from $N_k^i$ onto $N_0$ or $N_{\infty}$ is
either a conformal isomorphism or a 2-to-1 branched covering.
\end{enumerate}
\end{prop}

\begin{pf}
(a) Let $N_k^i$ and $N_l^j$ be two distinct nuclei which intersect. By \propref{skeleton}(c), we have $k\neq l$. Without loss of generality, we assume that $k>l$ and the iterate $B^{\circ l}$ maps $N_l^j$ onto $N_0$. So
for every $z$ in the intersection $N_k^i \cap N_l^j$, $B^{\circ l}(z)$ will belong to $N_0$. This contradicts \lemref{chain}(b).

(b) Since by (a) distinct nuclei are disjoint, an orbit
$$N_k^i=N_k^{i_0}\stackrel{B}{\longrightarrow}N_{k-1}^{i_1}\stackrel{B}{\longrightarrow}
\cdots \stackrel{B}{\longrightarrow}N_1^{i_{k-1}}\stackrel{B}{\longrightarrow}N_0\ \ \mbox{or}\ \ N_{\infty} $$
can hit every critical point of $B$ at most once. Since the critical point $z=1$ of $B$ does not belong to any nucleus, the above orbit can only hit the pair of critical points $c$ and $1/\overline c$, with $|c|\neq 1$. By \lemref{chain}(a), both critical points cannot belong to the 
above orbit simultaneously. This means that $B^{\circ k}:N_k^i
\rightarrow N_0$ or $N_{\infty}$ is either a conformal isomorphism or a 2-to-1 branched covering.
\end{pf}

\figref{hypbla}, \figref{nucleus}, and \figref{blaext} show the Julia sets of some Blaschke products in $\BB$ for $\theta=(\sqrt{5}-1)/2$.
 
\goodbreak

\section{The Surgery}
\label{sec:surgery}

From now on, unless otherwise stated, {\it we assume that $\theta$ is an irrational number of bounded type.} We describe a surgery on degree $5$ Blaschke products in $\BB$ to obtain cubic polynomials in $\PP$. A similar surgery was done previously in the case of quadratic polynomials \cite{Douady2}
using the following theorem of Swiatek and Herman (see \cite{Swiatek} or \cite{Herman2}). Recall that a homeomorphism $h:\BBB R \rightarrow \BBB R$ is called $k$-{\it quasisymmetric} if 
$$0< k^{-1} \leq \frac{|h(x+t)-h(x)|}{|h(x)-h(x-t)|}\leq k < +\infty$$
for all $x$ and all $t>0$. We call $h$ quasisymmetric if it is $k$-quasisymmetric for some $k$. A homeomorphism $h:\BBB T \rightarrow \BBB T$ is $k$-quasisymmetric if its lift to $\BBB R$ is such a homeomorphism.
\realfig{hypbla}{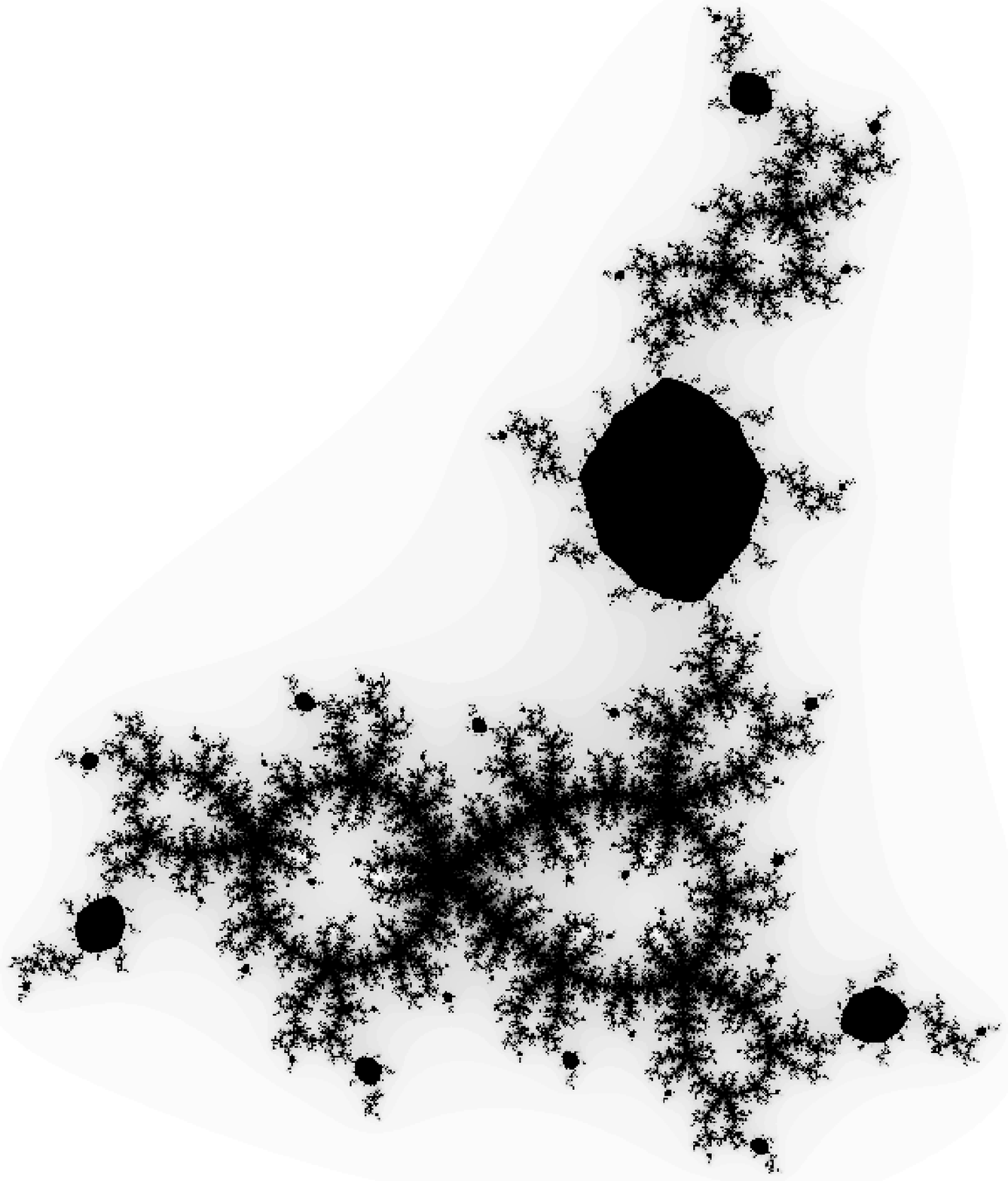}{{\sl The Julia set of a Blaschke product in $\BB$ for $\theta=(\sqrt{5}-1)/2$. There are two symmetric attracting cycles in the nuclei $N_0$ and $N_{\infty}$. The topological disks in black form the basin of attraction of these two cycles. After surgery this Blaschke product becomes a hyperbolic-like cubic in $\PP$.}}{8cm}

\realfig{nucleus}{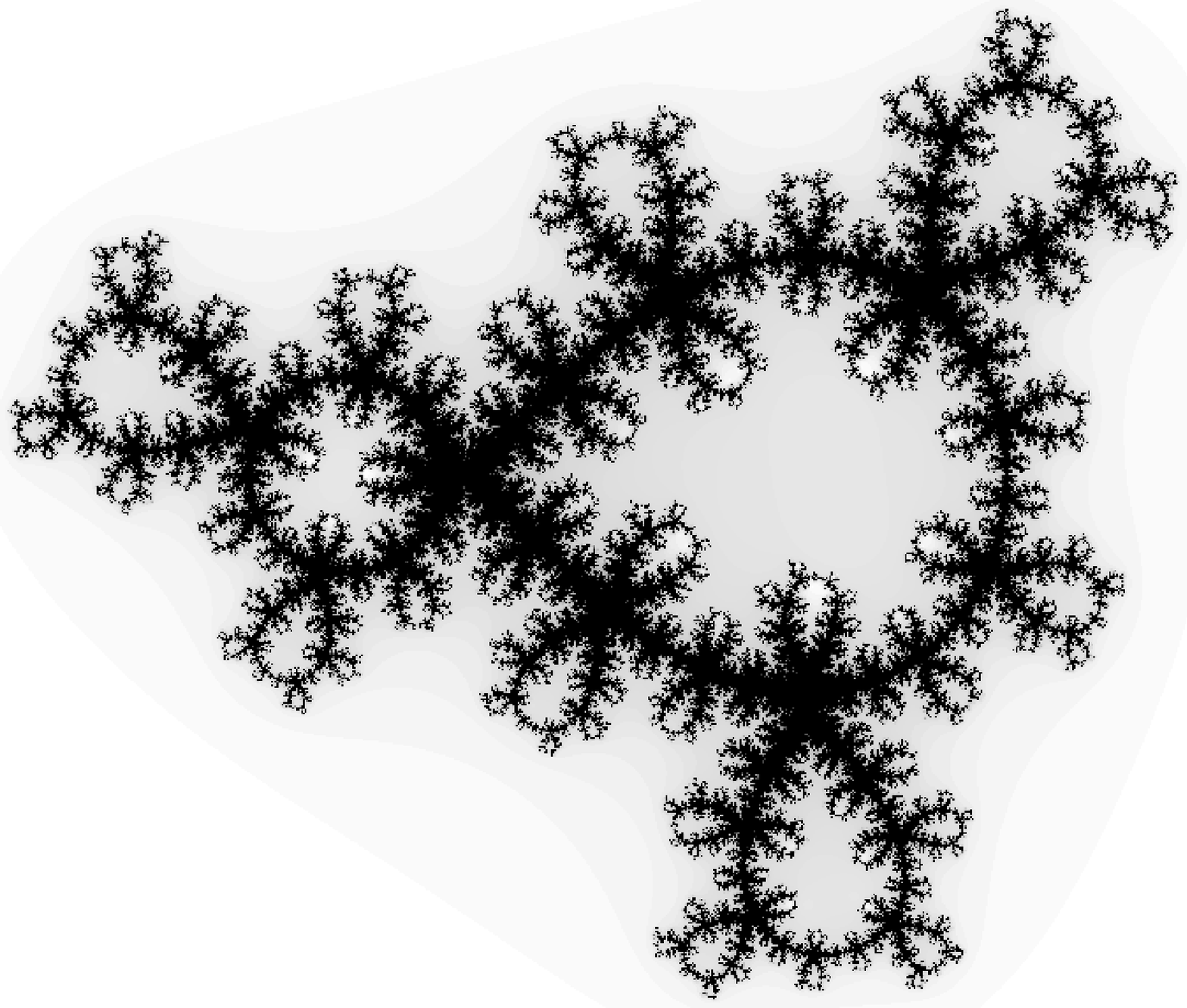}{{\sl Another example of the Julia set of a Blaschke product in $\BB$ for $\theta=(\sqrt{5}-1)/2$. There is a critical point in the nucleus of the large $1$-drop attached to the unit disk at $z=1$ which maps into $N_0$. Hence this nucleus contains the zeros $p$ and $q$. However, after surgery this Blaschke product becomes a capture cubic in $\PP$.}}{8cm}

\begin{thm}[Linearization of Critical Circle Maps]
\label{ccmap}
Let $f:\BBB T \rightarrow \BBB T$ be a real-analytic homeomorphism with finitely many critical points and rotation number $\theta$. Then there exists a quasisymmetric homeomorphism 
$h:\BBB T \rightarrow \BBB T$ which conjugates $f$ to the rigid rotation $R_{\theta}:t\mapsto t+\theta $ (mod 1) if and only if $\theta$ is an irrational number of bounded type. Moreover, if $f$ belongs to a compact family of real-analytic homeomorphisms with rotation number $\theta$, then $h$ is $k$-quasisymmetric, where the constant $k$ only depends on the family and not on the choice of $f$.
\end{thm}

Let us briefly sketch what this surgery does on a Blaschke product $B \in \BB$.
By \propref{realhom}, the restriction $B|_{\BBB T}$ is a real-analytic  homeomorphism with one (or two) critical point(s). When the rotation number of this circle map is of bounded type, by \thmref{ccmap} one can find a unique $k$-quasisymmetric homeomorphism $h:\BBB T\rightarrow \BBB T$ with $h(1)=1$ such that the following diagram commutes:
$$\begin{array}{rlcrl}
  & \BBB T & \stackrel{B}{\longrightarrow} & \BBB T & \\
h & \downarrow &  & \downarrow & h\\
  & \BBB T & \stackrel{R_{\theta}}{\longrightarrow}  & \BBB T &
\end{array}$$
Moreover, the family $\{ B|_{\BBB T} \}_{B \in \BB}$ is compact (compare \thmref{q=3}), hence $h$ is in fact $k(\theta)$-quasisymmetric, where the constant $k(\theta)$ only depends on $\theta$.
We can extend $h$ to a $K(\theta)$-quasiconformal homeomorphism $H:\BBB D \rightarrow \BBB D$ whose dilatation depends only on $\theta$. Possible extensions are given by the theorem of Beurling and Ahlfors \cite{Ahlfors} or Douady and Earle \cite{Douady-Earle} (which has the advantage of being conformally invariant). Define a {\it modified Blaschke product} $\tilde B$ as follows:
\begin{equation}
\label{eqn:modbla}
\tilde{B}(z)=  \left \{ \begin{array}{ll} 
B(z) &  |z|\geq 1 \\
(H^{-1}\circ R_{\theta}\circ H)(z) & |z|<1
\end{array}
\right.
\end{equation}
\realfig{blaext}{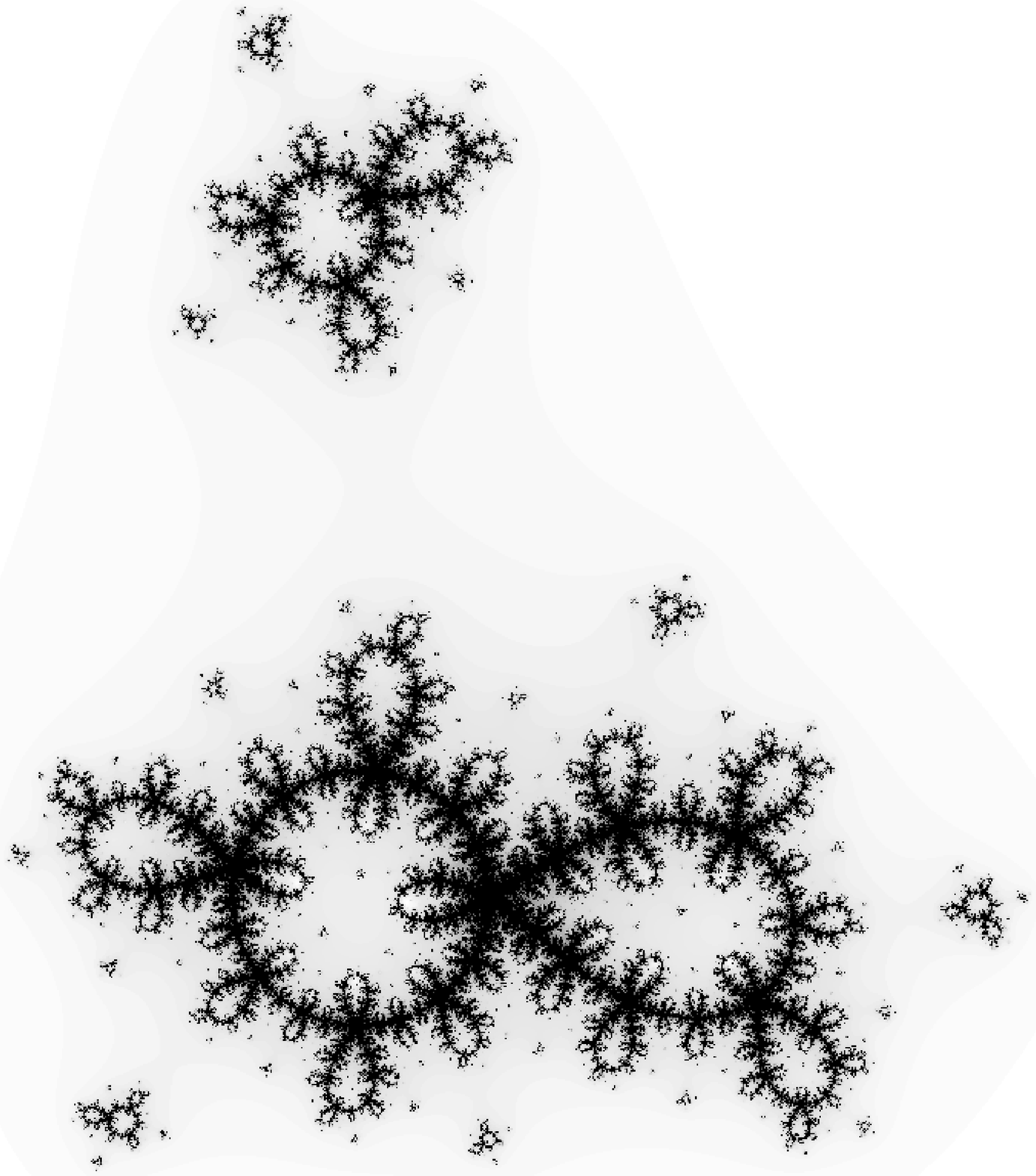}{{\sl The Julia set of a Blaschke product in $\BB$ for $\theta=(\sqrt{5}-1)/2$ outside the connectedness locus $\CC$ (see Section \ref{sec:newcon}). Surgery makes this Blaschke product into a cubic in $\Omega_{ext}$.}}{8cm}

\noindent
This amounts to cutting out the unit disk and gluing in a Siegel disk instead.
Note that the two definitions match along $\BBB T$ by the above commutative diagram. Now define a conformal structure $\sigma$ on the plane as follows: On $\BBB D$, let $\sigma$ be the pull-back $H^{\ast}\sigma_0$ of the standard conformal structure $\sigma_0$. Since $R_\theta$ preserves $\sigma_0$, $\tilde B$ will preserve $\sigma$ on $\BBB D$. For every $k\geq 1$, pull $\sigma|_{\BBB D}$ back by $\tilde {B}^{\circ k}=  
B^{\circ k}$ on $B^{-k}(\BBB D)\smallsetminus \BBB D$ (which consists of all the {\it maximal} $k$-drops of $B$; see Section \ref{sec:newcon}). Since $B^{\circ k}$ is 
holomorphic, this does not increase the dilatation of $\sigma$. Finally, let 
$\sigma=\sigma_0$ on the rest of the plane. By the construction, $\sigma$
has bounded dilatation and is invariant under 
$\tilde B$. Therefore, by the Measurable Riemann Mapping Theorem of 
Ahlfors and Bers, we can find a quasiconformal homeomorphism $\varphi:\BBB C\rightarrow \BBB C$ such that $\varphi^{\ast}\sigma_0=\sigma$. Set
\begin{equation}
\label{eqn:pmodb} 
P=\varphi \circ \tilde{B} \circ \varphi^{-1}.
\end{equation}
Then $P$ is a quasiregular self-map of the sphere which preserves $\sigma_0$, hence it is holomorphic. Also $P$ is proper of degree $3$ since $\tilde B$ has the same properties. Therefore
$P$ is a cubic polynomial.

Now the action of $P$ on $\varphi (\BBB D)$ is quasiconformally conjugate to a rigid rotation, hence $\varphi (\BBB D)$ is contained in a Siegel disk for $P$ with rotation number $\theta$. Since $\varphi (1)$ is a critical point for $P$, it follows that the entire orbit $\{ P ^{\circ k}(\varphi (1))\}_{k\geq 1}$ lives on the boundary of this Siegel disk. But $\{ P ^{\circ k}(\varphi (1))\}_{k\geq 1}$ is dense on $\varphi (\BBB T)$, so $\varphi (\BBB T )$ is exactly the boundary of this Siegel disk, which is a quasicircle passing through the critical point $\varphi (1)$ of $P$.

To mark the critical points of $P$, hence getting an element of $\PP$, we must normalize $\varphi$ carefully. Recall from  Section \ref{sec:blapar} that $\BB$ is uniformized by the parameter $\mu\in \BBB C^{\ast}$ as follows: If $|\mu| \geq 1$, $B_{\mu}$ has marked critical points at $\{ 0, \infty, c_1=\mu, 1/\overline {\mu}, c_2=1 \}$, while for $|\mu| \leq 1$, $B_{\mu}$ has marked critical points at $\{ 0, \infty, c_1=1, c_2=1/\mu, \overline {\mu} \}$. In the first case, we normalize $\varphi$ such that $\varphi(H^{-1}(0))=0$ and $\varphi(1)=1$. Call $\varphi(\mu)=c$ and mark the critical points of $P$ by declaring $P=P_c$ as in  Section \ref{sec:cubpar}. In the case $|\mu|\leq 1$, we normalize $\varphi$ similarly by putting $\varphi(H^{-1}(0))=0$ and $\varphi(1/\mu)=1$, but this time we call $\varphi(1)=c$ and set $P=P_c$. It is easy to see that when $|\mu|=1$, both normalizations produce the same critically marked cubic polynomial in $\PP$.

Let us denote the polynomial $P$ constructed this way by ${\cal S}_H(B)$, i.e., the cubic obtained by performing surgery 
on a Blaschke product $B$ using a quasiconformal extension $H$. The first question we would like to address is the following:
\begin{enumerate}
\item[]
``Given a $B\in \BB$, what cubic polynomials of the form ${\cal S}_H(B)$ can we obtain as the result of this surgery by choosing different quasiconformal extensions $H$?"
\end{enumerate}

We will see that for two quasiconformal extensions $H$ and $H'$, the cubics 
${\cal S}_H(B)$ and ${\cal S}_{H'}(B)$ are quasiconformally conjugate and the conjugacy is conformal everywhere except on the grand orbit of the Siegel disk centered at the origin. When ${\cal S}_H(B)$ is capture, we can certainly end up with two different cubics if we choose the extensions arbitrarily. In fact, let $k$ be the first moment the orbit of the critical point $c$ of $B$ hits the unit disk, and let $w=B^{\circ k}(c)$.
Then for two quasiconformal extensions $H$ and $H'$, the captured images of the critical points of ${\cal S}_H(B)$ and ${\cal S}_{H'}(B)$ have the same conformal position in their corresponding Siegel disks if and only if $H(w)=H'(w)$. It follows that ${\cal S}_H(B)\neq {\cal S}_{H'}(B)$ as soon as we choose two different extensions $H,H'$ with $H(w)\neq H'(w)$. 

The following proposition has a very nontrivial content in case the result of the surgery is a cubic whose Julia set has positive measure (say, in a queer component). It is the Bers Sewing Lemma which makes the proof work.  
\begin{prop}
\label{independent}
Let $P={\cal S}_H(B)$ and $H'$ be any other quasiconformal extension of the circle homeomorphism $h$ which linearizes $B|_{\BBB T}$. Then, if $P$ is not capture, ${\cal S}_H(B)={\cal S}_{H'}(B)$. On the other hand, when $P$ is capture, ${\cal S}_H(B)={\cal S}_{H'}(B)$ if and only if $H(w)=H'(w)$, where $w\in \BBB D$ is the captured image of the critical point of $B$.   
\end{prop} 

\begin{pf}
Let $Q={\cal S}_{H'}(B)$ and $\varphi_H$ and $\varphi_{H'}$ denote the quasiconformal homeomorphisms which satisfy $P=\varphi_H \circ \tilde{B}_H \circ \varphi_H^{-1}$ and $Q=\varphi_{H'} \circ \tilde{B}_{H'} \circ \varphi_{H'}^{-1}$ as in (\ref{eqn:pmodb}). The homeomorphism $\varphi$ defined by
$$\varphi(z)= \left \{ 
\begin{array}{ll}
(\varphi_{H'}\circ \varphi_H^{-1})(z) & z\in \BBB C\smallsetminus GO(\Delta_P)\\
(\varphi_{H'}\circ B^{-k}\circ {H'}^{-1}\circ H \circ B^{\circ k}\circ \varphi_H^{-1})(z) & z\in P^{-k}(\Delta_P)
\end{array} \right. $$
is quasiconformal and conjugates $P$ to $Q$. By \lemref{ext}, one can find a quasiconformal conjugacy $\psi:\BBB C \rightarrow \BBB C$ between $P$ and $Q$ which is conformal on the grand orbit $GO(\Delta_P)$ and agrees with $\varphi$ everywhere else. By the Bers Sewing Lemma, $\overline{\partial} \psi=\overline{\partial}\varphi$ almost everywhere on $\BBB C\smallsetminus GO(\Delta_P)$. But the latter generalized partial derivative vanishes almost everywhere on $\BBB C\smallsetminus GO(\Delta_P)$ because the surgery does not change the conformal structures outside $\bigcup_{k\geq 0}B^{-k}(\BBB D)$. Hence $\overline{\partial} \psi=0$ almost everywhere on $\BBB C$,
which means $\psi$ is conformal. This shows $P=Q$. 
\end{pf}
\vspace{0.15in}
\noindent
{\bf Convention.} For the rest of this paper, we always choose the Douady-Earle extension of circle homeomorphisms to perform surgery. By the above proposition, this is really a ``choice'' only in the capture case. We can therefore neglect the dependence on $H$ and call 
$${\cal S}:\BB \rightarrow \PP$$
the {\it surgery map}.\\

As an immediate corollary of normalization of $\varphi$ and construction of $\cal S$, we have the following:

\begin{cor}
\label{switch}
Let $\mu \in \BBB C^{\ast}$ and $P_c={\cal S}(B_{\mu})$ be the cubic obtained 
by performing the above surgery.
\begin{enumerate}
\item[$\bullet$]
If $|\mu|>1$, then $1 \in \partial \Delta_c$ and $c \notin \partial \Delta_c$.
\item[$\bullet$]
If $|\mu|<1$, then $c \in \partial \Delta_c$ and $1 \notin \partial \Delta_c$.
\item[$\bullet$]
If $|\mu|=1$, then both $c$ and $1 \in \partial \Delta_c$.\ $\Box$
\end{enumerate}
\end{cor}

\vspace{0.17in}

\section{The Blaschke Connectedness Locus ${\cal C}_5(\theta)$}
\label{sec:newcon}

Suggested by the case of cubic polynomials, we define the {\it Blaschke Connectedness Locus} $\CC$ by 
$$\CC= \{ B\in \BB:\ \mbox{The Julia set $J(B)$ is connected} \}.$$
The following theorem provides a useful characterization of $\CC$ in terms of the critical orbits. 
\begin{thm}
\label{hitD}
$B\in \CC$ if and only if one of the following holds:
\begin{enumerate}
\item[$\bullet$]
The orbit of $c$, the critical point of $B$ in $\BBB C \smallsetminus \BBB D$ other than  $1$, eventually hits $\overline{\BBB D}$.
\item[$\bullet$]
The orbit of $c$ never hits $\overline{ \BBB D}$, but remains bounded.
\end{enumerate}
\end{thm}

The proof of this theorem depends on an alternative dynamical description for the Julia set of a Blaschke product in $\BB$ which is obtained by taking pull-backs
along certain type of drops called the maximal drops. This description will be useful later in the proof of \thmref{ABconj}.\\ \\
{\bf Definition.} Let $D_k^i$ be a $k$-drop of $B\in \BB$. 
We call $D_k^i$ a {\it maximal drop} if $D_k^i=\BBB D$, or if $D_k^i\cap \BBB D=\emptyset$ and $D_k^i$ is not contained in any other $l$-drop of $B$ for $l\geq 1$.\\

It follows in particular that maximal drops of $B$ are disjoint.
\begin{prop}
\label{aux} 
Let $B\in \BB$ and let $P=\cal S (B)=\varphi \circ \tilde{B} \circ \varphi^{-1}$ as in $($\ref{eqn:pmodb}$)$. Then
\begin{enumerate}
\item[(a)]
$D_k^i$ is a maximal drop of $B$ if and only if $\varphi(D_k^i)$ is a Fatou component of $P$ which eventually maps to the Siegel disk $\Delta_P$.
\item[(b)]
$\varphi$ maps the nucleus $N_{\infty}$ of $B$ onto $\overline{\BBB C}\smallsetminus \overline{GO(\Delta_P)}$.
\item[(c)]
The boundary of the immediate basin of attraction of infinity for $B$ is precisely the closure of the union of the boundaries of all the maximal drops of $B$.
Under $\varphi$ this set maps to the Julia set $J(P)$.
\end{enumerate}
\end{prop}

\begin{pf}
(a) and (b) are easy consequences of the definitions. For (c), just note that under $\varphi$, the boundary of the immediate basin of attraction of infinity for $B$ corresponds to the similar boundary for $P$, and the closure of the union of the boundaries of all the maximal drops of $B$ corresponds to the Julia set $J(P)$ by (a). 
\end{pf}

\begin{lem}[Alternative description for Julia Sets]
\label{alt}
Let $B\in \BB$ and let $J_0$ be the boundary of the immediate basin of attraction of infinity for $B$. Define a sequence of compact sets $J_n=J_n(B)$ inductively by
\begin{equation}
\label{eqn:jn}
J_n= \overline{ \bigcup_{D_k^i\ maximal} B^{-k}(IJ_{n-1}\cap \BBB D)\cap D_k^i },
\end{equation}
Then
\begin{equation}
\label{eqn:jb}
J(B)=\overline{\bigcup_{n\geq 0}J_n}.
\end{equation}
\end{lem}
\begin{pf}
Each $J_n$ is compact and contained in $J(B)$. By \lemref{aux}(c), $J_0\subset J_1$ and it follows by induction on $n$ that $J_n\subset J_{n+1}$ for $n\geq 0$. Put
$$J_{\infty}=\overline{\bigcup_{n\geq 0}J_n}.$$
Clearly $J_{\infty}$ is compact and contained in the Julia set $J(B)$, and it is not hard to see that it is invariant under the reflection $I$. We will show that $J_{\infty}$ is totally invariant under $B$, i.e., $B^{-1}(J_{\infty})=J_{\infty}$. This will prove that $J_{\infty}=J(B)$.

First we prove that $J_{\infty}$ is forward invariant. For any $n$, it follows from (\ref{eqn:jn}) that $B(J_n \smallsetminus \BBB D)\subset J_n \subset J_{\infty}$. On the other hand, $B(J_n\cap \overline{\BBB D})= B(IJ_{n-1}\cap \overline{\BBB D})=
IB(J_{n-1}\smallsetminus \BBB D)\subset IJ_{\infty}=J_{\infty}$. These two inclusions show that $B(J_n)\subset J_{\infty}$, hence $B(J_{\infty})\subset J_{\infty}$. 

To prove backward invariance, first note that for any $n$, $B^{-1}(J_n)\smallsetminus \BBB D \subset J_n \subset J_{\infty}$ by (\ref{eqn:jn}). To obtain the same kind of inclusion for $B^{-1}(J_n)\cap \overline{\BBB D}$, we distinguish two cases: First, $B^{-1}(J_n\cap \overline{\BBB D})\cap \overline{\BBB D}= B^{-1}(IJ_{n-1}\cap \overline{\BBB D}) \cap \overline{\BBB D}\subset I(B^{-1}(J_{n-1}\smallsetminus \BBB D)) \subset IJ_{n-1}\cup J_n \subset J_{\infty}$. Second, $B^{-1}(J_n\smallsetminus \BBB D)\cap \overline{\BBB D}=I(B^{-1}(IJ_n\cap \overline{\BBB D})\smallsetminus \BBB D) \subset I(B^{-1}(J_{n+1})\smallsetminus \BBB D)\subset IJ_{n+1}\subset J_{\infty}$. Altogether, these three inclusions show that $B^{-1}(J_n)\subset J_{\infty}$ for all $n$. Hence $B^{-1}(J_{\infty})\subset J_{\infty}$ and this proves (\ref{eqn:jb}).
\end{pf}
\vspace{0.5cm}
\noindent
{\bf Remark.} In terms of the modified Blaschke product $\tilde{B}$ as defined in (\ref{eqn:modbla}), one can also define the sequence $J_n$ by $J_0=\varphi^{-1}(J(P))$ and
$$J_n=\overline{ \bigcup_{k\geq 0} \tilde{B}^{-k}(IJ_{n-1}\cap \BBB D) }.$$
Here $\tilde{B}^{-k}$ refers to any branch of $(\tilde{B}^{\circ k})^{-1}$ of the form $\tilde{B}^{-1} \circ \ldots \circ \tilde{B}^{-1}$ ($k$ times) where each branch of $\tilde{B}^{-1}$ satisfies $\tilde{B}^{-1}(\BBB D)\cap \BBB D =\emptyset$.\\ \\  
{\it Proof of \thmref{hitD}.} One direction is quite easy to see: If the orbit of $c$ never hits the closed unit disk and escapes to infinity, one can easily show that $J(B)$ is disconnected exactly like the polynomial case by considering the B\"{o}ttcher map of the immediate basin of attraction of $\infty$ for $B$ (see for example \cite{Milnor1}, Theorem 17.3).

Conversely, suppose that the orbit of the critical point $c$ either hits $\overline {\BBB D}$ or stays bounded in $\BBB C \smallsetminus 
\overline {\BBB D}$. Then the Julia set $J(P)$ is connected, where $P={\cal S}(B)$. Consider the sequence of compact sets $J_n$ in (\ref{eqn:jn}). By \propref{aux}(c), $J_0$ is connected and it follows by induction on $n$ that each $J_n$ defined by (\ref{eqn:jn}) is connected. Therefore (\ref{eqn:jb}) shows that $J(B)$ is connected. Hence $B\in \CC$. $\Box$\\
 
In what follows, we prove that the connectedness locus $\CC$ is compact. Other facts, e.g.,  having only two complementary components, or connectivity, will be proved later using surgery (see \corref{c5con} and \corref{c5full}). We would like to remark that unlike the case of cubic polynomials, it is often difficult to prove anything about the topology of the Blaschke connectedness locus, partly because of the complicated way these Blaschke products depend on their critical points, but more importantly because of the fact that the family
$\mu\mapsto B_{\mu}$ does not depend holomorphically on $\mu$. 

\begin{lem}
\label{hn}
Let $\{ B_n \}$ be an arbitrary sequence of Blaschke products in $\BB$ and $h_n:\BBB T \rightarrow \BBB T$ be the unique normalized quasisymmetric homeomorphism which conjugates $B_n|_{\BBB T}$ to the rigid rotation $R_{\theta}$. Let $H_n$ denote the Douady-Earle extension of $h_n$. Then the sequence $\{ H_n \}$ has a subsequence which converges locally uniformly to a quasiconformal homeomorphism of $\BBB D$.
\end{lem}

It follows that the sequence $\{ H_n^{-1}(0) \}$ stays in a compact subset of the unit disk.

\begin{pf}
Regarding $\BBB T$ as the quotient $\BBB R / \BBB Z$, we can lift each $h_n$ to a $k(\theta)$-quasisymmetric homeomorphism $\tilde{h}_n: \BBB R \rightarrow \BBB R$ which fixes $0$ and satisfies $\tilde{h}_n(x+1)=\tilde{h}_n(x)+1$ for all $x$. The space of all uniformly quasisymmetric normalized homeomorphisms of the real line is compact (\cite{Lehto}, Lemma 5.1), hence a subsequence of $\{ \tilde{h}_n \}$ converges uniformly to a $k(\theta)$-quasisymmetric homeomorphism $\tilde{h}:\BBB R \rightarrow \BBB R$. This homeomorphism descends to a quasisymmetric homeomorphism $h:\BBB T \rightarrow \BBB T$, which is the uniform limit of the corresponding subsequence of $\{ h_n \}$. On the other hand, the Douady-Earle extension depends continuously on the circle homeomorphism \cite{Douady-Earle}. It follows that the corresponding subsequence of $\{ H_n \}$ converges locally uniformly on $\BBB D$ to the extension of $h$.
\end{pf}

\begin{cor}
\label{normal}
Let $B \in \BB$ and $\varphi_B: \BBB C \rightarrow \BBB C$ be the quasiconformal homeomorphism which conjugates the modified Blaschke product $\tilde B$ to the cubic $P={\cal S}(B)$ as in (\ref{eqn:pmodb}): $P=\varphi_B \circ \tilde{B} \circ \varphi_B^{-1}$. Then the family ${\cal F}= \{ \varphi_B \}_{B \in \BB}$ is normal.
\end{cor}

\begin{pf}
By the surgery construction as described in Section \ref{sec:surgery}, $\cal F$ is uniformly quasiconformal. Let $\varphi_n=\varphi_{B_n}$ be a sequence in $\cal F$. Let $B_n=B_{\mu_n}$ and choose a subsequence, still denoted by $B_n$, such that $|\mu_n| \geq 1$ for all $n$ (the case $|\mu_n| \leq 1$ is similar). 
By the way we normalized $\varphi_n$,
$$\varphi_n(H_n^{-1}(0))=0, \ \ \ \varphi_n(1)=1, \ \ \ \varphi_n(\infty)=\infty.$$
But $\{ H_n^{-1}(0) \}$ lives in a compact subset of $\BBB D$. Hence the three points $H_n^{-1}(0)$, $1$ and $\infty$ has mutual spherical distance larger than some positive constant independent of $n$. This implies equicontinuity of $\{ \varphi_n \}$ by a standard theorem on quasiconformal mappings (\cite{Lehto}, Theorem 2.1).
\end{pf} 

Now we show that the surgery map constructed in the previous section is proper.

\begin{prop}
\label{proper}
The surgery map ${\cal S}:\BB \rightarrow \PP$ is proper.
\end{prop}
\begin{pf}
Let the sequence $\{ B_n \}$ leave every compact set in $\BB$ and consider the corresponding cubics $P_n={\cal S}(B_n)=\varphi_n \circ \tilde B_n \circ \varphi_n^{-1}$. To be more specific, let us assume that $B_n=B_{\mu_n}$ as in Section \ref{sec:blapar}, and the critical point $\mu_n$ tends to infinity. Clearly $P_n=P_{c_n}$, where $c_n=\varphi_n(\mu_n)$. Since $\{ \varphi_n \}$ is normal by the above corollary, we simply conclude that $c_n \rightarrow \infty$.
\end{pf}
 
\begin{prop}
\label{c5com}
$\CC$ is compact and invariant under the inversion $\mu \mapsto 1/\mu$. As a result, there exists an unbounded component $\Lambda_{ext}$ of $\BBB C^{\ast}\smallsetminus \CC$ which contains a punctured neighborhood of $\infty$ and a corresponding component $\Lambda_{int}$ which is mapped to it by $\mu \mapsto 1/\mu$.
\end{prop}
\begin{pf}
The invariance follows from the definition of $\BB$ and its identification with $\BBB C^{\ast}$. Note that the unit circle $\BBB T\subset \BB$ is contained in $\CC$ by \thmref{hitD}. So $\Lambda_{ext}$ and $\Lambda_{int}$ are actually distinct components of $\BBB C^{\ast}\smallsetminus \CC$. 

$\CC$ is clearly closed by \thmref{hitD}. Let us prove it is bounded. Assuming not, there is a sequence $B_{\mu_n}\in \CC$ with $\mu_n \rightarrow \infty$ as in the above proof. It follows from \propref{aux}(c) and \thmref{hitD} that the corresponding polynomials $P_{c_n}={\cal S}(B_{\mu_n})=\varphi_n \circ \tilde{B}_{\mu_n} \circ \varphi_n^{-1}$ have connected Julia sets. By \propref{m3com}, $1/30 \leq |c_n| \leq 30$. This contradicts properness of $\cal S$.
\end{pf}
\vspace{0.17in}

\section{Continuity of the Surgery Map}
\label{sec:continuity}

This section is devoted to the proof of continuity of the surgery map $\cal S$. This is by no means trivial, and in fact, as we will see, depends strongly on the cubic parameter space being one-dimensional. The fact that the cubics on the boundary of the connectedness locus $\MM$ are quasiconformally rigid is the most crucial step in the proof, and it is exactly this fact which makes the generalization of this work to higher degrees difficult. We would like to point out that the situation is similar to Douady-Hubbard's proof of the continuity of the ``straightening map'' in their study of the space of quadratic-like maps \cite{Douady-Hubbard2}. One additional difficulty here is the lack of complete information on quasiconformal conjugacy classes in the non-holomorphic family $\BB$ (the analogue of \thmref{qcclass}; see however \thmref{qcpath}). 

The idea of the proof is as follows: Given a sequence $B_n=B_{\mu_n}\in \BB$ such that $B_n\rightarrow B=B_{\mu}$, we prove that there exists a subsequence $\{ B_{n(j)} \}$ such that ${\cal S}(B_{n(j)})\rightarrow {\cal S}(B)$ in $\PP$. The topology of the parameter space $\PP$ is the uniform topology which respects the marking of the critical points. The same is true for $\BB$ with one exception (compare Section \ref{sec:surgery}): If $\mu$ has absolute value $1$, i.e., if $B$ has two double critical points on the unit circle, then $B_n\rightarrow B$ means that every subsequence of $\{ B_n \}$ has a further subsequence which either converges to $B$ or to its conjugate $R_{\mu}^{-1}\circ B\circ R_{\mu}$. From the construction of ${\cal S}$ it is easy to see that ${\cal S}(B)={\cal S}(R_{\mu}^{-1}\circ B\circ R_{\mu})$. Therefore, in order to prove continuity of ${\cal S}$, all we have to show is that $B_n\rightarrow B$ locally uniformly on $\BBB C$ (respecting the convergence of the marked critical points) implies
that for some subsequence $\{ B_{n(j)} \}$, ${\cal S}(B_{n(j)})\rightarrow {\cal S}(B)$ locally uniformly on $\BBB C$ (again, respecting the convergence of the marked critical points).    

So consider the sequence $\{ B_n|_{\BBB T} \}$ and let $h_n$ and $h$ be the unique $k(\theta)$-quasisymmetric homeomorphisms
which fix $z=1$ and conjugate $ B_n|_{\BBB T}$ and $B|_{\BBB T}$ to the rigid rotation $R_{\theta}$. It is easy to see that $h_n\rightarrow h$ uniformly
on $\BBB T$. Consider the Douady-Earle extensions $H_n$ and $H$, which are $K(\theta)$-quasiconformal homeomorphisms of the unit disk. By the construction of these extensions, $H_n$ and $H$ are real-analytic in $\BBB D$ and $H_n\rightarrow H$ locally uniformly in $C^{\infty}$ topology \cite{Douady-Earle}. In particular, the partial derivatives $\partial H_n$ and $\overline{\partial}H_n$ converge locally uniformly in $\BBB D$ to the corresponding derivatives $\partial H$ and $\overline{\partial}H$. This shows that $\sigma_n|_{\BBB D}\rightarrow \sigma|_{\BBB D}$ locally uniformly, where $\sigma_n$ and $\sigma$ are the conformal structures we constructed in the course of surgery for
$B_n$ and $B$ (see  Section \ref{sec:surgery}).

At this point, the main problem is to prove that $B_n\rightarrow B$
and $\sigma_n|_{\BBB D}\rightarrow \sigma|_{\BBB D}$ implies $\sigma_n\rightarrow \sigma$ in the $L^1$-norm on $\BBB C$, for this would show
that the normalized solutions $\varphi_n=\varphi_{H_n}$ of the Beltrami equations $\varphi_n^{\ast} \sigma_0=\sigma_n$ converge locally uniformly on $\BBB C$ to the normalized solution $\varphi$ of the equation $\varphi^{\ast} \sigma_0=\sigma$. This would simply mean that ${\cal S}(B_n)\rightarrow {\cal S}(B)$ as $n\rightarrow \infty$. 

Unfortunately, we cannot prove $\sigma_n\rightarrow \sigma$ in $L^1(\BBB C)$ in all cases. So, following \cite{Douady-Hubbard2}, we take a slightly different approach by splitting the argument into two cases depending on whether $\cal S (B)$ is quasiconformally rigid or not. In the first case, we show continuity directly using the rigidity. In the latter case, however, we prove $\varphi_n\rightarrow \varphi$ using the fact that $\cal S (B)$ admits nontrivial deformations.
\begin{thm}
\label{contin}
The surgery map $\cal S : \BB \rightarrow \PP$ is continuous.
\end{thm}
\begin{pf}
Consider $B_n,B\in \BB$ and start with the same construction as above to get a sequence $\{ \sigma_n\}$ of conformal structures on the plane with uniformly bounded dilatation and the corresponding sequence $\{ \varphi_n \}$
of normalized solutions of $\varphi_n^{\ast} \sigma_0=\sigma_n$. Since $\{ \varphi_n \}$ is a normal family by \corref{normal}, it has a subsequence, still denoted by $\{ \varphi_n \}$, which converges locally uniformly to a quasiconformal homeomorphism $\psi:\BBB C\rightarrow \BBB C$.
  
Set $P_n=\varphi_n\circ \tilde{B}_n\circ \varphi_n^{-1}={\cal S}(B_n)$, $P=\varphi\circ \tilde{B}\circ
\varphi^{-1}={\cal S}(B)$, and $Q=\psi \circ \tilde{B} \circ \psi^{-1}$.
All these maps are cubic polynomials in $\PP$. Also $P$ is quasiconformally conjugate to $Q$, and $P_n\rightarrow Q$ as $n \rightarrow \infty$. We will show that $P=Q$ and this will prove continuity at $B$. 

For the rest of the argument, we distinguish two cases: If $P=\cal S (B)$ is quasiconformally rigid, then automatically $P=Q$ and we are done. (By \thmref{qcclass} this case corresponds to the points on the boundary of $\MM$ or the centers of hyperbolic-like or capture components.) Otherwise, $P$ is not rigid, so the quasiconformal conjugacy class of $P$ is a nonempty open set $U\subset \PP$ by \corref{qcopen}. Assume by way of contradiction that $P\neq Q$. Since $P_n\rightarrow Q$ as $n\rightarrow \infty$, $P_n\in U$ for large $n$. Hence $P_n$ is quasiconformally conjugate to $P$ for large $n$, i.e., there exists a normalized quasiconformal homeomorphism $\eta_n:\BBB C\rightarrow \BBB C$ such that $\eta_n\circ P=P_n \circ \eta_n$. Observe that the dilatation of $\eta_n$ is uniformly bounded, since by \thmref{qcpar} the dilatation of $(\psi\circ \varphi^{-1})\circ \eta_n^{-1}$ goes to $1$ as $n$ goes to $\infty$ (see \figref{arrow}).
By ``lifting'' $\eta_n$, we can find a quasiconformal conjugacy $\xi_n=\varphi_n^{-1}\circ \eta_n \circ \varphi$ between the modified Blaschke products $\tilde B$ and $\tilde {B}_n$, i.e., 
\begin{equation}
\label{eqn:lift}
\xi_n\circ \tilde B=\tilde {B}_n \circ \xi_n.
\end{equation}
Again, note that the dilatation of $\xi_n$ is uniformly bounded.
\realfig{arrow}{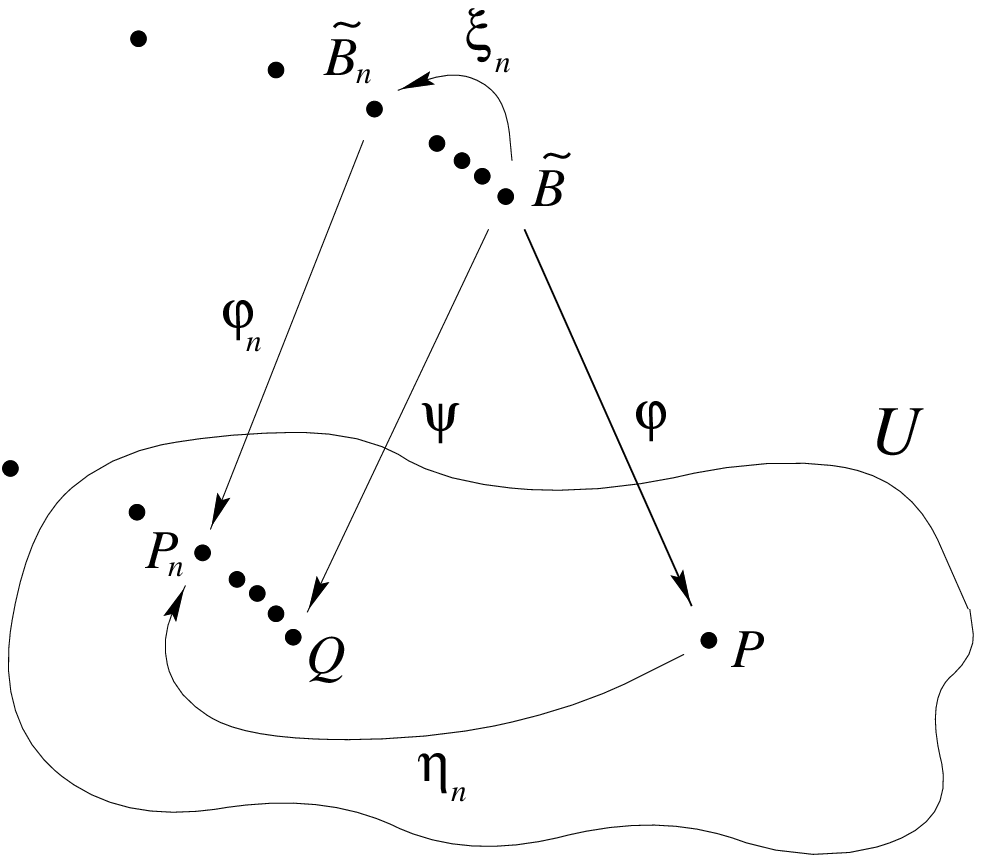}{{\sl Sketch of the proof of continuity of $\cal S$.}}{7 cm}

We prove that the sequence of conformal structures $\{ \sigma_n \} $ converges 
in $L^1(\BBB C)$ to $\sigma$. This, by a standard theorem on quasiconformal mappings (see for example \cite{Lehto}, Theorem 4.6), will show that $\varphi_n\rightarrow \varphi$ locally uniformly,
hence $P_n\rightarrow P$, hence $P=Q$, which contradicts our assumption.

To this end, we introduce the following sequences of conformal structures (where, as usual, we identify a conformal structure with its associated Beltrami differential):
$$\sigma_n^k(z)= \left \{ 
\begin{array}{ll}
\sigma_n(z) & \mbox{when}\ z\in \bigcup_{i=0}^k \tilde{B}_n^{-i}(\BBB D)\\
0           & \mbox{otherwise}
\end{array}
\right. $$
and
$$\sigma^k(z)= \left \{ 
\begin{array}{ll}
\sigma(z) & \mbox{when}\ z\in \bigcup_{i=0}^k \tilde{B}^{-i}(\BBB D)\\
0           & \mbox{otherwise}
\end{array}
\right. $$
Note that  $\sigma^k\rightarrow \sigma $ in $L^1(\BBB C)$ as $k\rightarrow \infty$ and for every fixed $k$, 
$\sigma_n^k\rightarrow \sigma^k$ in $L^1(\BBB C)$ as $n\rightarrow \infty$.
\begin{lem}
\label{area}
The $L^1$-norm $\| \sigma_n -\sigma \|_1$ goes to zero as
$n\rightarrow \infty$ if the area of the open set $\bigcup_{i=k}^{\infty}\tilde{B}_n^{-i}(\BBB D)$ goes to zero uniformly in $n$
as $k\rightarrow \infty$.
\end{lem}
\begin{pf}
For a given $\epsilon >0$, take $k_0$ so large that $k>k_0$
implies $area(\bigcup_{i=k}^{\infty}\tilde{B}_n^{-i}(\BBB D))<\epsilon $ for
all $n$. Then for a fixed large $k>k_0$ and $n$ large enough,
$$\begin{array}{rl}
\| \sigma_n -\sigma \|_1 \leq & \| \sigma_n -\sigma_n^k \|_1 + \| \sigma_n^k -\sigma^k \|_1 + \| \sigma^k -\sigma \|_1 \\
                         \leq & \| \sigma_n -\sigma_n^k \|_1 +2\epsilon \\
                          =   & \displaystyle{\int}_{\bigcup_{i=k+1}^{\infty} \tilde{B}_n^{-i}(\BBB D)}|\sigma_n -\sigma_n^k |\ dxdy + 2\epsilon \\
                           <  & 4\epsilon.
\end{array}$$
This completes the proof of the lemma.
\end{pf}

So it remains to prove that the area of $\bigcup_{i=k}^{\infty}\tilde{B}_n^{-i}(\BBB D)$ goes to zero uniformly in $n$
as $k\rightarrow \infty$. Clearly $area(\bigcup_{i=k}^{\infty}\tilde{B}^{-i}(\BBB D))\rightarrow 0$ as $k\rightarrow \infty$. Since $\{ \xi_n \}$ is uniformly quasiconformal, there
is a constant $C>0$ such that 
$$C^{-1}\ area(E) \leq area(\xi_n(E))\leq C\ area(E)$$
for any measurable set $E$. By (\ref{eqn:lift}), 
$$\bigcup_{i=k}^{\infty}\tilde{B}_n^{-i}(\BBB D)=\xi_n(\bigcup_{i=k}^{\infty}\tilde{B}^{-i}(\BBB D)),$$
so $area(\bigcup_{i=k}^{\infty}\tilde{B}_n^{-i}(\BBB D))\leq C\ area(\bigcup_{i=k}^{\infty}\tilde{B}^{-i}(\BBB D))$ and this proves that the left side goes to zero uniformly in $n$.
\end{pf}

\vspace{0.17in}

\section{Renormalizable Blaschke Products}
\label{sec:renbla}

Here we consider those Blaschke products in $\BB$ out of which one can ``extract'' the standard degree $3$ Blaschke product $\FT$ to be defined below. The importance of this particular Blaschke product comes from the fact that it provides a model for the dynamics of the quadratic polynomial $\QQ$. It will be convenient to define renormalizable Blaschke products in $\BB$ as ones which after the surgery give rise to renormalizable cubics in $\PP$ (see  Section \ref{sec:rencub}). In what follows we will have to work with a {\it symmetrized} version of the notion of a quadratic-like map in order to show that any renormalizable Blaschke product is quasiconformally conjugate near the Julia set of its renormalization to the standard map $\FT$. The proof of this fact resembles the proof of \cite{Douady-Hubbard2} that every hybrid class of polynomial-like maps contains a polynomial.

First we include the following simple fact for completeness. 
\begin{prop}
\label{bla3}
Let $0<\theta<1$ be a given irrational number and $f:\overline{\BBB C}\rightarrow\overline{\BBB C}$ be a degree $3$ Blaschke product with a superattracting fixed point at the origin and a double critical point at $z=1$. Let the rotation number of $f|_{\BBB T}$ be $\theta$. Then there exists a unique $0<t(\theta)<1$ such that
\begin{equation}
\label{eqn:ftet}
f(z)=\FT (z)=e^{2 \pi i t(\theta)} z^2 \left ( \frac{z-3}{1-3z} \right ).
\end{equation}
\end{prop}
\begin{pf}
Clearly $f(z)=e^{2 \pi i t} z^2 \displaystyle{ \left ( \frac{z-a}{1-\overline{a}z} \right ) }$, with $|a|>1$ and $0<t<1$. The fact that $f'(1)=0$ implies $a=3$. The rotation number of $f|_{\BBB T}$ as a function of $t$ is continuous and strictly monotone at all irrational values \cite{Katok}. Hence there exists a unique $t(\theta)$ for which this rotation number is $\theta$. 
\end{pf}

\realfig{ft}{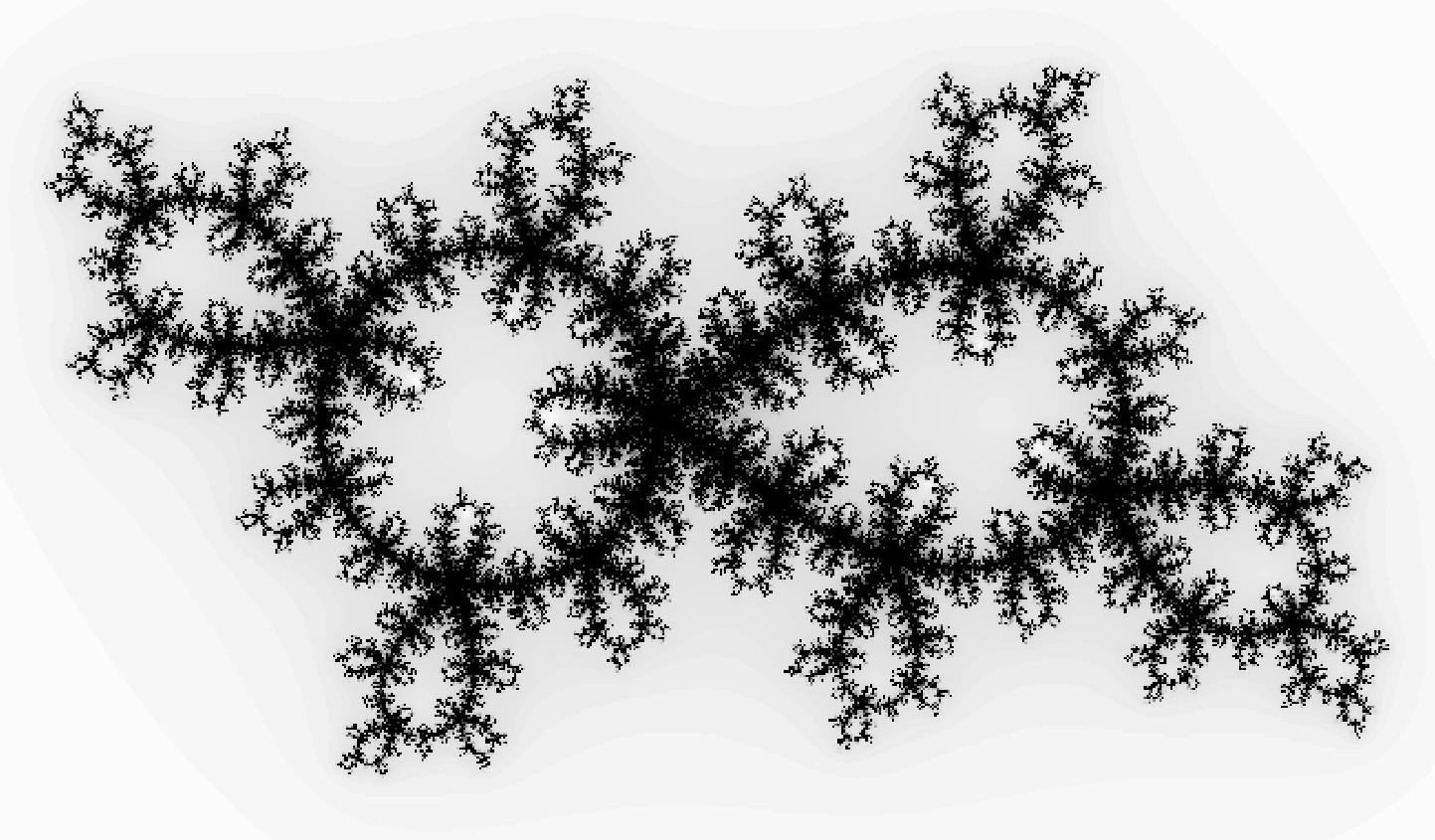}{{\sl The Julia set of $\FT$ for $\theta=(\sqrt{5}-1)/2$.}}{9 cm}

\noindent 
{\bf Remark.} Computer experiments give the value $t(\theta)=0.613648$ for the golden mean $\theta=(\sqrt{5}-1)/2$. \figref{ft} shows the Julia set of $\FT$ for this value of $\theta$. This standard degree $3$ Blaschke product was introduced by Douady, Ghys, Herman and Shishikura as a model for the quadratic $\QQ$ in the case $\theta$ is irrational of bounded type \cite{Douady2}. 
It was also used in \cite{Petersen} to prove that the Julia set of $Q_{\theta}$ is locally connected and has measure zero. \\ \\
{\bf Definition.} A Blaschke product $B\in \BB$ is called {\it renormalizable} if ${\cal S}(B)\in \PP$ is a renormalizable cubic, as defined in  Section \ref{sec:rencub}.
\begin{thm}
\label{b3like}
Let $B\in \BB$ be renormalizable. Then there exists a pair of annuli $W'\Subset W$, both containing the unit circle and symmetric with respect to it, and a quasiconformal homeomorphism $\varphi_B:\BBB C\rightarrow \BBB C$ such that:
\begin{enumerate}
\item[(a)]
$B:\partial W'\rightarrow \partial W$ is a degree $2$ covering map,
\item[(b)]
$\varphi_B\circ I=I\circ \varphi_B,$
\item[(c)]
$(\varphi_B\circ B)(z)=(\FT\circ \varphi_B)(z)$ for all $z\in W'$.
\end{enumerate}
Moreover, $\varphi_B$ can be chosen to be conformal (i.e., $\overline{\partial} \varphi_B =0$) on $K(B)=\bigcap_{n\geq 0}B^{-n}(W')$.
\end{thm}
\begin{pf}
Consider the cubic $P={\cal S}(B)=\varphi \circ \tilde{B}\circ \varphi^{-1}\in \PP$ which is renormalizable. Consider the quadratic-like restriction $P|_U:U\rightarrow V$ and the corresponding regions $U_1=\varphi^{-1}(U)$ and $V_1=\varphi^{-1}(V)$. Clearly $U_1\Subset V_1$ and both contain the closed unit disk. Define the symmetrized regions
$$W'=U_1\cap I(U_1),\ \ \ \ \ W=V_1\cap I(V_1)$$
which are topological annuli with $W'\Subset W$. Note that $B$ sends $\partial W'$ to $\partial W$ in a 2-to-1 fashion. 

Now extend $B|_{W'}$ to the whole complex plane by gluing it to the polynomial $z\mapsto z^2$ near $0$ and $\infty$ as follows: Let $r>1$ and  $\omega :\BBB C\smallsetminus W'\rightarrow \BBB C\smallsetminus \BBB A (r^{-1},r)$ be a diffeomorphism such that
$$\begin{array}{ll}
\omega\circ I=I\circ \omega, & \\
\omega (B(z))=\omega(z)^2, & z\in \partial W'.
\end{array}$$
Define the extension of $B|_{W'}$ by
$$F(z)= \left \{ 
\begin{array}{ll}
B(z) & z\in W'\\
\omega^{-1}(\omega(z)^2) & z\notin W'
\end{array}
\right. $$
Note that $F$ is a quasiregular degree $3$ self-map of the sphere, $F\circ I=I\circ F$, and every point outside $W'$ will converge to $0$ or $\infty$ under the iteration of $F$.  

Define a conformal structure $\sigma$ on the plane as follows: Put $\sigma=\omega^{\ast}\sigma_0$ on $\BBB C\smallsetminus W'$, and pull it back by $F^{\circ n}$
to all the components of $F^{-n}(\BBB C\smallsetminus W')\cap W'$. Finally, on $K(B)$ set $\sigma=\sigma_0$. It is easy to see that $\sigma$ has bounded dilatation on the plane, is symmetric with respect to the unit circle, and $F^{\ast}(\sigma)=\sigma$. By the Measurable Riemann Mapping Theorem of Ahlfors and Bers, there exists a unique quasiconformal homeomorphism $\varphi_B$ of the plane which fixes $0, 1, \infty$, such that $\varphi_B^{\ast}(\sigma_0)=\sigma$. The conjugate map $f=\varphi_B\circ F\circ \varphi_B^{-1}$ is easily seen to be a degree $3$ rational map on the sphere. The quasiconformal homeomorphism $I\circ \varphi_B\circ I$ also fixes $0,1,\infty$ and pulls $\sigma_0$ back to $\sigma$ because $\sigma$ is symmetric with respect to $\BBB T$. By uniqueness, $\varphi_B=I\circ \varphi_B\circ I$. This implies that $f$ commutes with $I$, hence it is a Blaschke product. By \propref{bla3}, $f=\FT$, and we are done. 
\end{pf}

While the above theorem establishes a direct connection between some Blaschke products in $\BB$ and $\FT$, it is curious to note the following entirely different relation:

\begin{thm}
\label{q=3}
Let $B_n=B_{\mu_n}$ be any sequence in $\BB$ such that $\mu_n \rightarrow \infty$ as $n \rightarrow \infty$. Then $B_n\rightarrow \FT$ locally uniformly on $\BBB C$ as $n \rightarrow \infty$.
\end{thm}

In other words, $\FT$ can be regarded as the point at infinity of the parameter space $\BB$. 

\begin{pf}
As in Section \ref{sec:blapar}, let
$$ B_n :z\mapsto e^{2 \pi i t_n} z^3 \left ( \frac{z-p_n}{1-\overline{p}_n z} \right ) \left ( \frac{z-q_n}{1-\overline{q}_nz} \right ) .$$
The logarithmic first derivative $B_n'/B_n$ and second derivative $(B_n B_n'' -(B_n')^2)/(B_n)^2$ both vanish at $z=1$. A straightforward computation shows that these two conditions translate into
\begin{equation}
\label{eqn:yek}
\frac{|p_n|^2-1}{|p_n-1|^2}+\frac{|q_n|^2-1}{|q_n-1|^2}=3,
\end{equation}
and  
\begin{equation}
\label{eqn:do}
\frac{(p_n-\overline{p}_n)(|p_n|^2-1)}{|p_n-1|^4}+\frac{(q_n-\overline{q}_n)(|q_n|^2-1)}{|q_n-1|^4}=0.
\end{equation}
Let us write $a_n \leadsto a$ when $a$ is an accumulation point of the sequence $a_n$. Since $\mu_n \rightarrow \infty$, both $p_n$ and $q_n$ cannot stay bounded. Hence one of them, say $p_n$ gets arbitrarily large, or $p_n \leadsto \infty$. Then (\ref{eqn:yek}) shows that $(|q_n|^2-1)/|q_n-1|^2 \leadsto 2$, or equivalently, $|q_n-2| \leadsto 1$ but $q_n$ stays away from $z=1$. On the other hand, (\ref{eqn:do}) shows that $(q_n-\overline{q}_n)(|q_n|^2-1)/|q_n-1|^4 \leadsto 0$, hence $(q_n-\overline{q}_n)/|q_n-1|^2 \leadsto 0$. Since $q_n$ does not accumulate on $z=1$, this implies that $(q_n-\overline{q}_n) \leadsto 0$. Near the circle $|z-2|=1$ this can happen only if $q_n \leadsto 3$. 

We have shown that there exists a subsequence $B_{n(j)}$ such that $p_{n(j)} \rightarrow \infty$ and $q_{n(j)} \rightarrow 3$ as $j\rightarrow \infty$. Since the rotation number depends continuously on the circle map, it is easy to see that this implies $B_{n(j)} \rightarrow \FT$ locally uniformly on $\BBB C$. 
\end{pf}

Consider a sequence $B_n=B_{\mu_n}$ going off to infinity as in the previous theorem. Consider the cubics $P_n=P_{c_n}={\cal S}(B_n)=\varphi_n \circ \tilde{B}_n \circ \varphi_n^{-1}$ as in (\ref{eqn:pmodb}). By the previous theorem, $B_n\rightarrow \FT$, so $\tilde{B}_n\rightarrow \tilde{\FT}$. Since $\{ \varphi_n \}$ is normal by \corref{normal}, by passing to a subsequence if necessary, $\varphi_n$ converges to a quasiconformal homeomorphism $\varphi$. Since the surgery map is proper by \propref{proper}, $c_n \rightarrow \infty$. By examining the normal form (\ref{eqn:normform}), we see that $P_n \rightarrow Q$, where $Q:z\mapsto \lambda z (1-1/2 z)$ is affinely conjugate to $\QQ$. Hence, $Q=\varphi \circ \tilde{\FT} \circ \varphi^{-1}$ and we recover the surgery introduced by Douady and others. We conclude that the surgery map ${\cal S}: \BB \rightarrow \PP$ extends continuously to the points at infinity of both parameter spaces, and the extension is also a surgery.\\  

The next theorem is the analogue of \thmref{qcpar} for Blaschke products. It will be more convenient to formulate it for a general Blaschke product since we would like to use it for $\FT$ as well as the elements of $\BB$.
\begin{thm}[Paths of QC Conjugacies]
\label{qcpath}
Let $A$ and $B$ be two Blaschke products of degree $d$ and let $\Phi$ be a quasiconformal homeomorphism which fixes $0,1,\infty$ such that $\Phi \circ I=I\circ \Phi$ and $\Phi \circ A=B\circ \Phi$. Then there exists a path $\{ \Phi_t \}_{0\leq t\leq 1}$ of quasiconformal homeomorphisms, with $\Phi_0=id$ and $\Phi_1=\Phi$, such that $A_t=\Phi_t\circ A\circ \Phi_t^{-1}$ is a Blaschke product for every $0\leq t\leq 1$. In particular, either $A$ is quasiconformally rigid or its conjugacy class is nontrivial and path-connected.
\end{thm}
\begin{pf}
The proof is almost identical to that of \thmref{qcpar}. 
Consider $\sigma=\Phi^{\ast}\sigma_0$, which is invariant under $A$, and take the {\it real} perturbations $\sigma_t=t\sigma$, $0\leq t\leq 1$. Let $\Phi_t$ be the unique quasiconformal homeomorphism which fixes $0,1,\infty$ and satisfies $\Phi_t^{\ast}\sigma_0=\sigma_t$. The map $A_t=\Phi_t\circ A\circ \Phi_t^{-1}$ is easily seen to be a degree $d$ rational map. By uniqueness, $I\circ \Phi_t \circ I=\Phi_t$ since the left-hand side also pulls $\sigma_0$ back to $\sigma_t$ and fixes $0,1,\infty$. Hence $A_t$ commutes with $I$. So it is a Blaschke product.
\end{pf}  

We will need the next lemma in the proof of \thmref{inj}.
\begin{lem}[Rigidity on the Julia Set]
\label{rig}
Let $\psi$ be a quasiconformal homeomorphism defined on an open annulus containing the Julia set $J(\FT)$ of the Blaschke product $\FT$ defined in $($\ref{eqn:ftet}$)$. Suppose that $\psi$ commutes with $I$ and conjugates $\FT$ to itself. Then $\psi|_{J(\FT)}$ is the identity.
\end{lem}
\begin{pf}
Extend $\psi$ to a quasiconformal homeomorphism $\BBB C\rightarrow \BBB C$ which commutes with $I$ and conjugates $\FT$ to itself. By the previous theorem, there exists a path $t\mapsto \psi_t$ of quasiconformal homeomorphisms, with $0\leq t\leq 1$ and $\psi_0=id, \psi_1=\psi$, such that $\psi_t\circ \FT \circ \psi_t^{-1}$ is a degree $3$ Blaschke product quasiconformally conjugate to $\FT$. By \propref{bla3}, this Blaschke product has to be $\FT$ itself, so $\psi_t$ commutes with $\FT$. 

Now for any periodic point $z\in J(\FT)$ of period $n$, $t\mapsto \psi_t(z)$ is a continuous path in the finite set of all period-$n$ points in $J(\FT)$. Since $\psi_0(z)=z$, we must have $\psi(z)=z$. Since such points $z$ are dense in the Julia set, $\psi|_{J(\FT)}$ must be the identity.
\end{pf}
\vspace{0.17in}

\section{On Injectivity of the Surgery Map}
\label{sec:inject}

In this section we prove that the surgery map $\cal S: \BB \rightarrow \PP$ is injective on the set of Blaschke products which map to $\BBB C^{\ast}\smallsetminus \MM$ or to hyperbolic-like cubics. The proof of this fact is based on the combinatorics of drops and their nuclei as developed in  Section \ref{sec:blapar}. Here is the outline of the proof: If
$\cal S (A)=\cal S (B)$ for some $A,B\in \BB$, there exists a quasiconformal homeomorphism of the plane which conjugates the modified Blaschke products $\tilde A$ and $\tilde B$, which is conformal everywhere except on the union of the maximal drops. A careful analysis will then show that when $\cal S(A)$ is not capture, one can redefine this homeomorphism on all the drops of the two Blaschke products to get a conjugacy between $A$ and $B$ 
everywhere. A pull-back argument together with the Bers Sewing Lemma at each step shows that this conjugacy is conformal away from the Julia sets (\thmref{ABconj}). When ${\cal S}(A)$ is hyperbolic-like or has disconnected Julia set,
one can use the renormalization scheme in  Section \ref{sec:renbla} and the rigidity on the Julia sets (\lemref{rig}) to conclude that the conjugacy between $A$ and $B$ is in fact conformal (\thmref{inj}). The main \thmref{main} and some corollaries on the connectedness locus $\CC$ will follow immediately.
\begin{thm}
\label{ABconj}
Let $A,B\in \BB$ and $\cal S (A)=\cal S (B)=P$. Suppose that $P$ is not capture. Then there exists a quasiconformal homeomorphism $\Phi: \overline{\BBB C} \rightarrow \overline{\BBB C}$ which fixes $0,1,\infty$, commutes with $I$, and conjugates $A$ to $B$. Moreover, $\Phi$ is conformal on the Fatou set $\overline{\BBB C} \smallsetminus J(A)$.
\end{thm}
\begin{pf}
Following the notation of (\ref{eqn:pmodb}), we assume that $P=\varphi \circ \tilde A \circ \varphi^{-1}=\varphi' \circ \tilde B \circ {\varphi'}^{-1}$ for some quasiconformal homeomorphisms $\varphi$ and $\varphi'$. Consider the quasiconformal homeomorphism $\Phi_0={\varphi'}^{-1}\circ \varphi$ which conjugates $\tilde A$ to $\tilde B$ on the entire plane and is conformal (i.e., $\overline{\partial}\Phi_0=0$) everywhere except on $\bigcup_{k\geq 0}\tilde A^{-k}(\BBB D)$. 

Note that by \propref{aux}(b) the open set $\overline{\BBB C} \smallsetminus \overline{\bigcup_{k\geq 0}\tilde A^{-k}(\BBB D)}$ is precisely the nucleus $N_{\infty}$ as defined in  Section \ref{sec:blapar}. Also, $\bigcup_{k\geq 0}\tilde A^{-k}(\BBB D)$ is the {\it disjoint} union of the maximal drops of $A$ (which by \propref{aux}(a) correspond to the bounded Fatou components of $P$ which map to the Siegel disk $\Delta_P$). Similar correspondence holds for the open set $\bigcup_{k\geq 0}\tilde B^{-k}(\BBB D)$. Therefore, corresponding to any maximal $k$-drop $D_k^i(A)$, there exists a unique maximal
$k$-drop $D_k^i(B)=\Phi_0 (D_k^i(A))$. Finally, note that for any such maximal drops, $A^{\circ k}:D_k^i(A)\rightarrow \BBB D$ and $B^{\circ k}:D_k^i(B)\rightarrow \BBB D$ are conformal isomorphisms since by our assumption $P$ is not capture.  

In what follows we construct a sequence of quasiconformal homeomorphisms $\Phi_n$ which preserve the unit circle $\BBB T$ and another sequence $\Upsilon_n$ by {\it symmetrizing } each $\Phi_n$:
$$\Upsilon_n(z)=  \left \{ \begin{array}{ll} 
\Phi_n(z) &  |z|\geq 1 \\
(I\circ \Phi_n\circ I)(z) & |z|<1
\end{array}
\right. $$
We have already constructed $\Phi_0$, hence $\Upsilon_0$. Consider the sequences of compact sets $\{ J_n(A) \}$ and $\{ J_n(B) \}$ as in \lemref{alt}. Note that $\Phi_0 \circ A=B\circ \Phi_0$ on $J_0(A)$. The next step is to define $\Phi_1$: Let $\Phi_1=\Upsilon_0$ everywhere except on the maximal drops of $A$. On any maximal $k$-drop $D_k^i(A)$ we define $\Phi_1:D_k^i(A)\rightarrow D_k^i(B)$ by $B^{-k}\circ \Upsilon_0 \circ A^{\circ k}$. (When $k=0$, the only maximal $0$-drop is $\BBB D$ and by this definition $\Phi_1|_{\BBB D}=\Upsilon_0|_{\BBB D}$.) Observe that the two definitions match along the common boundary. Hence $\Phi_1$ is in fact a quasiconformal homeomorphism by the Bers Sewing Lemma.  Note that $\Phi_1|_{J_0(A)}=\Phi_0|_{J_0(A)}$ and by definition of $J_1(A)$ in (\ref{eqn:jn}), $\Phi_1 \circ A=B\circ \Phi_1$ on $J_1(A)$. The homeomorphism $\Upsilon_1$ is then obtained by symmetrizing $\Phi_1$. 

Continuing inductively, we define $\Phi_n$ to be equal to $\Upsilon_{n-1}$ everywhere except on the maximal drops of $A$ and then on the maximal drops we define it by taking pull-backs. In other words, $\Phi_n:D_k^i(A)\rightarrow D_k^i(B)$ will be defined by $B^{-k}\circ \Upsilon_{n-1} \circ A^{\circ k}$.
\begin{lem}
\label{induct}
The sequence of quasiconformal homeomorphisms $\{ \Phi_n \}$ has the following properties:
\begin{equation}
\label{eqn:ind1}
\Phi_n|_{J_{n-1}(A)}=\Phi_{n-1}|_{J_{n-1}(A)},
\end{equation}
and
\begin{equation}
\label{eqn:ind2} 
(\Phi_n \circ A)(z)=(B\circ \Phi_n)(z)\hskip 1cm  z\in J_n(A).
\end{equation}
\end{lem}
\begin{pf}
Both properties follow by induction on $n$. Let us prove (\ref{eqn:ind1}) first. We have already seen (\ref{eqn:ind1}) for $n=1$. Assume (\ref{eqn:ind1}) is true and let $z\in J_n(A)$. We distinguish three cases:\\

$\bullet${\it Case 1:} $z\in J_n(A)\cap \overline{\BBB D}$. Then $I(z)\in J_{n-1}(A)$ and we have $\Phi_{n+1}(z)=\Upsilon_n(z)=(I\circ \Phi_n \circ I)(z)=(I\circ \Phi_{n-1} \circ I)(z)$ by the induction hypothesis. The latter is clearly equal to $\Upsilon_{n-1}(z)=\Phi_n(z)$.\\
 
$\bullet${\it Case 2:} $z\in J_n(A)\smallsetminus \overline{\BBB D}$ and $A^{\circ k}(z)\in \overline{\BBB D}$ for some $k\geq 1$. $A^{\circ k}(z)\in IJ_{n-1}$ and hence $(I\circ A^{\circ k})(z)\in J_{n-1}(A)$. So $\Phi_{n+1}(z)=(B^{-k}\circ \Upsilon_n \circ A^{\circ k})(z)=(B^{-k}\circ I\circ \Phi_n \circ I \circ A^{\circ k})(z)=(B^{-k}\circ I\circ \Phi_{n-1} \circ I \circ A^{\circ k})(z)$ by the induction hypothesis. Again, the latter is equal to $(B^{-k}\circ \Upsilon_{n-1} \circ A^{\circ k})(z)=\Phi_n(z)$.\\

$\bullet${\it Case 3:} $z\in J_n(A)\smallsetminus \overline{\BBB D}$ and $z$ is accumulated by points of the form {\it Case 2}. Then, clearly, $\Phi_{n+1}(z)=\Phi_n(z)$ by continuity.\\ 

Altogether the three steps show that $\Phi_{n+1}|_{J_n(A)}=\Phi_n|_{J_n(A)}$, which completes the induction step and the proof of (\ref{eqn:ind1}). 

To prove (\ref{eqn:ind2}) we have to work a little bit more. We have already seen (\ref{eqn:ind2}) for $n=1$. Assume (\ref{eqn:ind2}) is true and let $z\in J_{n+1}(A)$. We split the induction step into the following cases:\\ 

$\bullet${\it Case 1:} $z\in J_{n+1}(A)\smallsetminus \overline{\BBB D}$ and $A(z)
\notin \overline{\BBB D}$. Then $(\Phi_{n+1} \circ A)(z)=(B\circ \Phi_{n+1})(z)$ automatically since $\Phi_{n+1}$ is defined by pull-backs.\\

$\bullet${\it Case 2:} $z\in J_{n+1}(A)\smallsetminus \overline{\BBB D}$ but $A(z)
\in \overline{\BBB D}$. Then $(\Phi_{n+1}\circ A)(z)=(\Upsilon_n\circ A)(z)=
(B\circ B^{-1}\circ \Upsilon_n\circ A)(z)=(B\circ \Phi_{n+1})(z).$\\

$\bullet${\it Case 3:} $z\in J_{n+1}(A)\cap \overline {\BBB D}$ and $A(z)\in \overline {\BBB D}$. 
Then $(\Phi_{n+1}\circ A)(z)=(\Upsilon_n\circ A)(z)=(I\circ \Phi_n\circ I)(A(z))=(I\circ \Phi_n \circ A)(I(z))$. But $I(z)\in J_n(A)$ so by the induction hypothesis, $(I\circ \Phi_n\circ A)(I(z))=(I\circ B\circ \Phi_n)(I(z))=(B\circ I\circ \Phi_n)(I(z))=(B\circ \Upsilon_n)(z)=(B\circ \Phi_{n+1})(z)$.\\
 
$\bullet${\it Case 4:} $z\in J_{n+1}(A)\cap \overline {\BBB D}$ but $A(z)\notin \overline {\BBB D}$. Then $I(z)\in J_n(A)$. Let $w=A(z)$. Since $A(I(z))=I(w)\in \BBB D$, we have $I(w)\in IJ_{n-1}(A)$, hence $w\in J_{n-1}(A)$. By (\ref{eqn:ind1}), one has $\Phi_{n+1}(w)=\Phi_n(w)=\Phi_{n-1}(w)=\Upsilon_{n-1}(w)=(I\circ \Upsilon_{n-1}\circ I)(w)=(I\circ \Phi_n\circ I)(w)=(I\circ \Phi_n\circ I)(A(z))=(I\circ \Phi_n\circ A)(I(z))=(I\circ B\circ \Phi_n)(I(z))$ by the induction hypothesis. The latter is equal to $(B\circ I\circ \Phi_n)(I(z))=(B\circ \Upsilon_n)(z)=(B\circ \Phi_{n+1})(z)$.
\end{pf} 
  
Back to the proof of \thmref{ABconj}. By the Bers Sewing Lemma, the symmetrization $\Phi_n \longrightarrow \Upsilon_n$ does not increase the dilatation. On the other hand, the modification $\Upsilon_n \longrightarrow \Phi_{n+1}$ achieved by pull-backs along the maximal drops does not increase the dilatation either, simply because $A$ and $B$ are holomorphic. So we may assume that $\{ \Phi_n \}$ is uniformly quasiconformal. Since all the $\Phi_n$ fix $0,1,\infty$, it follows that some subsequence $\Phi_{n(j)}$ converges locally uniformly to a quasiconformal homeomorphism $\Phi$. \lemref{alt} and \lemref{induct} imply that $\Phi \circ A=B\circ \Phi$ on $J(A)$. 

In particular, this shows that $\Phi$ sends {\it all} the drops of $A$ bijectively to the drops of $B$ (before we only had a correspondence between the {\it maximal} drops of $A$ and $B$). 

It is easy to check that $\Phi$ obtained this way is conformal on the union $N=\bigcup_{i,k} N_k^i(A)$ of all the nuclei of drops of $A$ at all depths as defined in Section (\ref{sec:blapar}) and in fact conjugates $A$ to $B$ there. Since $N$ is clearly disjoint from the Julia set $J(A)$ by (\ref{eqn:inv}), it remains to show that every Fatou component of $A$ is contained in $N$. 

Consider a component $U$ of the Fatou set of $A$. Under the iteration of $A$, $U$ visits both $\BBB D$ and $\BBB C \smallsetminus \overline{\BBB D}$ either finitely many times or infinitely often. In the first case, $U$ has to eventually map into the nucleus $N_0(A)$ or $N_{\infty}(A)$, hence it has to be contained in $N$. We prove that the second case cannot occur. In fact, suppose that the orbit of $U$ visits $\BBB D$ and $\BBB C \smallsetminus \overline{\BBB D}$ infinitely often. According to Sullivan \cite{Sullivan1}, $U$ eventually maps to a periodic Fatou component of $A$ which is either an attracting or parabolic basin or a Siegel disk or a Herman ring. It follows that this cycle of periodic Fatou components intersects both $\BBB D$ and $\BBB C \smallsetminus \overline{\BBB D}$, so in either case a critical point of $A$ has to enter $\BBB D$ and escapes from it infinitely often, which is impossible since ${\cal S}(A)$ is not a capture. This shows that $N=\overline{\BBB C}\smallsetminus J(A)$ and proves that $\Phi$ is a conjugacy between $A$ and $B$ everywhere and is conformal on $\overline{\BBB C}\smallsetminus J(A)$. It is easy to see that $\Phi$ constructed this way commutes with $I$.
\end{pf}
\begin{thm}
\label{inj}
Let $A,B \in \BB$ and ${\cal S}(A)={\cal S}(B)$. If ${\cal S}(A)$ is hyperbolic-like or has disconnected Julia set, then $A=B$.
\end{thm}
\begin{pf}
$A$ and $B$ are renormalizable by \thmref{ren}. Consider the quasiconformal homeomorphism $\Phi$ given by \thmref{ABconj}. By \thmref{b3like}, there exists a pair of annuli $W'_A\Subset W_A$ (resp. $W'_B\Subset W_B$) and a quasiconformal homeomorphism $\varphi_A$ (resp. $\varphi_B$) which conjugates $A$ (resp. $B$) to $\FT$ on $W'_A$ (resp. $W'_B$). Since ${\cal S}(A)={\cal S}(B)$, we can assume that $W'_B=\Phi (W'_A)$ and $W_B=\Phi (W_A)$.
The quasiconformal homeomorphism $\psi=\varphi_B \circ \Phi \circ \varphi_A^{-1}: \varphi_A(W'_A)\rightarrow \varphi_B(W'_B)$ is a self-conjugacy of $\FT$ near its Julia set which commutes with $I$. By \lemref{rig}, we must have $\psi|_{J(\FT)}=id$. It follows from the Bers Sewing Lemma that the $\overline{\partial}$-derivative of $\psi$ is zero almost everywhere on $J(\FT)$. Since by \thmref{b3like}(b) $\varphi_A$ (resp. $\varphi_B$) has zero $\overline{\partial}$-derivative on $K(A)$ (resp. $K(B)$), we conclude that $\overline{\partial}\Phi=0$ almost everywhere on $K(A)$. But, as in the proof of \corref{measure}, up to a set of measure zero, 
$J(A)=\bigcup_{n\geq 0} A^{-n}(K(A))$. Therefore, $\overline{\partial}\Phi$ has to be zero almost everywhere on the Julia set $J(A)$. Hence $\Phi$ is conformal, so $A=B$.
\end{pf}

\noindent
{\bf Remark.} We believe that the surgery map is a homeomorphism, at least outside of the capture components where it might have branching. This would imply that the connectedness loci $\CC$ and $\MM$ are actually homeomorphic, a conjecture that is strongly supported by computer experiments.
 
\begin{cor}
\label{c5con}
The surgery map $\cal S$ restricts to a homeomorphism 
$\Lambda_{ext}\iso \Omega_{ext}$. Similar conclusion holds for $\Lambda_{int}$ and $\Omega_{int}$. In particular, the connectedness locus $\CC$ is connected.
\end{cor}
\begin{pf}
Clearly $\cal S$ maps $\Lambda_{ext}$ into $\Omega_{ext}$ injectively by the previous theorem. Since $\cal S$ is a proper map by \propref{proper}, it extends to a continuous injection $\Lambda_{ext}\cup \{ \infty \} \hookrightarrow \Omega_{ext}\cup \{ \infty \}$. We claim that this injection is onto. To this end, it suffices to show that for any sequence $B_n \in \Lambda_{ext}$ which converges to the boundary of the connectedness locus $\CC$, the sequence $P_n={\cal S}(B_n)\in \Omega_{ext}$ converges to the boundary of $\MM$. If not, there is a subsequence of $B_n$ which converges to $B\in \partial \CC$ but the corresponding subsequence of $P_n$ converges to some $P\in \Omega_{ext}$. By continuity, $P={\cal S}(B)$. But $B$ has connected Julia set while $J(P)$ is disconnected. This is impossible by \thmref{hitD}.
\end{pf}
\begin{cor}
\label{c5full}
The connectedness locus $\CC$ has only two complementary
components $\Lambda_{ext}$ and $\Lambda_{int}$.
\end{cor}
\begin{pf}
Let $U$ be a bounded component of $\BBB C^{\ast}\smallsetminus \CC$ which is not $\Lambda_{int}$. Without loss of generality, we assume that $U$ maps into $\Omega_{ext}$ by $\cal S$. Take $A\in U$. By the previous corollary, there exists a $B\in\Lambda_{ext}$ such that ${\cal S}(A)={\cal S}(B)$. By \thmref{inj}, $A=B$ and this is a contradiction.
\end{pf}
\begin{cor}
\label{surj} 
The surgery map ${\cal S}:\BB\rightarrow \PP$ is surjective.
\end{cor}
\begin{pf}
Compactify $\BB$ and $\PP$ by adding points at $0$ and $\infty$ to get topological 2-spheres. $\cal S$ extends to a continuous map between these spheres by \propref{proper}. This map has topological degree $\neq 0$ because it is a homeomorphism $\Lambda_{ext}\iso \Omega_{ext}$ and ${\cal S}^{-1}(\Omega_{ext})=\Lambda_{ext}$. Therefore it has to be surjective. 
\end{pf}

Since the boundary of the Siegel disk of a cubic which comes from the surgery is a quasicircle passing through some critical point, we have proved the following:

\begin{thm}[Bounded type cubic Siegel disks are quasidisks]
\label{main}
Let $P$ be a cubic polynomial which has a fixed Siegel disk $S$ of rotation number $\theta$. Let $\theta$ be of bounded type. Then the boundary of $S$ is a quasicircle which contains one or both critical points of $P$.
\end{thm}

By a recent theorem of Graczyk and Jones \cite{GJ}, we have the following corollary:

\begin{cor}
\label{HD}
Under the assumptions of \thmref{main}, the boundary of the Siegel disk $S$ has Hausdorff dimension greater than $1$.
\end{cor}

Now it is possible to show that despite all the bifurcations taking place near the boundary of the connectedness locus $\MM$ which give rise to discontinuity of the Julia sets, the boundaries of the Siegel disks move continuously.

\begin{thm}[Boundary of Siegel disks move continuously]
\label{move}
The boundary $\partial \Delta_c$ of the Siegel disk of $P_c\in \PP$ centered at $0$ is a continuous function of $c\in \BBB C^{\ast}$ in the Hausdorff topology.
\end{thm}

\begin{pf}
Let us fix some $P\in \PP$. If $P\notin  \partial \MM$, \thmref{unstable} shows that $J(P)$, hence $\partial \Delta_P$, moves holomorphically in a neighborhood of $P$ and continuity at $P$ is obvious. So let us assume that $P\in \partial \MM$ and consider a sequence $P_n\in \PP$ which converges to $P$ as $n\rightarrow \infty$. Since the surgery map is surjective, there exists a sequence $B_n\in \BB$ such that ${\cal S}(B_n)=P_n$. By properness (\propref{proper}), some subsequence which we still denote by $B_n$ converges to some $B\in \BB$, which by continuity maps to $P$. Now consider the representations
$P_n=\varphi_n \circ \tilde{B}_n\circ \varphi_n^{-1}$ as in (\ref{eqn:pmodb}). Then the boundary $\partial \Delta_{P_n}$ is just the image $\varphi_n(\BBB T)$. Since $\{ \varphi_n \}$ is normal by \corref{normal}, some further subsequence, still denoted by $\{ \varphi_n \}$, converges to a quasiconformal homeomorphism $\psi$. The map $Q=\psi\circ \tilde{B} \circ \psi^{-1}\in \PP$ is quasiconformally conjugate to $P$. Since $P$ is rigid by \thmref{qcclass}, $P=Q$. Now, as $n\rightarrow \infty$, $\partial \Delta_{P_n}=\varphi_n(\BBB T)$ converges in the Hausdorff topology to $\psi(\BBB T)=\partial \Delta_Q=\partial \Delta_P$. 
\end{pf} 

\vspace{0.17in}

\section{Siegel Disks with Two Critical Points on Their Boundary}
\label{sec:twocrit}

In this section we characterize those cubics in $\PP$ which have both critical points on the boundary of their Siegel disk. In \thmref{gamma1} we will prove that the set of all such cubics is a Jordan curve $\Gamma$ in $\PP$. The proof of this theorem will use the fact that the quasiconformal conjugacy classes in $\BB$ are path-connected (\thmref{qcpath}). We then show that when there are no queer components, $\Gamma$ is in fact the common boundary of $\Omega_{ext}$ and $\Omega_{int}$ (\thmref{gamma2}).

Consider the set $\Gamma$ which consists of all cubics $P\in \PP$ such that both critical points of $P$ belong to the boundary of the Siegel disk $\Delta_P$. \figref{Gamma2} shows this set in the parameter space $\PP$.

Since the surgery map ${\cal S}:\BB \rightarrow \PP$ is surjective by \corref{surj}, every $P\in \Gamma$ is of the form ${\cal S}(B_{\mu})$ with $B_{\mu}$ having two double critical points on the circle. \corref{switch} shows that $\mu$ must belong the unit circle $\BBB T\subset \BBB C^{\ast}\simeq \BB$. Therefore, we simply have
$$\Gamma={\cal S}(\BBB T).$$
In particular, $\Gamma$ is a closed path in $\PP \simeq {\BBB C}^{\ast}$. Suggested by \figref{Gamma2}, we want to prove that $\Gamma$ is a Jordan curve. This would follow immediately if we could prove that ${\cal S}|_{\BBB T}$ is injective.
However, I have not been able to show this. In fact, I do not know how to prove that Blaschke products on the boundary of the connectedness locus $\CC$ are quasiconformally rigid. So we take a slightly different approach by showing that the fibers of ${\cal S}|_{\BBB T}:\BBB T\rightarrow \Gamma$ are connected.

\begin{lem}
\label{fiber}
Let $A,B\in \BB$ and ${\cal S}(A)={\cal S}(B)=P$. Suppose that $P$ is not capture. Then there exists a path $t\mapsto A_t \in \BB$ of Blaschke products for $0\leq t\leq 1$, with $A_0=A,\ A_1=B$, such that ${\cal S}(A_t)=P$ for all $t$.
\end{lem}

\begin{pf}
Since $P$ is not capture, by \thmref{ABconj} there exists a quasiconformal homeomorphism $\Phi$ which conjugates $A$ to $B$, which is conformal away from the Julia set $J(A)$. By \thmref{qcpath} there exists a path $\{ \Phi_t \} _{0\leq t\leq 1}$ connecting the identity map to $\Phi$ and a corresponding path $\{ A_t=\Phi_t \circ A\circ \Phi_t^{-1} \} _{0\leq t\leq 1}$ of elements of $\BB$ connecting $A$ to $B$. Note that by the definition of $\Phi_t$, these quasiconformal homeomorphisms are all conformal away from $J(A)$.

It remains to show that ${\cal S}(A_t)=P$ for all $0\leq t\leq 1$. Consider the Douady-Earle extension $H:\BBB D\rightarrow \BBB D$ used in the definition of ${\cal S}(A)$ in  Section \ref{sec:surgery}. Recall that $H|_{\BBB T}$ conjugates $A|_{\BBB T}$ to the rigid rotation $t\mapsto t+\theta$ (mod 1). Hence, the quasiconformal homeomorphism $H_t=H\circ \Phi_t^{-1}:\BBB D\rightarrow \BBB D$ will conjugate $A_t|_{\BBB T}$ to the rigid rotation as well. Note that $H_t$ is not in general the Douady-Earle extension of the linearizing homeomorphism $h_t:\BBB T \rightarrow \BBB T$ for $A_t$. Nevertheless, ${\cal S}_{H_t}(A_t)={\cal S}(A_t)$ by \propref{independent}. Consider the modified Blaschke products
$$ \tilde{A}(z)=  \left \{ \begin{array}{ll} 
A(z) &  |z|\geq 1 \\
(H^{-1}\circ R_{\theta}\circ H)(z) & |z|<1
\end{array}
\right.$$
and
$$ \tilde{A_t}(z)=  \left \{ \begin{array}{ll} 
A_t(z) &  |z|\geq 1 \\
(H_t^{-1}\circ R_{\theta}\circ H_t)(z) & |z|<1
\end{array}
\right.$$
Note that $\Phi_t \circ \tilde{A}=\tilde{A_t}\circ \Phi_t$. 

Define the corresponding conformal structures $\sigma$ and $\sigma_t$ as in  Section \ref{sec:surgery}. It is easy to see that
\begin{equation}
\label{eqn:conf} 
\sigma=\Phi_t^{\ast}\sigma_t.
\end{equation}
Here we use that fact that $\Phi_t$ is conformal away from $J(A)$. Consider the normalized solutions $\varphi$ and $\varphi_t$ of the Beltrami equations
$$\varphi^{\ast}\sigma_0=\sigma,\ \ \ \ \varphi_t^{\ast}\sigma_0=\sigma_t.$$
By (\ref{eqn:conf}) and uniqueness, we have 
$$\varphi_t=\varphi \circ \Phi_t^{-1}.$$
Hence, by \propref{independent},
$$\begin{array}{rl}
{\cal S}(A_t) & =\varphi_t\circ \tilde{A_t}\circ \varphi_t^{-1} \\
              & =\varphi \circ \Phi_t^{-1}\circ \tilde{A_t}\circ \Phi_t                         \circ \varphi^{-1}\\
              & =\varphi \circ \tilde{A}\circ \varphi^{-1}\\
              & ={\cal S}(A).
\end{array}$$
This completes the proof of the lemma.
\end{pf}

\begin{cor}
\label{fiber2}
The fibers of ${\cal S}|_{\BBB T}: \BBB T\rightarrow \Gamma$ are connected.
\end{cor}

\begin{pf}
Let $A,B \in \BBB T \subset \BB$ and ${\cal S}(A)={\cal S}(B)$. Apply the previous lemma to $A,B$. Note that $A_t\in \BBB T$ for all $0\leq t\leq 1$, since $A_t$ is quasiconformally conjugate to $A$, hence has two double critical points on the unit circle. 
\end{pf}

\begin{lem}
\label{moore}
Let $\sim$ be an equivalence relation on the unit circle $\BBB T$ such that every equivalence class is closed and connected. Suppose that the whole circle is not an equivalence class. Then the quotient space $\BBB T/\sim$ is also homeomorphic to $\BBB T$.
\end{lem}

\realfig{Gamma2}{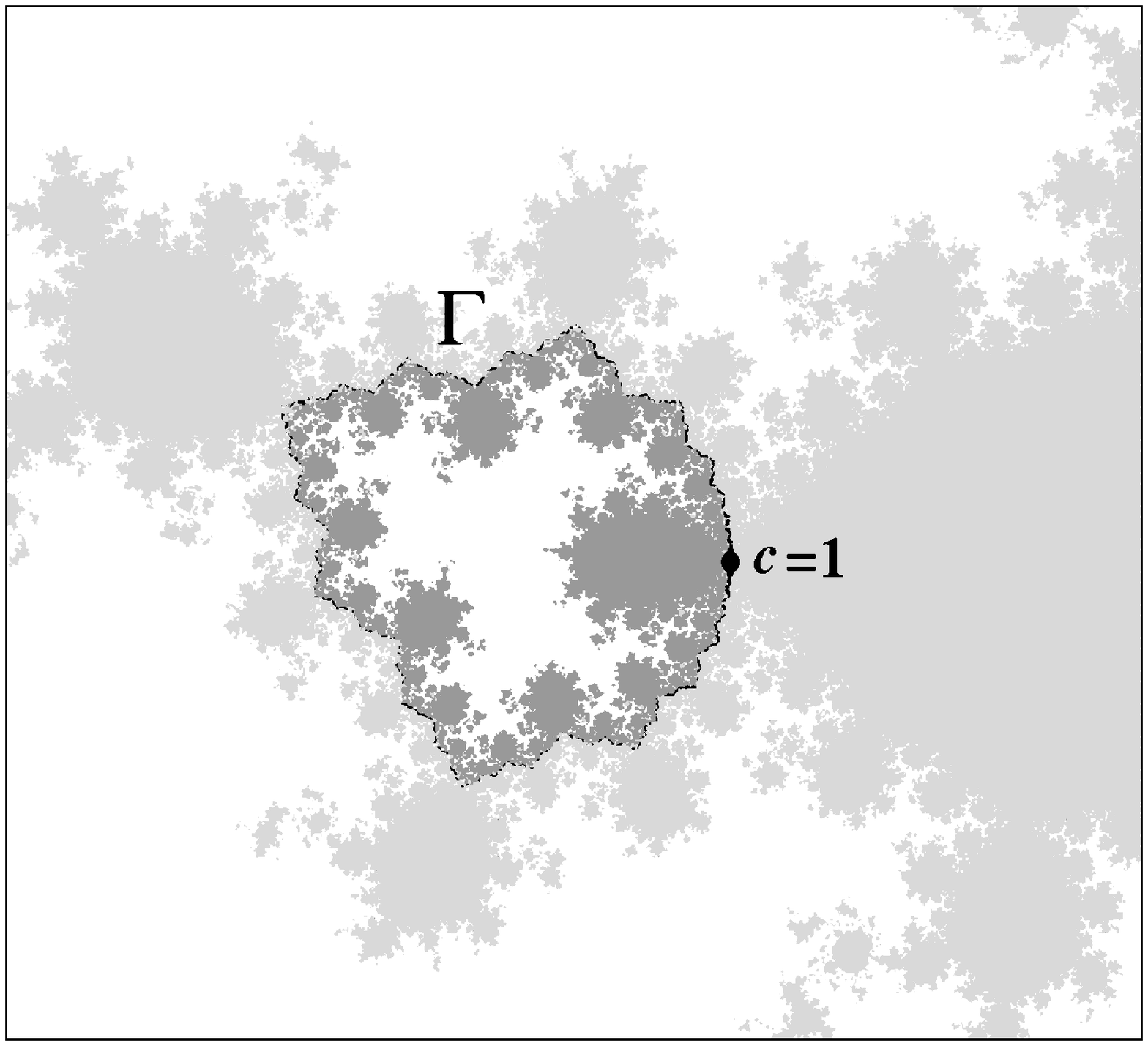}{{\sl The Jordan curve $\Gamma$. This is the locus of all critically marked cubics in $\PP$ which have both critical points on the boundary of their Siegel disk. Topologically it can be described as the common boundary of the complementary regions $\Omega_{ext}$ and $\Omega_{int}$. Note that $\Gamma$ is invariant under the inversion $c\mapsto 1/c$. In particular, it passes through $c=1$.}}{12cm}

\begin{pf}
One can easily construct the homeomorphism as follows: Identify $\BBB T$ with $\BBB R/\BBB Z$, and let $\{ S_i \}_{i\in \BBB N}$ be the collection of nontrivial equivalence classes of $\sim$. (In case this collection is empty or finite, the lemma is clear.) Each $S_i$ can be regarded as a closed interval in $(0,1]$, and we may assume that the right endpoint of $S_1$ is $1$. Note that there is a natural order $<$ on the collection $\{ S_i \}$. Define a function $f:[0,1]\rightarrow [0,1]$ by putting $f(0)=0$, $f|_{S_1}\equiv 1$, and $f|_{S_2}\equiv 1/2$ and proceed inductively as follows. Suppose that $n\geq 2$ and $f$ is already defined on $S_1 \cup \cdots \cup S_n$. Consider $S_{n+1}$ and let $S_i$ and $S_j$ be its two neighbors with $S_i < S_{n+1} < S_j$ and $1\leq i,j \leq n$. 
Define $f|_{S_{n+1}}\equiv 1/2(f(S_i)+f(S_j))$. (In case $S_{n+1}$ has no left neighbor, simply set $f|_{S_{n+1}}\equiv 1/2f(S_j)$.) This defines $f$ inductively on $\bigcup S_i$. By the construction, $f$ extends continuously to the closure $\overline{\bigcup S_i}$. Interpolate $f$ linearly on each open interval in $(0,1)\smallsetminus \overline{\bigcup S_i}$. 

It is easy to check that $f$ constructed this way is continuous, increasing, and the preimage of every point is either a single point in $[0,1]\smallsetminus \bigcup S_i$ or an interval $S_i$. Clearly such a function induces a homeomorphism between $\BBB T$ and $\BBB T/\sim$.
\end{pf}

\noindent
{\bf Remark.} This simple lemma should be thought of as the one-dimensional (baby) version of the following deep theorem of R. L. Moore \cite{Moore}: Let $\sim$ be a closed equivalence relation on the 2-sphere $S^2$ such that every equivalence class is closed and connected and nonseparating. Then $S^2/\sim $ is also homeomorphic to $S^2$.

\begin{thm}
\label{gamma1}
$\Gamma$ is a Jordan curve.
\end{thm}

\begin{pf}
Consider ${\cal S}|_{\BBB T}:\BBB T\rightarrow \Gamma$ whose fibers are closed and connected by \lemref{fiber}. By general topology, $\Gamma$ is homeomorphic to $\BBB T/\sim$, where $A\sim B$ means ${\cal S}(A)={\cal S}(B)$. By \lemref{moore}, $\BBB T/\sim$ is homeomorphic to the circle.
\end{pf}

Finally, we find a topological characterization of $\Gamma$ in $\PP$ under the assumption that there are no queer components in the interior of $\MM$.

\begin{thm}[Topological characterization of $\Gamma$]
\label{gamma2}
$\Gamma$ is a subset of the boundary $\partial \MM$ which contains $\partial \Omega_{ext} \cap \partial \Omega_{int}$. If there are no queer components in the interior of $\MM$, then $\Gamma=\partial \Omega_{ext} \cap \partial \Omega_{int}$.
\end{thm}

\begin{pf}
First let us show that $\partial \Omega_{ext} \cap \partial \Omega_{int}\subset \Gamma $. Let $P\in \partial \Omega_{ext} \cap \partial \Omega_{int}$ and assume that $P\notin \Gamma$. Choose $B_{\mu}\in \BB$ such that ${\cal S}(B_{\mu})=P$. We can assume without loss of generality that $|\mu|>1$. Choose a sequence $P_n\in \Omega_{int}$ converging to $P$ and a sequence $B_n\in \Lambda_{int}$ such that ${\cal S}(B_n)=P_n$. By passing to a subsequence we may assume that $B_n\rightarrow B_{\mu '}$ as $n\rightarrow \infty$. By continuity, ${\cal S}(B_{\mu '} )=P$ so we must have $|\mu '|<1$. Since $P$ is not capture by \corref{capop}, \lemref{fiber} shows that there is a path $t\mapsto B_t$ of quasiconformally conjugate Blaschke products connecting $B_{\mu}$ to $B_{\mu '}$ all of which are mapped to $P$. Since this path must intersect $\BBB T$ somewhere, we conclude that $P\in \Gamma$ which is a contradiction. 
   
Now we prove that $\Gamma \subset \partial \MM$. Fix some $P\in \Gamma$. Since $P$ has both critical points on $\partial \Delta_P$, it cannot belong to any hyperbolic-like or capture component. Also, $P$ cannot be in a queer component $U$ of the interior of $\MM$, since otherwise every $Q\in U$ would have to be quasiconformally conjugate to $P$ by \thmref{qcclass}, which would imply that $Q$ has two critical points on $\partial \Delta_Q$, which would show $U\subset \Gamma$. But this is evidently impossible because $U$ is open and $\Gamma$ is a Jordan curve. Therefore, $P$ has to lie in $\partial \MM=\partial \Omega_{ext} \cup \partial \Omega_{int}$. 
\realfig{DD}{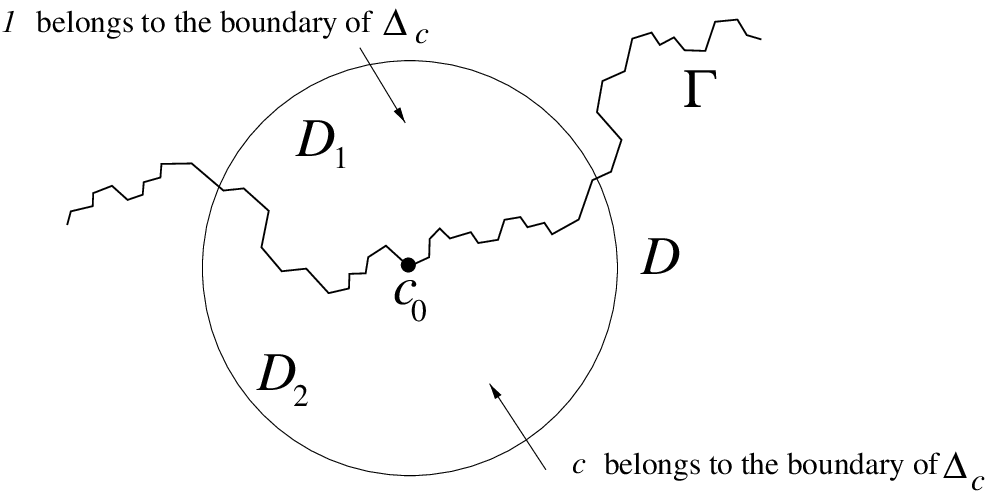}{}{8cm}

Now assume that there are no queer components in the interior of $\MM$. To show that $\Gamma=\partial \Omega_{ext} \cap \partial \Omega_{int}$, let $P=P_{c_0}$ and assume by way of contradiction that $c_0\in \partial \Omega_{ext} \smallsetminus 
\partial \Omega_{int}.$ Since $c_0$ has positive distance from $\Omega_{int}$, for all $c$ in a neighborhood $D$ of $c_0$ the sequence $\{ P_c^{\circ n}(1) \}$ has to be normal. Assuming that $D$ is a small disk, the Jordan curve $\Gamma$ cuts $D$ into two topological disks $D_1$ and $D_2$ such that for every $c\in D_1$, $1 \in \partial \Delta_c$ and $c\notin \partial \Delta_c$, and for every $c\in D_2$, $c\in \partial \Delta_c$ and $1 \notin \partial \Delta_c$ (see \figref{DD}). 

Clearly $D_2\cap \partial \Omega_{ext}=D_2 \cap \partial \Omega_{int}=\emptyset$. So $D_2$ has to be a subset of a component $U$ of the interior of $\MM$. Since there are no queer components by the assumption, $U$ is either hyperbolic-like or capture.

For every $c\in D_1$, we have $1 \in \partial \Delta_c$ and the restriction $P_c|_{\partial \Delta_c}$ is conjugate to the rigid rotation by angle $\theta$. Therefore, $P_c^{\circ q_n}(1)\rightarrow 1$ for all $c\in D_1$, where the $q_n$ are the denominators of the rational approximations of $\theta$. Since $\{ P_c^{\circ n}(1) \} $ is normal in $D$, for a subsequence $\{ q_{n(j)} \} $ we must have $P_c^{\circ q_{n(j)}}(1)\rightarrow 1$ throughout $D$. In particular, if $c\in D_2$, the critical point $1$ of $P_c$ must be recurrent. This is impossible if $U$ is hyperbolic-like or capture, since over $D_2$, $c\in \partial \Delta_c$ and hence $1$ either gets attracted to the attracting cycle or eventually maps to the Siegel disk $\Delta_c$.   
\end{pf}

\end{document}